\documentclass[12pt]{amsart}
\usepackage{amsmath,amssymb,array,longtable, MnSymbol, yfonts}
\usepackage[mathscr]{euscript}
\usepackage[colorlinks=true,urlcolor=black,citecolor=black,linkcolor=black,%
pdftitle={Seshadri constants, Diophantine approximation, and Roth's theorem for arbitrary varieties},%
pdfauthor={David McKinnon and Mike Roth},%
pdfsubject={Local version of Bombieri-Lang phenomena},%
pdfkeywords={Algebraic geometry, Diophantine approximation, Seshadri constants}]{hyperref}
\usepackage{pstricks}    
\usepackage{pst-plot}

\newtheorem{theorem}{Theorem}[section]

\newtheorem{proposition}[theorem]{Proposition}
\newtheorem{lemma}[theorem]{Lemma}
\newtheorem{corollary}[theorem]{Corollary}
\newtheorem{definition}[theorem]{Definition}

\newcounter{figcount}[section]
\renewcommand{\thefigcount}{\thesection\alph{figcount}}
\newcounter{eqcount}[section]
\renewcommand{\theeqcount}{\thesection.\arabic{eqcount}}

\newcommand{\Qbar}{\overline{k}}        
\newcommand{\Fbar}{\overline{F}}        
\newcommand{\mult}{\operatorname{mult}} 
\newcommand{\Xtil}{\widetilde{X}}       
\newcommand{\Xbar}{\overline{X}}        
\newcommand{\Ytil}{\widetilde{Y}}       
\newcommand{\Ybar}{\overline{Y}}        
\newcommand{\Ztil}{\widetilde{Z}}       
\newcommand{\Ptil}{\widetilde{\PP}}     
\newcommand{\Pic}{\operatorname{Pic}}   
\newcommand{\sepv}{d_v}                 
\newcommand{\sepvp}{d_{v'}}             
\newcommand{\sepw}{d_w}                 
\newcommand{\squareplus}{\boxplus}      
\newcommand{\imX}{\textfrak{X}}         
\newcommand{\imXt}{\widetilde{\imX}}    
\newcommand{\imXv}{\textfrak{X}_{v}}    
\newcommand{\Okv}{\widehat{\Osh}_{k,v}}   
\newcommand{\OFv}{\widehat{\Osh}_{F,v}}   
\newcommand{\Spec}{\operatorname{Spec}} 
\newcommand{\Osh}{\mathcal{O}}          
\newcommand{\Esh}{\mathcal{E}}          
\renewcommand{\AA}{\mathbb{A}}          
\newcommand{\CC}{\mathbb{C}}            
\newcommand{\PP}{\mathbb{P}}            
\newcommand{\QQ}{\mathbb{Q}}            
\newcommand{\RR}{\mathbb{R}}            
\newcommand{\bpf}{\noindent{\em Proof:} \/}  
\newcommand{\epf}{$\Box$}               
\newcommand{\xseq}{\{x_i\}}             
\newcommand{\xseqtox}{\{x_i\}\rightarrow x}  
\newcommand{\yseq}{\{y_i\}}             
\newcommand{\zseq}{\{z_i\}}             
\newcommand{\qseq}{\{q_i\}}             
\newcommand{\xsubseq}{\{x_i'\}}         
\newcommand{\yseqn}{\{y_{i,n}\}}        
\newcommand{\ynseq}{\{y_{n}\}}          
\newcommand{\st}{\,\,\rule[-0.2cm]{0.02cm}{0.6cm}\,\,}   
\newcommand{\Kv}{K^{(v)}}     
\newcommand{\Xv}{X^{(v)}}               
\newcommand{\Xvtil}{\widetilde{X}^{(v)}}
\newcommand{\Lvz}{L^{(v)}_{0}}          
\newcommand{\Lvg}{L^{(v)}_{\gamma}}     
\newcommand{\Lvgj}{L^{(v)}_{\gamma_{v,j}}} 
\newcommand{\Lvgjpo}{L^{(v)}_{\gamma_{v,j+1}}} 
\newcommand{\Ev}{E^{(v)}}                  
\newcommand{\ep}{\epsilon}                 
\newcommand{\ephat}{\hat{\ep}}             
\newcommand{\ephatet}{\hat{\ep}^{\et}}     
\newcommand{\etale}{\'etale }              
\newcommand{\et}{\,\mbox{\scriptsize\'et}}  
\newcommand{\hhat}{\hat{h}}                
\newcommand{\Ar}{\beta}                    
\newcommand{\Arhat}{\hat{\Ar}}             
\newcommand{\Arhatet}{\hat{\Ar}^{\et}}     
\newcommand{\geff}{\gamma_{\mbox{\tiny eff}}}     
\newcommand{\geffx}{\gamma_{\mbox{\rm\tiny eff},x}}  
\newcommand{\geffxv}{\gamma_{\mbox{\tiny eff},x_v}}
\newcommand{\Vol}{\operatorname{Vol}}    
\newcommand{\NS}{\operatorname{NS}}      
\newcommand{\Whn}{W_{\!\circ}}        
\newcommand{\muo}{\mu_{\circ}}        
\newcommand{\Gal}{\operatorname{Gal}} 
\renewcommand{\geq}{\geqslant}        
\renewcommand{\leq}{\leqslant}        
\newcommand{\Cnrm}[1]{|#1|}         
\newcommand{\nrm}[1]{||#1||}        
\newcommand{\avec}{\mathbf{a}}    
\newcommand{\bvec}{\mathbf{b}}    
\newcommand{\IshE}{\mathscr{I}_{E}}      
\newcommand{\np}{\medskip\noindent}      

\providecommand{\binom}[2]{{#1\choose#2}}

\renewenvironment{equation}{\medskip\noindent\refstepcounter{eqcount}\makebox[0pt][l]{\rm ({\bf\theeqcount})}\begin{minipage}[b]{\textwidth}$$}{$$\end{minipage}\medskip\noindent}

\newcommand{\fig}{\refstepcounter{figcount}Figure \thefigcount.}

\newcounter{dummycounter}
\newcommand{\forcelabel}[2]{
\renewcommand{\thedummycounter}{#2}
\refstepcounter{dummycounter}
\label{#1}
}

\newcommand{\forceandshowlabel}[2]{\forcelabel{#1}{#2}\thedummycounter}

\newcommand{\comment}[1]{}

\makeatletter
\newcommand{\xdashrightarrow}[2][]{\ext@arrow 0359\rightarrowfill@@{#1}{#2}}
\def\rightarrowfill@@{\arrowfill@@\relax\relbar\rightarrow}
\def\arrowfill@@#1#2#3#4{%
  $\m@th\thickmuskip0mu\medmuskip\thickmuskip\thinmuskip\thickmuskip
   \relax#4#1
   \xleaders\hbox{$#4#2$}\hfill
    #3$%
}
\makeatother
\newcommand{\longdashrightarrow}{\xdashrightarrow{\rule{0.49cm}{0cm}}}

\begin{document}
\title[Seshadri constants and Roth's theorem for arbitrary varieties]{Seshadri constants, diophantine approximation, and Roth's theorem for arbitrary varieties}

\author{David McKinnon}
\address{Department of Pure Mathematics \\
University of Waterloo \\
Waterloo, ON\ \  N2L 3G1 \\
Canada}
\email{dmckinnon@math.uwaterloo.ca}
\indent
\thanks{David McKinnon was partially supported by an NSERC research grant.}

\author{Mike Roth}
\address{Dept.\ of Mathematics and Statistics, Queens University, Kingston,
Ontario, Canada}
\email{mikeroth@mast.queensu.ca}
\indent
\thanks{Mike Roth was partially supported by an NSERC research grant.}


\begin{abstract}
In this paper, we associate an invariant $\alpha_{x}(L)$ to an algebraic point $x$ on an algebraic variety $X$
with an ample line bundle $L$.  The invariant $\alpha$ measures how well $x$ can be approximated by rational
points on $X$, with respect to the height function associated to $L$.
We show that this invariant is closely related to the Seshadri constant $\ep_{x}(L)$ measuring local positivity
of $L$ at $x$, and in particular that Roth's theorem on $\PP^1$ generalizes as an inequality between these
two invariants valid for arbitrary projective varieties.
\end{abstract}


\maketitle

\section{Introduction}

Let $k$ be a number field, and $X$ an irreducible projective variety over $\Spec(k)$.
The Bombieri-Lang conjecture predicts
that if $X$ is of general type then the $k$-points of $X$ are contained in a proper closed subset of $X$.
We view this as a statement that a global fact about the canonical bundle of $X$ (that it is ``generically positive'',
where positivity is used in a broad sense) implies a global fact about the accumulation of rational points.
Following a well-established principle in geometry
one should study the local influence of positivity on the local accumulation of rational points.  To do this
we need local measures of both these phenomena.

Let $L$ be an ample line bundle on $X$, and $x$ a point of $X(\Qbar)$.
By slightly modifying the usual definition of approximation exponent on $\PP^1$ (and inspired by a definition
from \cite{McK} by the first author) we define a new invariant $\alpha_{x}(L)\in (0,\infty]$ which measures how quickly
rational points accumulate around $x$, from the point of view of the line bundle $L$ and a fixed place $v$ of $k$.

The central theme of this paper is the interrelations between $\alpha_x(L)$ and the Seshadri constant $\ep_x(L)$,
an invariant defined by Demailly \cite{Dem} which measures local positivity of a line bundle $L$ near a point $x$.
The two share common formal properties, and this similarity is even more evident when $\alpha_x$ is interpreted
through Arakelov theory.
Moreover, the classic approximation results on $\PP^1$ --- the theorems of
Liouville and Roth --- generalize as inequalities between $\alpha_x$ and $\ep_x$ valid for arbitrary projective
varieties.  This general version of Roth's theorem admits further generalizations to simultaneous approximation and
improvements via \etale covers.

In order to motivate our results we first quickly review approximation on the line, and
to simplify this part of the discussion we assume that $k=\QQ$ and that the place is archimedean.

\medskip
\noindent
{\bf Approximation on $\AA^1$.}  For a point $x\in \RR$ the {\em approximation exponent} of $x$
is defined as the smallest real number $\tau_{x}$ such that for any $\delta>0$
the inequality
$$\left|{x-\frac{a}{b}}\right| \leq \frac{1}{b^{\tau_{x}+\delta}}$$
has only finitely many solutions $a/b\in \QQ$ (see \cite[Part D]{HS}).  The approximation exponent
measures a certain tension between our ability to closely approximate $x$ by rational numbers (the $|x-a/b|$ term)
and the complexity (the $1/b$ term) of the rational number needed to make this approximation.

If $x\in \QQ$ then it is easy to see that $\tau_{x}=1$.  In 1842 Dirichlet proved his famous approximation theorem:
if $x\in \RR\setminus\QQ$ then $\tau_{x}\geq 2$. One therefore seeks upper bounds on $\tau$.
In 1844 Liouville showed that
if $x\in \RR$ is algebraic of degree $d$ over $\QQ$ then $\tau_{x}\leq d$,
and used this to give concrete examples of transcendental numbers.  Further improvements
in the upper bound were obtained by Thue (1909), Siegel (1921), and Dyson and Gelfand (1947), culminating in the
1955 theorem of Roth: for $x\in \RR$ algebraic over $\QQ$, $\tau_{x}\leq 2$.  Thus the theorems of Dirichlet
and Roth give $\tau_{x}=2$ for irrational algebraic $x\in \RR$.

\noindent
{\bf The invariant $\alpha_{x}(L)$.}
In \S\ref{sec:rat-approx} we generalize the approximation exponent to arbitrary projective varieties
$X$ defined over a number field $k$.  To do this we replace the function $|x-a/b|$ by a distance function
$\sepv(x,\cdot)$ depending on a place $v$ of $k$, and measure the complexity of a rational point via
a height function $H_{L}(\cdot)$ depending on an ample line bundle $L$.
The one essential change in our definition is to move the exponent from the height to the distance.
As a result, as Proposition \ref{prop:equiv-alpha} shows,
for $x\in \RR = \AA^1(\RR)\subset \PP^1(\RR)$ we have $\alpha_{x}(\Osh_{\PP^1}(1))=\frac{1}{\tau_x}$.
The choice of moving the exponent is justified by Proposition \ref{prop:alpha}(a,b) which shows that this form
is more natural when we vary $L$, and by the resulting similarities with the Seshadri constant.

In particular, for $x\in \RR\setminus\QQ$, algebraic of degree $d$ over $\QQ$, and $L=\Osh_{\PP^1}(1)$
the theorems of Liouville and Roth become the lower bounds $\alpha_{x}(L)\geq \frac{1}{d}$ and
$\alpha_{x}(L)\geq \frac{1}{2}$ respectively.  One of the main goals of this paper is to generalize these
statements to lower bounds for $\alpha_{x}(L)$ on an arbitrary variety $X$.

\noindent
{\bf Examples.} \label{ex:alpha-intro}
Here are three examples of lower bounds on $\alpha$ given by previously known results on Diophantine approximation.
We work over an arbitrary number field $k$.

\begin{enumerate}
\item If $X=\PP^1_k$, $x\in X(\Qbar)$, $L=\Osh_{\PP^1}(1)$, then $\alpha_{x}(L)\geq \frac{1}{2}$.
\item If $X=\PP^n_k$, $x\in X(\Qbar)$, $L=\Osh_{\PP^n}(1)$, then either $\alpha_{x}(L)\geq \frac{n}{n+1}$
or there is a smaller linear space $Z\cong \PP^m_k \subset \PP^n_k$,  with $m<n$, and $x\in Z(\Qbar)$.
\item If $X$ is an abelian variety, $x\in X(\Qbar)$, and $L$ any ample line bundle then $\alpha_{x}(L)=\infty$.
\end{enumerate}

Example (a) is Roth's theorem for a general number field (and place $v$), and example (b) follows from the Schmidt
subspace theorem.  In both of these cases by using a Dirichlet-type argument \cite{CV} 
one obtains exact values for $\alpha_x$.
In the case of $\PP^1$, if $x\in \PP^1(\Qbar\cap k_v)\setminus \PP^1(k)$ then $\alpha_{x}(\Osh_{\PP^1}(1))=\frac{1}{2}$.
In the case of $\PP^n$, if $x\in \PP^n(\Qbar\cap k_v)\setminus \PP^n(k)$, and
$m$ is the smallest value so that there exists a linear subspace $\PP^m_k\subset \PP^n_k$  with $x\in \PP^m(\Qbar)$,
then $\alpha_{x}(\Osh_{\PP^n}(1))=\frac{m}{m+1}$.
Finally example (c) is \cite[p. 98; second theorem]{Se}.

The basic interpretation of $\alpha_{x}(L)$ is as the cost in complexity
required to get closer to $x$.  When $\alpha_{x}$ is finite this indicates that the complexity
has polynomial growth in the reciprocal of the distance, with $\alpha_{x}$ as the exponent.
In example (c) the complexity grows roughly exponentially in the reciprocal of the distance
(see \cite[p.\ 98 again]{Se}) and thus $\alpha_{x}=\infty$.

\noindent
{\bf The invariant $\mathbf \ep_{x}(L)$}.
The definition and elementary properties of the Seshadri constants are given in \S\ref{sec:seshadri}. We
list two of these properties, and the corresponding properties for $\alpha$, here in order to emphasise the similarity
between $\alpha_x$ and $\ep_x$, and to use one of the properties in the discussion below.
Both $\alpha$ and $\ep$ make sense for $\QQ$-bundles. 
Fix $x\in X(\Qbar)$, then

\begin{enumerate}
\item for any ample $\QQ$-bundle $L$, and any $m\in \QQ_{>0}$,
$$\alpha_{x}(mL) = m\alpha_{x}(L)\,\,\,\mbox{and}\,\,\, \ep_{x}(mL) = m\ep_{x}(L);$$

\item $\alpha_x$ and $\ep_x$ are concave functions of the line bundle.  For any ample $\QQ$-bundles $L_1$
and $L_2$, and any $a,b\in \QQ_{\geq_0}$,
$$ \alpha_{x}(aL_1+bL_2) \geq a\,\alpha_{x}(L_1) + b\,\alpha_{x}(L_2)\,\,\,\mbox{and}\,\,\,
\ep_{x}(aL_1+bL_2) \geq a\,\ep_{x}(L_1) + b\,\ep_{x}(L_2).$$
\end{enumerate}

\noindent
These and other parallel properties appear in Propositions \ref{prop:alpha} and \ref{prop:ex}.

\noindent
{\bf Examples.}

\begin{enumerate}
\item If $X=\PP^n$, $x\in X(\CC)$, and $L=\Osh_{\PP^n}(1)$ then $\ep_{x}(L)=1$, and so $\ep_{x}(\Osh_{\PP^n}(e))=e$
for all $e> 0$.
\item If $X$ is a smooth cubic surface, and $L=\Osh_{\PP^3}(1)|_{X}$ then
\[
\ep_{x}(L) =
\begin{cases}
1 & \text{if $x$ is on a line} \\
\frac{3}{2} & \text{otherwise.}
\end{cases}
\]
\end{enumerate}

If $X$ is a variety with a transitive group action, such as $\PP^n$ or an abelian variety, then the
value of $\ep_{x}(L)$ is independent of $x\in X(\CC)$.
One thesis of this paper is that $\ep_{x}$ affects approximation results.  On varieties where
$\ep_{x}$ does not depend on the point this effect is essentially invisible since it becomes a global
property of the line bundle.  On arbitrary varieties however one can expect more precise
approximation theorems by taking the differing values of $\ep$ into account.  This will be a feature
of the results below.

\noindent
{\bf Roth theorems.}\footnote{All uses of ``Roth'' as an adjective in this paper are in homage to
the theorem proved by Klaus F.\ Roth and its later extensions by Ridout and Lang, and do not refer to the second named author of the paper.} 
If $X$ is a variety over $\Spec(k)$, and $x\in X(\Qbar)$ with field of definition $K$,
then for any ample line
bundle $L$ on $X$ we have $\alpha_{x}(L)\geq \frac{1}{d} \ep_{x}(L)$, where $d=[K\colon k]$.  On $\PP^1$, this
is the inequality $\alpha_{x}(\Osh_{\PP^1}(1))\geq \frac{1}{d}$, and hence we regard this as the general version
of Liouville's theorem.   This result follows from elementary properties of the height of the exceptional divisor
(see the end of \S\ref{sec:seshadri} or \cite[\S3]{McKR} for a proof).

Our main concern is proving general ``Roth'' theorems.  By this we mean lower bounds on $\alpha_{x}(L)$ that are:
(1) independent of the field of definition of $x$, and (2) (following the philosophy of this paper)
expressed in terms of $\ep_{x}(L)$.
The examples of $\PP^n$ and $\PP^1$ suggest two possible interpretations of this goal.

First, based on the example of $\PP^n$ one might hope for a theorem of the form:
for every $n\geq 1$ there is a constant $c_n$ so that for
every irreducible $n$-dimensional variety $X$, ample line bundle $L$ and $x\in X(\Qbar)$,
either $\alpha_{x}(L)\geq c_n \ep_{x}(L)$ or there is a proper subvariety $Z$, with $x\in Z(\Qbar)$,
such that $\alpha_{x}(L)=\alpha_{x,Z}(L|_{Z})$.

Second, one might seek to generalize the $\PP^1$ example:  there is a constant $c$ so that for every variety
$X$, every ample line bundle $L$, and every $x\in X(\Qbar)$ the inequality $\alpha_{x}(L)\geq c\,\ep_{x}(L)$ holds.
Considering varieties of the form $X=\PP^1 \times Y$ shows that $c=\frac{1}{2}$ is the best possible constant
(i.e., it does not help to have the constant vary with the dimension of $X$).

We establish versions of both of these statements;  here is our version of the first type. 

\noindent
{\bf Theorem (\ref{thm:RothII}, ``Schmidt type''):} {\em  Let $X$ be an irreducible $n$-dimensional
variety over $\Spec(k)$.  For any ample $\QQ$-bundle $L$ and any $x\in X(\Qbar)$ either

\begin{enumerate}
\item $\alpha_x(L) \geq \frac{n}{n+1}\ep_x(L)$

\hfill or \hfill\rule{0.1cm}{0.0cm}

\medskip
\item there exists a proper subvariety $Z\subset X$, irreducible over $\Qbar$,  with $x\in Z(\Qbar)$
so that $\alpha_{x,X}(L) = \alpha_{x,Z}(L|_{Z})$, i.e., ``$\alpha_x(L)$ is computed on a proper subvariety of $X$''.
\end{enumerate}
}

\medskip
\noindent
This theorem has an equivalent version expressed in more familiar terms.

\noindent
{\bf Theorem (\ref{thm:RothII}, alternate statement):} {\em
Let $L$ be any ample $\QQ$-bundle on $X$, and choose any $x\in X(\Qbar)$.
Then there is a proper subvariety $Z\subset X$ so that for each $\delta>0$
there are only finitely many solutions $y\in X(k)\setminus Z(k)$ to
$$\sepv(x,y) < H_{L}(y)^{-(\frac{n+1}{n\,\ep_{x}(L)}+\delta)}.$$
}

Theorem \ref{thm:RothII} generalizes the Schmidt subspace theorem, insofar as the Schmidt theorem
concerns approximating a point.
It is an important part of the Schmidt theorem that $Z$ be a union of linear spaces
so that the theorem may be applied inductively.   Since Theorem \ref{thm:RothII}
applies to arbitrary varieties, the ability to apply induction of this type is automatic.
In particular, since the Seshadri constant is weakly increasing when restricting to a subvariety
(Proposition \ref{prop:ex}(c)), Theorem \ref{thm:RothIII} and  induction on dimension yield a
theorem of the second type.

\medskip
\noindent
{\bf Theorem (\ref{thm:RothIII}, ``Roth type''):} {\em For all varieties $X$ over $\Spec(k)$
(possibly reducible), all $x\in X(\Qbar)$ and all ample line bundles $L$, $\alpha_{x}(L)\geq \frac{1}{2}\ep_{x}(L)$. }

In order for equality to hold in Theorem \ref{thm:RothIII} the induction must have gone down to a one-dimensional
variety, and from this we deduce that if equality holds then there is a $k$-rational curve $C$ passing through $x$,
and unibranch at $x$, which also computes the Seshadri constant, i.e., $\ep_{x}(L) = \ep_{x,C}(L|_{C})$.
The exact statement and its converse appear as part of Theorem \ref{thm:RothIII}, as fully stated in
\S\ref{sec:Roth-theorems}.  This is one of the few examples we know of where an arithmetic condition about
approximation implies a geometric condition about $X$ (namely that there must be a rational curve passing through $x$).
If there is no rational curve passing through $x$ then the lower bound in Theorem \ref{thm:RothIII} may be improved;
see Corollary  \ref{cor:no-rat-curve}.

It is useful to state Theorem \ref{thm:RothIII} in an equivalent form closer to that of the usual
statement of Roth's theorem on $\PP^1$.

\noindent
{\bf Corollary (\ref{cor:one-place-Roth}):} {\em For any $\delta>0$ there are only finitely many $y\in X(k)$
such that
$$\sepv(x,y) < H_{L}(y)^{-\left({\frac{2}{\ep_{x}(L)}+\delta}\right)}.$$
}

\noindent
{\bf Heuristic explanation.}
Given an ample line bundle $L$ and $x\in X(\Qbar)$ consider the problem of finding an exponent $e$ so that
for all $\delta>0$ there are only finitely many solutions $y\in X(k)$ to $\sepv(x,y) < H_{L}(y)^{-(e+\delta)}.$
If $m$ is such that $mL$ is very ample, then embedding $X$ via $mL$, projecting on coordinates, and using
Roth's theorem for $\PP^1$ shows that the exponent $e=2m$ will do.  The smaller the value of $e$, the stronger
such a statement is, so we now ask the question: what is the smallest value of $m$ so that $mL$ is very ample?

If $A$ is a very ample line bundle, then $\ep_{x'}(A)\geq 1$ for all $x'\in X(\Qbar)$ 
(see Proposition \ref{prop:ex}(d)).
In particular, if $mL$ is very ample then we must have $m\ep_{x}(L) = \ep_{x}(mL) \geq 1$, and thus
$m\geq \frac{1}{\ep_{x}(L)}$.
In general $m=\frac{1}{\ep_{x}(L)}$ does not guarantee that $mL$ is very ample.
There are basically three problems.  (1) We need $\alpha_{x'}(mL)\geq 1$ for all $x'\in X$, and not just $x$.
(2) Even if the previous condition holds, this does not guarantee that $mL$ is very ample.  (3) With this
value of $m$, $mL$ may not be an integral (or, conjecturally, even a rational) line bundle.

As an example of two of these issues, let $X$ be a smooth cubic surface, $L=\Osh_{\PP^3}(1)|_{X}$,
and $x\in X(\Qbar)$ a point not on a line.
As stated above $\ep_{x}(L)=\frac{3}{2}$.  However $\frac{2}{3}L$ is not an integral line bundle
(it has degree $\frac{2}{3}$ on every line), nor is $\ep_{x'}(\frac{2}{3}L)\geq 1$ for points $x'$ on a line.

The essential point of Corollary \ref{cor:one-place-Roth} is that these concerns don't matter:  as long as we only care
about approximating $x$ the local estimate of amplitude $m=\frac{1}{\ep_{x}(L)}$ works.
This is a good illustration of the effects of local positivity on approximation.

\medskip
\noindent
{\bf Simultaneous approximation.}
As with Roth's theorem on $\PP^1$, our theorems admit generalizations to simultaneous approximation.
In order to indicate the nature of the results let us consider the two equivalent statements for a single place
given by Theorem \ref{thm:RothIII} and Corollary \ref{cor:one-place-Roth} above and see how they generalize.
In \S\ref{sec:rat-approx}, as part of defining $\alpha_{x}(L)$ we also define $\alpha_{x}(\xseq,L)$ for
any sequence $\xseq$ of $k$-points of $X$, and we will need this
notation to state our results below.  In particular, Theorem \ref{thm:RothIII} can be equivalently
stated as $\alpha_{x}(\xseq,L)\geq \frac{1}{2}\ep_{x}(L)$ for all
sequences $\xseq$ of $k$-points of $X$.

To set up the simultaneous approximation problem let $S$ be a finite set of places of $k$, each extended to $\Qbar$.
For each $v\in S$ let $\sepv(\cdot,\cdot)$ be the distance function computed with respect to $v\in S$ and
choose a point $x_v\in X(\Qbar)$.
To simplify notation, we set $\alpha_v$ to be $\alpha_{x_v}$ computed with respect to $\sepv$.

We are interested in understanding how well sequences of $k$-points can simultaneously approximate each $x_v$.
The generalizations of Theorem \ref{thm:RothIII} and Corollary \ref{cor:one-place-Roth} to simultaneous
approximation (see Corollary \ref{cor:alpha-hypersurface}) are respectively:

\begin{itemize}
\item[(1)] for any sequence $\xseq$ of $k$-points, $\sum_{v\in S} \frac{\ep_{x_v}(L)}{\alpha_{v}(\xseq,L)} \leq 2$, and

\medskip
\item[(2)] for any $\delta>0$ there are only finitely many $y\in X(k)$ such that
$$\prod_{v\in S} \sepv(x_v,y)^{\ep_{x_v}(L)} < H_{L}(y)^{-(2+\delta)}.$$
\end{itemize}

The other results (e.g., Theorems \ref{thm:RothI} and \ref{thm:RothII}) also have their simultaneous versions.
Full statements and further discussion appear in \S\ref{sec:simul-approx}.

\medskip
\noindent
{\bf Improvements via \etale covers.}
Given $X$ (which we assume normal to simplify the discussion),
an ample line bundle $L$ on $X$, and $x\in X(\Qbar)$ we define $\ephatet_{x}(L)$ by
$$\ephatet_{x}(L) = \sup_{y\in \varphi^{-1}(x)} \ep_{y}(\varphi^{*}L)$$
where the supremum is over all irreducible \etale covers $\varphi\colon Y\longrightarrow X$.
In \S\ref{sec:unram-bounds} we show that all the previous theorems, for one place or simultaneous places,
hold with $\ep_{x}(L)$ replaced by $\ephatet_{x}(L)$ (see Corollary \ref{cor:everything}).
Since $\ephatet$ is in general larger, this can be a significant strengthening  of the results.  For instance,
if $X$ is an abelian variety and $L$ an ample line bundle, then $\ep_{x}(L)$ is always finite, while
$\ephatet_{x}(L)=\infty$ (see the example on page \pageref{ex:ab-var-etale}).
Thus Theorem \ref{thm:RothIII} applied with $\ephatet$ in place of $\ep$ shows that $\alpha_{x}(L)=\infty$
on an abelian variety.

The results proved in \S\ref{sec:unram-bounds} are slightly more general (for instance, one can take the supremum
over irreducible unramified covers) and the reader is referred there for more detailed statements.

In his 1962 book {\em Diophantine Geometry}, Lang (\cite[p. 119]{La}) suggests three directions for future
progress on Roth's theorem.   The first is to make the result quantitative, and we seem to know as much about
this now as was known in 1962; the second is to deal with approximation in $\AA^n$ or $\PP^n$,  which has
been fully answered by the Schmidt subspace theorem; and the third (in paraphrase) is to generalize
Roth's theorem to projective varieties in a way which is compatible with unramified covers.
We feel that the results of this paper are a partial fulfillment of the third suggestion.  (We say partial
since Lang wanted a generalization of his ``geometric formulation'' of Roth's theorem, which applied to maps,
and since it is not completely clear to us what Lang intended by this suggestion.
Unfortunately we can no longer ask him.)

\medskip
\noindent
{\bf Other results.}
The proofs of the theorems (in particular Theorem \ref{thm:RothII}) hinge on a third invariant of a point and
ample line bundle $L$.  This invariant, $\Ar_{x}(L)$, is defined in \S\ref{sec:f-def} and further explored
in \S\ref{sec:more-about-Ar}.    This invariant is purely geometric in the sense that, like $\ep_{x}(L)$, it only
depends on the base change of $X$ to the algebraic closure.

This invariant is obtained by integrating a function $f(\gamma)$ which measures the ``relative asymptotic volume''
of the subspace of sections of $L$ vanishing to order $\geq \gamma$ at $x$.  One of the reasons for using $\Ar_{x}(L)$
is that the asymptotic behaviour of a line bundle is often better than any particular multiple.

In order to prove Theorem \ref{thm:RothII} we first prove an approximation result using $\Ar_{x}(L)$.

\noindent
{\bf Theorem (\ref{thm:RothI}):} {\em
Let $X$ be an irreducible variety over $\Spec(k)$.  Then
for any ample $\QQ$-bundle $L$ and any $x\in X(\Qbar)$ either

\begin{enumerate}
\item $\alpha_x(L) \geq \Ar_x(L)$

\hfill or \hfill\rule{0.1cm}{0.0cm}

\medskip
\item $\alpha_x(L)$ is computed on a proper subvariety of $X$.
\end{enumerate}
}

If $X$ is $n$-dimensional then there is an easy estimate $\Ar_{x}(L)\geq \frac{n}{n+1}\ep_{x}(L)$ (see
Corollary \ref{cor:beta-bound}) and so Theorem \ref{thm:RothI} immediately implies Theorem \ref{thm:RothII}.
It is interesting to study when $\Ar_{x}(L)=\frac{n}{n+1}\ep_{x}(L)$, i.e., when replacing $\Ar_{x}(L)$
by $\frac{n}{n+1}\ep_{x}(L)$ does not diminish the strength of the result.   Equivalent conditions
for this equality are given in Theorem \ref{thm:inequal-is-equal}.
The reader will also find a heuristic interpretation of $\Ar_{x}(L)$ in \S\ref{sec:more-about-Ar}.

Finally, we note that \S\ref{sec:unram-bounds} also proves that all theorems involving $\Ar_{x}(L)$ hold
with $\Ar_{x}(L)$ replaced by its limit $\Arhat_{x}(L)$ over unramified covers.

\medskip
\noindent
{\bf Remarks on the proof.}
The central motor of this paper, which largely implies the other approximation results,
is Theorem \ref{thm:simul-approx-I} to which \S\ref{sec:central-thm} is devoted.
This theorem is a simultaneous approximation theorem written in terms of
$\{\Ar_{x_v}(L)\}_{v\in S}$ where $S$ is a finite set of places of $k$, and $x_v\in X(\Qbar)$ for $v\in S$.
Theorem \ref{thm:simul-approx-I} is proved using the Faltings-W\"{u}stholz theorem
and the definition of $\Ar_{x}(L)$ has been chosen in order to optimize an estimate used in applying
that theorem.  The basic idea is explained at the beginning of the proof of Theorem
\ref{thm:simul-approx-I}, which appears at the end of \S\ref{sec:central-thm}.

The Faltings-W\"{u}stholz theorem implies Roth's theorem for $\PP^1$ and the Schmidt subspace theorem, and thus
the values of $\frac{n}{n+1}$ and $\frac{1}{2}$  when approximating on $\PP^n$ and $\PP^1$ respectively.
Our theorems (e.g., Theorem \ref{thm:RothII} applied to $\PP^n$) also produce these values, but we deduce them
from the Faltings-W\"{u}stholz theorem by a different method than their paper,
and it is worth commenting on this difference.

In the argument
of \cite[\S9]{FW} the value $\frac{n}{n+1}$ arises as the ratio of the dimension of the subspace of
$\Gamma(\PP^n,\Osh_{\PP^n}(1))$ vanishing at a point $x$, and the dimension of the entire space.
In our result the value $\frac{n}{n+1}$ arises as the integral
$\Ar_{x}(\Osh_{\PP^n}(1))=\int_{0}^{1} f(\gamma) \,d\gamma$
of the relative asymptotic volume function $f(\gamma)=1-\gamma^n$ for the line bundle $\Osh_{\PP^n}(1)$.
Thus --- as mentioned above as a motivation for $\Ar_{x}$ --- we deduce the constant $\frac{n}{n+1}$
from asymptotic properties of $\Osh_{\PP^n}(1)$ and not from its global sections.

\noindent
{\bf Organization of the paper.}
Sections \ref{sec:rat-approx}, \ref{sec:seshadri}, and \ref{sec:f-def} are devoted to the definitions
and basic properties of $\alpha_{x}(L)$, $\ep_{x}(L)$, and $\Ar_{x}(L)$ respectively.
In \S\ref{sec:central-thm} we prove Theorem 5.2, which will is used to prove all the other approximation
results in the paper.
In \S\ref{sec:Roth-theorems} we prove approximation results for a single
place, and in \S\ref{sec:simul-approx} we prove simultaneous approximation results for several places.
In \S\ref{sec:unram-bounds} we show that all of the previous theorems hold with $\Ar_x$ and $\ep_x$ replaced
by their suprema $\Arhat_{x}$ and $\ephat_{x}$ over unramified covers.  In \S\ref{sec:more-about-Ar} we
provide some complementary material about $\Ar_{x}(L)$, and finally in \S\ref{sec:Vojta} we give an elementary
application of our theorems to establish some previously unknown special cases of Vojta's main conjecture.

\noindent
{\bf Notation and Conventions.}  Unless otherwise specified we work over a fixed number field $k$.
By ``variety over $\Spec(k)$'' we mean a (possibly reducible, possibly singular) projective variety over $\Spec(k)$,
i.e., a reduced projective scheme over $\Spec(k)$.
We use additive notation for line bundles since this is in line with the behaviour
of $\alpha_x$, $\ep_x$, and $\Ar_{x}$.  On a product $X\times Y$ we therefore use  $L_1\squareplus L_2$
instead of $L_1\boxtimes L_2$
for a line bundle of the form $pr_{X}^{*}L_1 + pr_{Y}^{*}L_2$, with $pr_{X}$ and $pr_{Y}$ being the projections.

If $X$ is a variety over $\Spec(k)$, a point $x\in X(\Qbar)$ is a map $\Spec(\Qbar)\longrightarrow X$ of $k$-schemes.
Such a point gives rise to a point of $X\times_{k}\Qbar$, and a closed point (the image of this map) of $X$.
The symbol $\kappa(x)$ denotes the residue field of this closed point of $X$, called the field of definition of $x$.
We say that ``$x$ is defined over $K$'' if $\kappa(x)$ is a subfield of $K$ (this inclusion may be implicit).
A sequence of $k$-points of $X$ (or a sequence in $X(k)$) means an infinite sequence of distinct points of $X(k)$.
We denote such a sequence by $\xseq$ rather than $\xseq_{i\geq 0}$.

The absolute values are normalized with respect to $k$: if $v$ is a finite place of $k$, $\pi$ a uniformizer
of the corresponding maximal ideal, and $\kappa$ the residue field then $\nrm{\pi}_{v}=1/\#\kappa$;  if
$v$ is an infinite place corresponding to an embedding $i\colon k\hookrightarrow \CC$ then $\nrm{x}_{v} =
|i(x)|^{m_v}$ for all $x\in k$, where $m_v=1$ or $2$ depending on whether $v$ is real or complex.

Two real-valued functions $g$ and $g'$ with the same domain are called {\em equivalent} if there are positive
real constants $c \leq C$ so that $cg\leq g' \leq Cg$ for all values of the domain.  We will apply this terminology
in three situations: to distance functions $\sepv(\cdot,\cdot)$, to
height functions $H_{L}(\cdot)$, and to partially evaluated distance functions $\sepv(x,\cdot)$.
Typical domains are $X(\Qbar)\times X(\Qbar)$, $X(\Qbar)$, and Zariski open subsets or $v$-adically compact subsets
of these.

\noindent
{\bf Acknowledgements.}
We thank Chris Dionne, Laurence Ein, Robert Lazarsfeld, Victor Lozovanu, and Damien Roy for helpful discussions.
We are also extremely grateful to the referees of this paper for pointing out several mathematical and expositional
errors in the initial versions, and for their suggestions on how to correct them.
Finally, we wish to acknowledge an intellectual debt to Michael Nakamaye who has long advocated the point of view
that Seshadri constants are diophantine.

\section{Approximation by rational points}\label{sec:rat-approx}

Let $k$ be a number field, and $X$ a projective variety over $\Spec(k)$.
We begin by discussing the distance functions in the archimedean and non-archimedean cases.

\noindent
{\bf Distance Functions: Archimedean case.}
Fix an archimedean place $v_0$ of $k$, and an extension of $v_0$ to $\Qbar$, which we denote by $v$.
We choose a distance function
on $X(\Qbar)$ by choosing an embedding $X\hookrightarrow\PP^{r}_{k}$  and
pulling back (via $v$) the function on $\PP^{r}(\CC)\times \PP^{r}(\CC)$ given by the formula

\begin{equation}\label{eqn:distance-formula-arch}
\sepv(x,y) =
\left({1-\frac{|\sum_{i=0}^{r}{x_i\overline{y_i}}|^2}{(\sum_{i=0}^{r} |x_i|^2)
(\sum_{j=0}^{r} |y_j|^2)}}\right)^{[k_v:\RR]/2}
\end{equation}

\noindent
where $x=[x_0\colon\cdots\colon x_r]$, and $y=[y_0\colon\cdots\colon y_r]$ are points of $\PP^r(\CC)$, and
$|\cdot|$ is the absolute value on $\CC$ extending the usual absolute value on $\RR$, i.e,. such that
$|3+4\sqrt{-1}|=5$.

Note that if $k_v=\CC$ then this function does not satisfy the triangle inequality, nonetheless we continue to
call it a distance function.  (To see that this function does satisfy the triangle inequality if $k_v=\RR$
see \cite[Proposition 2.8.18]{BG}.)

\noindent
{\bf Distance Functions: Non-archimedean case.}
Fix a non-archimedean place $v_0$ of $k$, and an extension of $v_0$ to $\Qbar$, which we denote by $v$.
The place $v$ defines an absolute value $\nrm{\cdot}_v$ on $\Qbar$,
normalized according to our conventions in the introduction.
(This normalization agrees with the use of the symbol $\nrm{\cdot}_v$ in the books of Bombieri-Gubler
and Hindry-Silverman;  see \cite[1.3.6 and 1.4.3]{BG} and \cite[p.\ 171--172]{HS} respectively.)
We choose a distance function
on $X(\Qbar)$ by choosing an embedding $X\hookrightarrow\PP^{r}_{k}$  and
pulling back the distance function on $\PP^{r}(\Qbar)$ given by the formula

\begin{equation}\label{eqn:distance-formula-nonarch}
\sepv(x,y) =
\frac{\max_{0\leq i<j\leq r}(\nrm{x_iy_j-x_jy_i}_v)}{%
\max_{0\leq i\leq r}(\nrm{x_i}_v) \max_{0\leq j\leq r}(\nrm{y_j}_v)}
\end{equation}

\noindent
where $x=[x_0\colon\cdots\colon x_r]$, and $y=[y_0\colon\cdots\colon y_r]$ are points of $\PP^r(\Qbar)$.

\medskip
\noindent
{\bf Basic properties of distance functions.}
These definitions are somewhat opaque on first reading, but they are standard distance functions in Arakelov theory,
albeit normalized with respect to $k$, rather than $\QQ$.
(See for instance \cite[\S2.8]{BG} where a distance function $\delta_v(\cdot,\cdot)$ is defined for each place $v$;
the distance functions are related by the formula $\sepv(\cdot,\cdot)=\delta_v(\cdot,\cdot)^{[k:\QQ]}$.)
We will also briefly discuss the geometric meaning of $\sepv(\cdot,\cdot)$ for non-archimedean $v$ below.

We note two elementary properties of the distance function, whose proofs follow easily from the definitions.
\begin{proposition} Let $v$ be a place of $k$ extended to $\Qbar$, and $\sepv(\cdot,\cdot)$ the distance
function constructed by choosing an embedding $X\hookrightarrow \PP^{r}_k$.   Then

\begin{enumerate}
\item For all $x,y\in X(\Qbar)$ we have $\sepv(x,y)\in [0,1]$, with $\sepv(x,y)=0$ if and only if $x=y$.
\item If $K$ is a finite extension of $k$, then $\sepv(\cdot,\cdot)_{K}=\sepv(\cdot,\cdot)_{k}^{m_v}$, where
$m_v=[K_v:k_v]$ is the local degree.   (Here $\sepv(\cdot,\cdot)_K$ refers to the distance function defined
by using the same embedding and normalizing with respect to $K$
and $\sepv(\cdot,\cdot)_{k}$ the distance function normalized with respect to $k$, as above.)
\end{enumerate}
\end{proposition}

We will use the next result several times in proving equivalence of different types of distance functions.
Let $\CC_v$ be the completion of $\Qbar$ with respect to the place $v$.

\begin{lemma}\label{lem:equiv-on-compact-v-adic}
Let $Y$ be a variety over $\Spec(k)$, $U$ an affine open subset of $Y_{K}=Y\times _{k} K$ for some finite
extension $K/k$, and $u_1$,\ldots, $u_r$ and $u'_1$,\ldots, $u'_s$ two collections of elements of
$\Gamma(U,\Osh_Y)$ which generate the same ideal.   Then
the functions
$\max(\nrm{u_1(\cdot)}_v,\ldots, \nrm{u_r(\cdot)}_v)$ and
$\max(\nrm{u'_1(\cdot)}_v,\ldots, \nrm{u'_s(\cdot)}_v)$ are equivalent on any compact subset of $U(\CC_v)$.
\end{lemma}

\bpf
Since $u_1,\ldots, u_r$ and $u'_1$,\ldots, $u'_s$ generate the same ideal on $U$ there are functions
$f_{j,\ell}\in \Gamma(U,\Osh_{Y_{K}})$ such that
$u_j = \sum_{\ell=1}^{s} f_{j,\ell} u'_\ell$ for each $j=1$,\ldots, $r$.
Similarly there are functions $g_{\ell,j} \in \Gamma(U,\Osh_{Y_{K}})$
such that $u'_{\ell} = \sum_{j=1}^{r} g_{\ell,j} u_j$ for all $\ell=1$,\ldots, $s$.
On any compact subset $T$ of $U(\CC_v)$ the functions
$\nrm{f_{j,\ell}(\cdot)}_v$ and $\nrm{g_{\ell, j}(\cdot)}_v$ are bounded on $T$.
It follows that
the functions $\max(\nrm{u_1(\cdot )}_v,\ldots, \nrm{u_r(\cdot)}_v)$ and
$\max(\nrm{u'_1(\cdot)}_v,\ldots, \nrm{u'_s(\cdot)}_v)$ are equivalent on $T$.
\epf

\noindent
{\bf Remark.} \label{rem:same-local-degree}
Let $x$ be a point of $X(\Qbar)$ and let $K$ be the field of definition of $x$.
Throughout the paper we will be interested in approximating $x$ by points of $X(k)$.
If $K\not\subseteq k_v$, or equivalently, $K_v\neq k_v$ then it will be impossible to find a sequence of
points of $X(k)$ converging (in terms of $\sepv$) to $x$ (e.g., when $v$ is archimedean this happens when
$k_v=\RR$ and $K_v=\CC$).
Thus, in all cases we can approximate $x$ by points of $X(k)$ we may assume that $K_v=k_v$.

\begin{lemma}\label{lem:local-equiv-dist}
Let $V$ and $W$ be vector spaces over $k$,  $j\colon X\hookrightarrow \PP(V^{*})$ and $j'\colon X\hookrightarrow
\PP(W^{*})$ embeddings,
$\PP(V^{*})\cong \PP^r$ and $\PP(W^{*})\cong \PP^s$ choices of coordinates, and $\sepv$ and $\sepv'$ the
induced distance functions on $X$.  Let $K/k$ be any finite extension.
Then for any point $x\in X(K_v)$ and any inclusion $V\hookrightarrow W$
of $k$-vector spaces so that the resulting rational map $f\colon\PP^s \longdashrightarrow \PP^r$  is defined
at $j(x)$ and such that $f\circ j=j'$ near $x$, there is a compact $v$-adic neighbourhood $T$ of $(x,x)$ in
$X(K_v)\times X(K_v)$ such that $\sepv$ and $\sepv'$ are equivalent on $T$.
\end{lemma}

\bpf
Change of basis by $k$-linear transformation only changes the distance function by bounded amount
(see \cite[Theorem 3]{CV} for this statement for $\delta_v$).
We may therefore change coordinates and assume that the map $f$ is given by dropping the last $s-r$ coordinates
on $\PP^s$.  In the non-archimedean case we are therefore reduced to comparing the behaviour of
$$
\frac{\max_{0\leq i<j\leq s}(\nrm{x_iy_j-x_jy_i}_v)}{%
\max_{0\leq i\leq s}(\nrm{x_i}_v) \max_{0\leq j\leq s}(\nrm{y_j}_v)}
\,\,\,\mbox{and}\,\,\,
\frac{\max_{0\leq i<j\leq r}(\nrm{x_iy_j-x_jy_i}_v)}{%
\max_{0\leq i\leq r}(\nrm{x_i}_v) \max_{0\leq j\leq r}(\nrm{y_j}_v)}
$$
near $(x,x)$.  Let $X_0$,\ldots, $X_s$ and $Y_0$,\ldots, $Y_s$ be the coordinates on $\PP^s\times \PP^s$.
Choose an affine open $U$ containing
$x$ such that the embedding line bundle is trivial on $U$ and so we may identify sections with functions.
Since $f$ is defined at $x$ the functions $X_0$,\ldots, $X_r$ have no common zero at $x$, and so by shrinking $U$
we may assume that they generate the unit ideal on $U$.  The same is therefore true for the larger collection
of functions $X_0$,\ldots, $X_s$.  By Lemma \ref{lem:equiv-on-compact-v-adic} we thus have that
$\max_{0\leq i\leq s}(\nrm{x_i}_v)$ and $\max_{0\leq i\leq r}(\nrm{x_i}_v)$ are equivalent on any compact neighbourhood
of $x\in X(K_v)\subset X(\CC_v)$.  (Note that in the non-archimedean case these compact neighbourhoods
in $X(K_v)$ will have empty interior as a subset of $X(\CC_v)$.  Indeed, $\CC_v$ is not locally compact, so there
are no compact neighbourhoods of $x$ in $X(\CC_v)$ at all.)
A similar statement applies to the coordinates $Y_i$.  The functions
$\{X_iY_j-X_jY_i\}_{0\leq i<j \leq s}$ and
$\{X_iY_j-X_jY_i\}_{0\leq i<j \leq r}$ generate the ideal of the diagonal on $\PP^s\times \PP^s$ and
$\PP^r\times \PP^r$ respectively, and so restricted to $U\times U$ generate the ideal of the diagonal there.
Applying Lemma \ref{lem:equiv-on-compact-v-adic} again, the functions
$\max_{0\leq i<j\leq s}(\nrm{x_iy_j-x_jy_i}_v)$ and $\max_{0\leq i<j\leq r}(\nrm{x_iy_j-x_jy_i}_v)$
are equivalent on a compact neighbourhood of $(x,x)$ in $U(K_v)\times U(K_v)$.
This proves the lemma in the non-archimedean case.  In the archimedean case one uses the same strategy, the identity
$$
1-\frac{|\sum_{i=0}^{r}{x_i\overline{y_i}}|^2}{(\sum_{i=0}^{r} |x_i|^2)(\sum_{j=0}^{r} |y_j|^2)}
=
\frac{\sum_{0\leq i<j\leq r}|{x_iy_j-x_jy_i}|^2}{%
(\sum_{0\leq i\leq r}|{x_i}|^2) (\sum_{0\leq j\leq r}|{y_j}|^2)},
$$
and the fact that $\max(\nrm{f_1}_v,\ldots, \nrm{f_r}_v)$ and $(|f_1|^2+\cdots + |f_r|^2)^{[k_v:\RR]/2}$ are
equivalent for functions $f_1$,\ldots, $f_r$ taking values in $\CC$.
\epf

We now check that distance functions coming from two different embeddings are equivalent.  We are indebted to
one of the referees of this paper for the following argument, which is substantially simpler and shorter
than our original one.

\begin{proposition}\label{prop:equiv-dist}
Let $\sepv$ and $\sepv'$ be two distance functions coming from different embeddings of $X$.
Then for all finite extensions $K/k$, $\sepv$ is equivalent to $\sepv'$ on $X(K_v)\times X(K_v)$. 
\end{proposition}

\bpf
It suffices to show that for each $x\in X(K_v)$ there is a compact neighbourhood of $(x,x)$ in $X(K_v)\times X(K_v)$
where $\sepv$ and $\sepv'$ are equivalent.  Since $X$ is projective, $X(K_v)\times X(K_v)$ can be covered by finitely
many such neighbourhoods which then proves the proposition.

We first observe that we may assume that each of the embeddings is by a complete linear series.
Suppose that $L$ is a very ample line bundle and $V\subset W=H^0(X,L)$ a basepoint-free subseries such that the
associated map $j\colon X\longrightarrow \PP(V^{*})$ is an embedding.  Then $j'\colon X\hookrightarrow
\PP(W^{*})$ is also an embedding and the rational map $\PP(W^{*})\longdashrightarrow \PP(V^{*})$ arising from
the inclusion $V\hookrightarrow W$ is defined everywhere along $j(X)$.
Thus the result we want follows from Lemma \ref{lem:local-equiv-dist}.

Now let $\sepv$ and $\sepv'$ be two distance functions coming from embeddings
$j\colon X\hookrightarrow \PP(V^{*})$ and $j'\colon X\hookrightarrow \PP(W^{*})$, with
$V=H^0(X,L)$ and $W=H^0(X,L')$ where $L$ and $L'$ are very ample line bundles.
Assume that $L-L'$ is basepoint free.  Then for any point $x\in X(K_v)$
there is a section $s$ of $L-L'$ that does not vanish at $x$. Multiplication by $s$ induces a rational map
$f_s\colon \PP(V^{*})\longdashrightarrow \PP(W^{*})$ that is defined at $j(x)$ and such that $f_s\circ j=j'$ near $x$.
Thus the neighbourhood we want is again guaranteed by Lemma \ref{lem:local-equiv-dist}.

Finally, for general $L$ and $L'$, we may replace $L$ with a multiple $mL$ large enough so that $mL-L$ and $mL-L'$
are both basepoint free, and the proposition immediately follows.
\epf

\np
We next turn to local descriptions of the distance function useful in computations.

\begin{lemma}
\label{lem:local-distance-Zariski-open}
Let $x$ be a point of $X(\Qbar)$ and $K$ any finite extension of $k$ over which $x$ is defined.
Then there exists an open affine subset $U$ of
of $X_{K}:=X\times_{k}K$   containing $x$, and elements
$u_1$, \ldots, $u_r$ of $\Gamma(U,\Osh_{X_K})$
which generate the maximal ideal of $x$ and positive real constants $c\leq C$ such that

\begin{equation}\label{eqn:local-distance-equiv}
c\,\sepv(x,y)\leq \min\left({1,\max\left(\nrm{u_1(y)}_v,\ldots, \nrm{u_m(y)}_v\right)\rule{0cm}{0.4cm}}\right)
\leq C\, \sepv(x,y)
\end{equation}
\noindent
for all $y\in U(K_v)$.  That is, on $U(K_v)$ the function
$\min(1,\max\left(\nrm{u_1(\cdot)}_v,\ldots, \nrm{u_r(\cdot)}_v\right))$ is equivalent to the function $\sepv(x,\cdot)$.
\end{lemma}

\bpf
We start with two reductions.
First, since the absolute value $\nrm{\cdot}_v$ and the distance function $\sepv(\cdot,\cdot)$ transform the same
way under field extensions, we may assume that $x$ is defined over $k$.  Second, by Proposition \ref{prop:equiv-dist}
we may choose whichever embedding of $X$ we wish when performing the calculation.

Given these reductions, choose an embedding $X\hookrightarrow\PP^r$ so that $x$ is sent to
$[1\colon0\colon \cdots \colon 0]$.
Let $Z_0$,\ldots, $Z_r$ be homogenous coordinates on $\PP^r$, and choose
the open set $U$ of $X$ to be the set $Z_0\neq 0$, and $u_i=Z_i/Z_0$ for $i=1$,\dots, $r$ as the generators of
the maximal ideal of $x$.  If $v$ is non-archimedean, then \eqref{eqn:distance-formula-nonarch} and
the fact that $x$ is sent to $[1\colon0\colon\cdots\colon0]$ give
$$\sepv(x,y) = \frac{\max(\nrm{Z_1(y)}_v,\nrm{Z_2(y)}_v,\ldots, \nrm{Z_r(y)}_v)}%
{\max(\nrm{Z_0(y)}_v,\nrm{Z_1(y)}_v,\ldots, \nrm{Z_r(y)}_v)}\,\,\,\mbox{for all $y\in X(\Qbar)$}.$$
For $y\in U(K_v)$, this is equal to $\min(1,\max\left(\nrm{u_1(y)}_v,\ldots, \nrm{u_r(y)}_v\right))$.

In the case that $v$ is archimedean, we may further assume that $k_v=\CC$, again using the fact that the
functions to be compared transform the same way under field extensions.  From \eqref{eqn:distance-formula-arch}
and the fact that $x$ is sent to $[1\colon0\colon\cdots\colon0]$ we obtain
$$\sepv(x,y) = 1-\frac{|y_0|^2}{|y_0|^2+\cdots + |y_r|^2}
= \frac{|y_1|^2+\cdots + |y_r|^2}{|y_0|^2+|y_1|^2+\cdots + |y_r|^2}
= \frac{\nrm{y_1}_v+\cdots + \nrm{y_r}_v}{\nrm{y_0}_v+\nrm{y_1}_v+\cdots + \nrm{y_r}_v}.
$$
For $y\in U(K_v)$, $y_0\neq 0$, and $u_j(y) = y_j/y_0$ for $j=1$,\ldots, $r$.
Thus $\sepv(x,y) = \frac{\nrm{u_1(y)}_v+\cdots+\nrm{u_r(y)}_v}{1+\nrm{u_1(y)}_v+\cdots+\nrm{u_r(y)}_v}$;
it is then elementary to check that
\eqref{eqn:local-distance-equiv} holds with $c=\frac{1}{r}$ and $C=2$.
\epf

We will need an extension of this lemma which applies to any affine open $U$ containing $x$, and
any choice $u_1$,\ldots, $u_r$ of generators the maximal ideal.  To do this we need to drop the requirement
that the result hold for all $y\in U(K_v)$, and restrict to points $v$-adically close to $x$; for our purposes
it will be sufficient to restrict to sequences converging to $x$.

\begin{lemma}
\label{lem:local-distance}
Let $x$ be a point of $X(\Qbar)$ and $K$ any finite extension of $k$ over which $x$ is defined.
Let  $U$ be any open affine subset
of $X_{K}:=X\times_{k}K$  containing $x$.
Let $u_1$, \ldots, $u_r$ be any elements of $\Gamma(U,\Osh_{X_K})$
which generate the maximal ideal of $x$.  Then for any sequence $\xseq$ of points of $U(K_v)$ such that
$\sepv(x,x_i)\to 0$ as $i\to \infty$ the functions $\sepv(x,\cdot)$ and
$\max(\nrm{u_1(\cdot )}_v,\ldots, \nrm{u_r(\cdot)}_v)$ are equivalent on $\xseq$.  In other words, there
are positive constants $c < C$ such that for all $i\geq 0$ we have
$$
c\,\sepv(x,x_i)\leq \max\left(\nrm{u_1(x_i)}_v,\ldots, \nrm{u_m(x_i)}_v\right) \leq C\, \sepv(x,x_i).
$$
\end{lemma}

\bpf
By Lemma \ref{lem:local-distance-Zariski-open} there is an affine open set $U'$ containing $x$
and functions $u'_1$,\ldots, $u'_{s'}$ generating the maximal ideal of $x$ (on $U'$) such that
$\min(1,\max\left(\nrm{u'_1(\cdot)}_v,\ldots, \nrm{u'_s(\cdot)}_v\right))$
is equivalent to $\sepv(x,\cdot)$ on $U'$.  By replacing $U$ and $U'$ with their
intersection and possibly omitting initial members of the sequence we may assume $U=U'$.
Since $\sepv(x,x_i)\to 0$
as $i\to \infty$, each of the $\nrm{u'_j(x_i)}_v\to 0$ as well, and so $\sepv(x,\cdot)$ is equivalent to
$\max\left(\nrm{u'_1(\cdot)}_v,\ldots, \nrm{u'_s(\cdot)}_v\right)$ on $\xseq$.   It therefore suffices
to show the equivalence of $\max(\nrm{u_1(\cdot)}_v,\ldots, \nrm{u_r(\cdot)}_v)$ and
$\max(\nrm{u'_1(\cdot)}_v,\ldots, \nrm{u'_s(\cdot)}_v)$.
Since $\sepv(x,x_i)\to 0$ as $i\to \infty$, for large enough $i$ the points $x_i$ are contained in a compact
neighbourhood of $x$ in $U(K_v)$.  Thus the equivalence follows by
Lemma \ref{lem:equiv-on-compact-v-adic}.
\epf

\noindent
One warning: Lemma \ref{lem:local-distance}, with the freedom to choose $U$ and $u_1$,\ldots, $u_r$ does not
hold in the generality of Lemma \ref{lem:local-distance-Zariski-open}, i.e., for all $y\in U(K_v)$.
As an elementary example, let $k=\QQ$,
$K=\QQ(\sqrt{2})$, $X=\PP^1$, $x=[-\sqrt{2}:1]$, and let $v$ be an extension of the archimedean absolute value
on $\QQ$.  Let $\AA^1_K=\Spec(K[t])$ be the open affine subset of $\PP^1$ obtained by removing $[1:0]$,
and $U$ the open subset of $\AA^1_K$ obtained by removing the point $[\sqrt{2}:1]$.
Then on $U$ the function $u_1=t^2-2$ generates the maximal ideal of $x$.
Pick a sequence of points $\xseq$ in $\AA^1(\QQ)$ converging to $[\sqrt{2}:1]$.
Then $\sepv(x,x_i)$ does not go to zero as $i\to \infty$, while $\nrm{u_1(x_i)}_v$ does, so the
two functions are not equivalent.

The issue is fairly clear.
The function $u'_1=t+\sqrt{2}$ also generates the maximal ideal of $x$ on $U$, and $\min(1,\nrm{u'_1(\cdot)}_v)$
is equivalent to the distance function $\sepv(x,\cdot)$ on $U(\RR)$.  The problem is that although
$u_1$ and $u'_1$ satisfy the relation $u'_1=\frac{1}{t-\sqrt{2}} u_1$ we cannot deduce that
$u'_1(x_i)\to 0$ as $i\to \infty$ from the fact that
$u_1(x_i)\to 0$ because the function $\frac{1}{t-\sqrt{2}}$ is unbounded on $\xseq$.

Lemma \ref{lem:local-distance-Zariski-open} gives one justification that the distance functions chosen are reasonable:
they are globally defined functions which locally, around any point $x$, behave like the standard $v$-adic
distance functions induced from an embedding into an affine space.
The following discussion connecting the non-archimedean distance functions to order of contact provides
another justification.  This discussion is not necessary for any of the arguments in the paper, but is
included to provide further geometric intuition behind the definition.

\noindent
{\bf Geometric meaning of non-archimedean distance.}
The formula in \eqref{eqn:distance-formula-nonarch} is a compact way of stating a very concrete notion of
$v$-adic distance: points $x$ and $y$ are close if the corresponding curves in an integral model of $X$ have high
order of contact at the place $v$.

To see this, we will define a distance function, $\sepv'(\cdot,\cdot)$ via order of contact, suitably normalized,
and show that it equals $\sepv(\cdot,\cdot)$.  Let  $\imX$ be the projective integral model of $X$ over
$\Spec(\Osh_{k})$ obtained by taking the closure of $X$ in $\PP^{m}_{\Osh_{k}}$, under the inclusions
$X\hookrightarrow  \PP^{m}_k\hookrightarrow \PP^{m}_{\Osh_{k}}$.  Let $\Okv$ be
the completion of $\Osh_{k}$ at the maximal ideal corresponding to $v$, and set
$\imXv$ to be the base-change of $\imX$ to $\Okv$.

Suppose that $x,y\in X(k)$. Then $x$ and $y$ give rise to sections $\sigma_x$ and $\sigma_y$ of $\imXv$ over
$\Spec(\Okv)$.  If $x=y$ we set $\sepv'(x,y)=0$. If $x\neq y$, then let $Z$ be the scheme of intersection of
$\sigma_x$ and $\sigma_y$ in $\imXv$.
The ring $\Gamma(Z,\Osh_{Z})$ of global sections of the structure sheaf of $Z$ has finitely many elements, and we
set $\sepv'(x,y) = 1/(\#\Gamma(Z,\Osh_{Z}))$ where $\#$ denotes the number of elements in the ring.  Note that if $Z$
is empty then $\Gamma(Z,\Osh_{Z})$ is the zero ring with a single element (namely $0$). I.e., if $Z=\emptyset$ then
$\sepv'(x,y)=1$.

In the general case that $x,y\in X(\Qbar)$, let $F/k$ be any finite extension so that $x$ and $y$ are defined over $F$
and set $\imX_{F,v}$ to be the base change of $\imX$ to $\Spec(\OFv)$,
where $\OFv$ is the completion of $\Osh_{F}$ at $v$.
As before, $x$ and $y$ give rise to sections $\sigma_{x}$ and $\sigma_{y}$ of $\imX_{F,v}$ over $\Spec(\OFv)$.
If $x=y$ then set $\sepv'(x,y)=0$.  Otherwise let $Z$ be the scheme of intersection, and set
$\sepv'(x,y) = 1/(\#\Gamma(Z,\Osh_{Z}))^{\frac{1}{[F_v:k_v]}}$.  If $F'/F$ is any finite extension, and $Z'$ the scheme
of intersection of the corresponding sections $\sigma'_{x}$ and $\sigma'_{y}$ of $\imX_{F',v}$ then
$\#\Gamma(Z',\Osh_{Z'})=\#\Gamma(Z,\Osh_{Z})^{[F'_v\colon F_v]}$.  It follows that $\sepv'$ is well defined.

To see that $\sepv'(\cdot,\cdot)$ is equal to $\sepv(\cdot,\cdot)$
we make the following observations: (1) Since both functions transform in the same way when extending the field,
we may assume that $x$ and $y$ are defined over $k$.  (2) The section $\sigma_x$ is obtained by multiplying the
coordinates of $x$ by an element in $k$ so that all coordinates are in $\Osh_{k,v}$ and such that at least one
coordinate is not in the maximal ideal corresponding to $v$.
After multiplying, we have $\max_{0\leq i\leq m} (\nrm{x_i}_v)=1$.  Similarly, we may assume that the section
$\sigma_y$ is given by $[y_0:\cdots:y_m]$ and that $\max_{0\leq j\leq m} (\nrm{y_j}_v)=1$.
(3) The diagonal of $\PP^m\times\PP^m$ is cut out by
the equations $\{X_iY_j-X_jY_i\}$ for $0\leq i< j\leq m$, where $X_0$,\ldots, $X_m$ and $Y_0$,\ldots, $Y_m$
are the coordinates on the product. Thus the ideal of $\OFv$ generated by
$\{x_iy_j-x_jy_i\}_{0\leq i<j\leq m}$ is the ideal of the scheme of intersection $Z$.
(4) Our normalization for $\nrm{\cdot}_v$ now shows that
$\max_{0\leq i< j\leq m}(\nrm{x_iy_j-x_jy_i}_v)=1/\#\Gamma(Z,\Osh_{Z})$, i.e., that $\sepv'(x,y)=\sepv(x,y)$.

\noindent
{\bf Height Functions.}
A height function is a function $H\colon X(\Qbar) \rightarrow \RR_{>0}$.
Two height functions $H$ and $H'$ are equivalent if there are positive real constants $c$ and $C$ with $0<c\leq C$ such
that $$ c\,H(x) \leq  H'(x) \leq C\,H(x)$$ for all $x\in X(\Qbar)$ (see also ``Notations and Conventions'' in the
introduction).  The set of height functions forms a group under multiplication and the group operation descends
to equivalence classes of height functions.

For any line bundle $L$ on $X$ we may associate a height function $H_{L}$, well defined up to
equivalence, in such a way that the map from $\Pic(X)$ to the equivalence classes of height functions is
a group homomorphism and the height function is functorial with respect to pullbacks.
For details on how to do this, see for example any one of
\cite[Chap.\ 2]{BG}, \cite[Part B]{HS}, \cite[Chap.\ III]{La}, or \cite[Chap.\ 2]{Se}.
One caveat: the normalizations used in these references are not all the same.  In this paper we normalize our height
functions so that for a point
${x}=[x_0:\cdots :x_r]\in\PP^r(k)$, the height with respect to $\Osh_{\PP^r}(1)$ is
\[H({x})=\prod_v \max(\nrm{x_0}_{v},\ldots, \nrm{x_r}_{v})\]
where the product ranges over all the places $v$ of $k$.
Unless otherwise specified all height functions in this paper are multiplicative and relative to $k$.

\noindent
{\bf Approximation Constants.}
We now define the main objects of study in this paper, inspired by similar definitions from \cite{McK}.
We fix a single place $v$, archimedean or non-archimedean, and a corresponding distance function $\sepv$
as described above.

\begin{definition}\label{seqappconst}
Let $X$ be a projective variety, $x\in X(\Qbar)$, $L$ a line bundle
on $X$.  For any sequence $\xseq\subseteq X(k)$ of distinct
points with $\sepv(x,x_i)\rightarrow 0$ (which we denote by $\xseqtox$), we set
$$A(\xseq, L) = \left\{{
\gamma\in\RR \st
\sepv(x,x_i)^{\gamma} H_{L}(x_i)\,\,\mbox{is bounded from above}
}\right\}.
$$
If $\xseq$ does not converge to $x$ then we set $A(\xseq,L)=\emptyset$.
\end{definition}

\noindent
{\bf Remarks.}
(a) It follows easily from the definition that if $A(\xseq, L)$ is nonempty then it is an
interval unbounded to the right, i.e., if $\gamma\in A(\xseq,L)$ then $\gamma+\delta\in A(\xseq,L)$ for any $\delta>0$.

(b) If $\xsubseq$ is a subsequence of $\xseq$ then $A(\xseq,L)\subseteq A(\xsubseq,L)$.

\begin{definition}
For any sequence $\xseq$ we set $\alpha_{x}(\xseq,L)$ to be the infimum of $A(\xseq,L)$
(in particular if $A(\xseq,L)=\emptyset$ then $\alpha_x(\xseq,L)=\infty$).
We call $\alpha_{x}(\xseq,L)$ the approximation constant of $\xseq$ with respect to $L$.
\end{definition}

It follows immediately from the definition that for any $\delta>0$,
$\sepv(x,x_i)^{\alpha_{x}(\xseq,L)+\delta}H_{L}(x_i)\to 0$ as $i\to\infty$ whenever $\alpha_x$ is finite.
We will frequently use this fact.
By remark (b) above, if $\xsubseq$ is a subsequence of $\xseq$ then $\alpha_{x}(\xsubseq,L)\leq \alpha_{x}(\xseq,L)$.

\begin{definition}\label{approx}
The approximation constant $\alpha_{x,X}(L)$ of $x$ with respect to
$L$ is defined to be the infimum of all approximation constants of
sequences of points in $X(k)$ converging to $x$.
If no such sequence exists, we set $\alpha_{x,X}(L)=\infty$.
\end{definition}

\noindent
{\bf Remarks.}
(a) The asymptotics of the approximation are unchanged if we replace the distance and height functions by
equivalent ones.
Since the approximation constant $\alpha_x$ is local to $x$,
we are also free to replace the distance function by one which is
only equivalent to $\sepv$ in some open set (in the analytic, $v$-adic, or Zariski topology)
around $x$ without changing $\alpha_x$.
In particular, by Proposition \ref{prop:equiv-dist}
the definition of $\alpha_x$ does not depend on the choice of projective embedding
used to define $\sepv$.

(b) Slightly more generally,
two height functions $H$ and $H'$ are called quasi-equivalent
if for every $\delta>0$ there exist $0<c<C$ (depending on $\delta$) so that
$$c\,H^{1-\delta} \leq H' \leq C\,H^{1+\delta}.$$
The definitions of $\alpha_{x}(\xseq,L)$ and $\alpha_{x}(L)$ only depend on the quasi-equivalence class of the height
function.  For ample $L$ and any $M\in \Pic^0(X)$, the heights $H_{L}$ and $H_{L\otimes M}$ are quasi-equivalent
(see \cite[p.\ 26]{Se};  the proof also applies to singular varieties).
For ample $L$, the functions $\alpha_{x}(L)$
and $\alpha_{x}(\xseq,L)$ therefore only depend on the class of $L$ in $\Pic(X)/\Pic^{0}(X)$, i.e., on the class
of $L$ in the N\'eron-Severi group.

(c) If $L$ is ample, then there exists $c>0$ so that $H_{L}(x_i)\geq c$ for all
$x_i\in X(k)$.  Thus if the sequence
$\sepv(x,x_i)^{\gamma}H_L(x_i)$ is bounded we must have $\gamma\geq 0$.  We therefore conclude that $\alpha_x(L)\geq 0$.
In Proposition \ref{prop:alpha}(d) we will show the slightly stronger statement  $\alpha_{x}(L)>0$ for ample $L$.
Similarly, if some multiple of $L$ is an effective divisor and $x$ a point outside the asymptotic base locus $Z$ of $L$, we can again
conclude that $\alpha_x(L)\geq 0$, since again $H_L(x_i)\geq c>0$ for all $x_i\in X(k)\setminus Z(k)$.

(d) When $L$ is ample, $H_{L}(x_i)$ is a proxy for how complicated the point $x_i$ is.  The number
$\alpha_{x}(\xseq,L)$ therefore measures the cost (in terms of the growth of complexity of the approximating points)
required to get closer and closer to $x$.
Thus under this definition (for ample $L$) smaller approximation constants correspond to
better approximating sequences. 

(e) It is possible that $\alpha_x(L)=\infty$. This occurs if either there is no sequence of points in $X(k)$ converging
to $x$, or, if for every such sequence $\xseq$ the set $A(\xseq,L)$ is empty.  It is also possible that
$\alpha_x(L)=-\infty$. This can occur in either of the ways suggested by the definition.
For instance there may be one sequence $\xseq$ so that
$A(\xseq,L) = (-\infty,\infty)$. Alternatively given any $C>0$, there may be a sequence $\xseq$ such that
$\alpha_{x}(\xseq,L) < - C$.  This happens, for instance on $\PP^{n}$ with $L=\Osh_{\PP^n}(-1)$ and $x\in \PP^n(k)$.
See later comments and examples for more on these extreme situations.

(f) The definition given above is different from the definition of the ``approximation
constant'' given in \cite{McK}, since it is the infimum of the set
described above rather than the minimum, as in \cite{McK}.  In
\cite{McK} this difference is not important
to the results, since in all examples that appear in that paper, the
minimum exists and is equal to the infimum.

More significantly, the distance function used in \cite{McK} is
computed with respect to all of the archimedean places of $k$, rather
than a single archimedean or non-archimedean place, and is not normalized by local degree.  Thus, when
$k=\QQ$ and we choose the archimedean place,  this is no difference at all, but in general the distance
functions will be different.  Where necessary, we will reprove results from \cite{McK}
using the new definitions.

\vspace{.1in}

We next give an alternate characterization of $\alpha_x(L)$, valid for those line bundles whose heights satisfy the
Northcott property, similar to the usual definition of the approximation constant on the affine line.  Recall that a
line bundle $L$ has the Northcott property if for any constant $c\in\RR$, the set of points $y\in X(k)$ such that
$H_L(y)\leq c$ is finite.  Note in particular that every ample line bundle has the Northcott property.

\begin{definition}\label{def:B-Northcott}
For any point $x\in X(\Qbar)$ and any line bundle $L$ we set
$$B_x(L) = \left\{{
\gamma\in \RR_{\geq 0} \left|{\,\,\,\,
\mbox{\begin{minipage}{8cm}
for all $C>0$ the number of $x_i\in X(k)$ such that $\sepv(x,x_i)^{\gamma}H_{L}(x_i) < C$
is finite.
\end{minipage}}}\right.\,\,
}\right\}$$
\end{definition}

\noindent
{\bf Remarks.}

\begin{itemize}
\item[(a)] $0\in B_{x}(L)$ if and only if $L$ has the Northcott property.
\item[(b)] $B_{x}(L)\neq\emptyset$ if and only if $L$ has the Northcott property.
\item[(c)] $B_{x}(L)$ (if nonempty) is an interval: if $\gamma\in B_{x}(L)$ then $\gamma-\delta\in B_{x}(L)$
for all $0\leq\delta\leq \gamma$.
\end{itemize}

Part (a) is obvious from the definition. For part (b), if $L$ has the Northcott property then $B_{x}(L)$
is nonempty by (a).
If $L$ does not have the Northcott property then there is a constant $C$ so that the number of $x_i\in X(k)$
with $H_L(x_i)<C$ is infinite.   Since
$\sepv(x,x_i)\leq 1$, for any $\gamma>0$ these infinitely many $x_i$ also satisfy
$$\sepv(x,x_i)^{\gamma} H_{L}(x_i) \leq  H_{L}(x_i) <  C $$
and therefore $\gamma\not\in B_{x}(L)$.  Thus $B_{x}(L)$ is empty.
Part (c) follows by again using the fact that $\sepv(x,x_i)$ is bounded.

We remark that there are line bundles which have the Northcott property but
which are not ample.  For instance, let $X$ be the blowup of $\PP^2$ at the base locus of a $k$-rational pencil of
plane curves of genus at least three.
There is a morphism $\pi\colon X\to\PP^1$ whose fibres are exactly the curves in the pencil.  If the pencil is
chosen so that the singular fibres all have a single nodal singularity and the curves in the pencil intersect
transversely at smooth points (as is the case for a generic pencil),
then every fibre of $\pi$ contains finitely many $k$-rational points, by Faltings' Theorem.
Thus, the height associated to the nef line bundle $\pi^*\Osh_{\PP^1}(1)$ satisfies the Northcott property,
but is not ample.

\begin{proposition}\label{prop:equiv-alpha}
Suppose that $L$ has the Northcott property.  Then $\alpha_x(L) = \sup(B_{x}(L))$.
\end{proposition}

\bpf
Set $\alpha_x=\alpha_x(L)$ and $b_x=\sup(B_{x}(L))$.  By definition of $\alpha_x$, for any $\delta>0$ there
exists a sequence $\xseq$ such that $\alpha_{x}(\xseq,L)<\alpha_x+\delta$ and hence (by the definition of
$\alpha_{x}(\xseq,L)$) we conclude that $\sepv(x,x_i)^{\alpha_x+\delta}H_{L}(x_i)$ is bounded.  Therefore
$\alpha_x+\delta\notin B_{x}(L)$ and so $\alpha_x+\delta\geq b_x$.  Letting $\delta$ go to zero we conclude
$\alpha_x\geq b_x$.

On the other hand, by the definition of $b_x$, for any $\delta>0$ there is a $C$ such that there are infinitely
many solutions $x_i\in X(k)$ to $\sepv(x,x_i)^{b_x+\delta}H_{L}(x_k)<C$.  Since $L$ has the Northcott property,
the set of heights $H_{L}(x_i)$ must be unbounded, and we can therefore choose a subsequence $\{x_i\}$ of these
points so that $H_{L}(x_i)\to\infty$ as $i\to\infty$.  By the boundedness of the product, we conclude that
$\sepv(x,x_i)\to 0$, and so $\{x_i\}$ converges to $x$.  But then
$$\sepv(x,x_i)^{b_x+2\delta}H_{L}(x_i) < C\cdot \sepv(x,x_i)^{\delta} \to 0$$
and so $b_x+2\delta\in A(\xseq,L)$.  Thus $b_x+2\delta\geq \alpha_x$, and letting $\delta$ go to zero
we conclude that $b_x\geq \alpha_x$ and so $\alpha_x=b_x$. \epf

\noindent
{\bf Remark.}  If $L$ has the Northcott property then $0\in B_{x}(L)$ and hence $\alpha_x(L)=\sup(B_{x}(L))\geq 0$
by Proposition \ref{prop:equiv-alpha}.  In particular this shows again that for ample bundles $\alpha_x(L)\geq 0$.

\vspace{.1in}

It will be useful to know how the approximation constant changes when we change the field $k$.   We use
the notation that for an extension field $K/k$, $\alpha_{x}(\xseq,L)_{K}$
(respectively $\alpha_{x}(L)_{K}$) denotes the
approximation constant of a sequence (resp.\ point $x$) computed with respect to $K$.   This means that when
computing $\alpha$, we use the height $H_{L}$ relative to $K$ and normalize $\sepv$ relative to $K$.
If $d=[K\colon k]$ and $m_v=[K_v\colon k_v]$ then this means simply that $H_{L}(x_i)_{K} = H_{L}(x_i)^{d}_{k}$
and $\sepv(x,x_i)_{K} = \sepv(x,x_i)^{m_v}_{k}$.

\begin{proposition}\label{prop:change-of-field}
Suppose $x\in X(\Qbar)$, $L$ a line bundle on $X$, and $\xseqtox$ a sequence of points in $X(k)$
approximating $x$.   Let $K$ be any finite extension of $k$.  Then $\xseqtox$ can also be considered
to be a set of points of $X(K)$ approximating $x$.   Set $m_v=[K_v\colon k_v]$, and let $d=[K\colon k]$.
Then
$$\alpha_x\left({\xseq, L}\right)_{K} =
\frac{d}{m_v} \alpha_x\left({\xseq, L}\right)_{k}.$$
In particular, we have the bound $\alpha_x(L)_{K} \leq \frac{d}{m_v} \alpha_x(L)_{k}$.
\end{proposition}

\bpf
The claim that
$\alpha_x\left({\xseq, L}\right)_{K} =
\frac{d}{m_v} \alpha_x\left({\xseq, L}\right)_{k}$ follows immediately from the equalities
$H_{L}(\cdot)_{K} = H_{L}(\cdot)_{k}^{d}$
and $\sepv(\cdot,\cdot)_K=\sepv(\cdot,\cdot)_k^{m_{v}}$.  The inequality
$\alpha_x(L)_{K} \leq \frac{d}{m_v} \alpha_x(L)_{k}$ then follows since the sequences of $k$-points
approximating $x$ are a subset of the sequences of $K$-points approximating $x$.
\epf

\vspace{.1in}

\noindent
{\bf Basic properties of $\alpha$.} We start by computing $\alpha$ when $x\in\PP^n(k)$.

\begin{lemma}\label{lem:projective}
Let $x$ be any $k$-point of $\PP^n$.    Then $\alpha_{x,\PP^n}(\Osh_{\PP^n}(1))=1$.
\end{lemma}

\noindent
{\it Proof:}  Without loss of generality, we may assume that
$x=[1:0:\ldots:0]$.  We first show that $\alpha_{x}(\xseq,\Osh_{\PP^n}(1))\geq 1$ for all sequences $\xseq$ of
$k$-points.
Let $Z_0$,\ldots, $Z_n$ be the coordinates on $\PP^n$ and
$\xseq$ a sequence of $k$-points converging to $x$.  Since $\sepv(x,x_i)\to0$ as $i\to \infty$
we conclude that $\nrm{Z_j(x_i)/Z_0(x_i)}_{v}\to 0$ for each $j=1$,\ldots, $n$.  By passing to a subsequence
of the $x_i$, which can only possibly lower the value of $\alpha$, we may assume that for all $i$ we have that
$\nrm{Z_0(x_i)}_v$ is the largest of the $\nrm{Z_j(x_i)}_v$  and that there is a
fixed $j\in \{1,\ldots, n\}$ so that $\max(\nrm{Z_1(x_i)}_v,\ldots,\nrm{Z_n(x_i)}_v) = \nrm{Z_j(x_i)}_v$.
By Lemma \ref{lem:local-distance} we then have $\sepv(x,x_i) = \nrm{Z_j(x_i)/Z_0(x_i)}_v$ for all $i$ (at least
up to equivalence).
Thus, for any $\gamma\geq 0$

\begin{eqnarray*}
\sepv(x,x_i)^{\gamma}H(x_i) & = &
\left(\frac{\nrm{Z_j(x_i)}_v}{\nrm{Z_0(x_i)}_v}\right)^{\gamma} \cdot \nrm{Z_0(x_i)}_v
\cdot \prod_{w\neq v} \max(\nrm{Z_0(x_i)}_w,\ldots, \nrm{Z_n(x_i)}_w) \\
& \geq &
\left(\frac{\nrm{Z_j(x_i)}_v}{\nrm{Z_0(x_i)}_v}\right)^{\gamma} \cdot \nrm{Z_0(x_i)}_v
\cdot \prod_{w\neq v} \nrm{Z_j(x_i)}_w = \left(\frac{\nrm{Z_0(x_i)}_v}{\nrm{Z_j(x_i)}_v}\right)^{1-\gamma},
\end{eqnarray*}

\noindent
where in the last step we have used the product formula.  If $\gamma <1$ then the lower bound above goes to
infinity as $i\to\infty$, and hence $\alpha_{x}(\xseq,\Osh_{\PP^n}(1))\geq 1$.

We next show that we can achieve $\alpha=1$.   Since we can always choose to approximate along a rational
line containing $x$ it suffices to treat the case $n=1$ and approximate the point $[1\colon 0]$.  We
will handle the archimedean and non-archimedean cases separately.

In the archimedean case embed $\Osh_k$ as a lattice in the Minkowski space $\prod_{w\,\scriptstyle{\rm arch}} k_w = \RR^r\times\CC^s$.  For any $D>0$, there are infinitely many elements of $\Osh_k$ that lie in the cylinder $\{b\in\Osh_k\mid \nrm{b}_w\leq D\,\mbox{for}\,w\neq v\}$.  These elements $b_i$ satisfy $H([b_i:1])\sepv([1:0],[b_i:1])\leq D^{r+2s}$, and so for the sequence $x_i=[b_i\colon 1]$ we conclude that $\alpha_{x}(\xseq,\Osh_{\PP^1}(1))\leq 1$, and therefore that $\alpha_{x}(\xseq,\Osh_{\PP^1}(1))=1$.

In the non-archimedean case, since the
ideal class group is finite some power of the maximal ideal corresponding to $v$ is principal, generated by
$b\in \Osh_{k}$.  Thus
we have $\nrm{b}_v < 1$ and $\nrm{b}_w=1$ for all other finite places $w$ of $k$.  After taking a further power of $b$, and multiplying by a suitably chosen unit, we may suppose in addition that $\nrm{b}_w>1$ for all infinite places $w$.
Set $x_i=[1\colon b^{i}]$ for $i\geq 0$.  Then $H(x_i) = \prod_{w\, \scriptstyle{\rm arch}} \nrm{b^i}_w =
1/\nrm{b}^i_v$, where the last equality follows from the product formula.  Since $\sepv(x,x_i) = \nrm{b}_v^i$,
it is clear that $\alpha_x(\xseq,\Osh_{\PP^1}(1))=1$ for this sequence. \epf

\vspace{.1in}

The next proposition collects some elementary properties of $\alpha$.

\begin{proposition}\label{prop:alpha}
Let $X$ and $Y$ be projective varieties over $\Spec(k)$, $x\in X(\Qbar)$, and $L$ a line bundle on $X$.
\begin{enumerate}
\item For any positive integer $m$, $\alpha_{x,X}(m\cdot L)=m\cdot
\alpha_{x,X}(L)$.  This allows an extension of the definition of $\alpha_{x,X}(L)$ to $\QQ$-divisors $L$.
\item  $\alpha_{x}$ is a concave function of $L$:
for any positive rational numbers $a$ and $b$, and any $\QQ$-divisors $L_1$ and $L_2$
(with the exception of the case that $\{\alpha_{x}(L_1),\alpha_{x}(L_2)\} = \{-\infty,\infty\}$) we have
$$\alpha_x(a L_1+ b L_2) \geq a\alpha_x(L_1)+b\alpha_x(L_2).$$
\item If $Z$ is a subvariety of $X$ then for any point $z\in Z(\Qbar)$
we have $\alpha_{z,Z}(L|_Z)\geq \alpha_{z,X}(L)$.
\item If $x\in X(k)$ and $L$ is very ample then $\alpha_{x,X}(L) \geq  1$; if $x\in X(\Qbar)$ and
$L$ is ample then $\alpha_{x}(L)>0$.
\item Let $L_X$ and $L_Y$ be line bundles on $X$ and $Y$ which are asymptotically base point free,
and $x\in X(\Qbar)$, $y\in Y(\Qbar)$.  If neither $x$ nor $y$ are defined over $k$, then
\[\alpha_{x\times y,X\times Y}(L_X\squareplus L_Y) \geq \alpha_x(L_X)+\alpha_y(L_Y)\]
If $x$ is defined over $k$ but $y$ is not, then
\[\alpha_{x\times y,X\times Y}(L_X\squareplus L_Y) = \alpha_{y,Y}(L_Y).\]
If $x$ and $y$ are both defined over $k$, then
\[\alpha_{x\times y,X\times Y}(L_X\squareplus L_Y) = \min\{\alpha_{x,X}(L_X), \alpha_{y,Y}(L_Y)\}.\]
\item Suppose that $X$ is reducible over $k$ and let $X_1$,\ldots, $X_r$ be the irreducible components (over $k$)
containing $x$.  Then $\alpha_{x,X}(L) = \min(\alpha_{x,X_1}(L|_{X_1}),\ldots, \alpha_{x,X_r}(L|_{X_r}))$.
\end{enumerate}
\end{proposition}

\noindent
\bpf Since (up to equivalence) $H_{mL} = H_{L}^{m}$, part (a) follows immediately.

To simplify notation in part (b) set $\alpha_1=\alpha_{x}(L_1)$ and $\alpha_2=\alpha_{x}(L_2)$.  We will first prove (b) under the assumption that both $\alpha_1$ and $\alpha_2$ are finite.  We further note that in light of part (a), we may assume that $a+b=1$.

Suppose that there is a sequence $\xseq$ with $\alpha_{x}(\xseq,aL_1+bL_2)< a\alpha_1+b\alpha_2$.  
Fix $\delta>0$ small enough so that $a\alpha_1+b\alpha_2-\delta > \alpha_{x}(\xseq,aL_1+bL_2)$.  Then

\begin{equation}\label{eqn:convex-delta}
\rule{1.50cm}{0cm}
\sepv(x,x_i)^{a\alpha_1+b\alpha_2-\delta}H_{aL_1+bL_2}(x_i)
=
\left({\sepv(x,x_i)^{\alpha_1-\delta}H_{L_1}(x_i)\rule{0cm}{0.4cm}}\right)^{a}
\left({\sepv(x,x_i)^{\alpha_2-\delta}H_{L_2}(x_i)\rule{0cm}{0.4cm}}\right)^{b}.
\end{equation}

\noindent
By definition of $\alpha_1$ the term $\sepv(x,x_i)^{\alpha_1-\delta}H_{L_1}(x_i)$ is unbounded.  Hence, by passing to 
a subsequence of the $x_i$ (which can only lower the value of $\alpha_{x}(\xseq,L)$), we can assume that 
$\sepv(x,x_i)^{\alpha_1-\delta}H_{L_1}(x_i)\to\infty$ as $i\to\infty$.  
By definition of $\alpha_2$ the term $\sepv(x,x_i)^{\alpha_2-\delta}H_{L_2}(x_i)$ is also unbounded, and hence the 
left side of \eqref{eqn:convex-delta} is unbounded as well.  This implies that 
$\alpha_{x}(\xseq,aL_1+bL_2)\geq \alpha_1+\alpha_2 - \delta$, in contradiction to the way $\delta$ was chosen.  
Hence, for all sequences $\xseq$ of $k$-points we have $\alpha_{x}(\xseq,aL_1+bL_2) \geq a\alpha_1+b\alpha_2$.  
Taking the infimum over all sequences we conclude that $\alpha_{x}(aL_1+bL_2)\geq a\alpha_1+b\alpha_2$, which is the 
inequality in (b).

When one or both of $\alpha_1$ and $\alpha_2$ are infinite, with the exception of the case $\{\alpha_1,\alpha_2\}=\{\infty,-\infty\}$ either the resulting statement is obvious (for instance if both $\alpha_1=\alpha_2=-\infty$ then the 
bound is $\alpha_{x}(\xseq,aL_1+bL_2)\geq -\infty$ which is automatically true) or a minor variation of the argument 
above works.  In the case that $\{\alpha_1,\alpha_2\}=\{\infty,-\infty\}$ then it is not possible to deduce an upper 
bound for $\alpha_{x}(\xseq,aL_1+bL_2)$ from the data given (and also not clear what the purported upper bound of the 
form ``$\infty-\infty$'' is supposed to mean).

Part (c) is simple: We may assume that the distance function on $Z$ is the restriction of the distance
function on $X$ and that the height function on $Z$ is the restriction of $H_{L}$ to $Z(\Qbar)$.
Then for any sequence $\zseq$ of points of $Z(k)$ converging to $z$ we have
$\alpha_{z,Z}(\zseq, L|_{Z})=\alpha_{z,X}(\zseq, L)$.
The statement in (c) then follows from the observation that the set of $k$-points of
$Z$ is a subset of the set of $k$-points of $X$, and so the infimum used to define $\alpha_{z,Z}(L|_{Z})$ is
over a subset of the sequences used to define $\alpha_{z,X}(L)$.

For (d), if $L$ is very ample then $L$ induces an embedding $X\hookrightarrow \PP^n$ in some projective space.
If $x\in X(k)$ then by part (c) and Lemma~\ref{lem:projective} we conclude that
$\alpha_{x,X}(L)\geq \alpha_{x,\PP^n}(\Osh_{\PP^n}(1))=1$.
If $L$ is ample then some multiple $mL$ is very ample, and so if $x\in X(k)$ then
$\alpha_{x}(L)\geq\frac{1}{m}$ by the first part of this statement and homogeneity.  Finally, if $x\in X(\Qbar)$ let $K$
be the field of definition of $x$.   We have just established that $\alpha_{x}(L)_{K}>0$,
hence by 
Proposition \ref{prop:change-of-field} we have
$\alpha_{x}(L)=\alpha_{x}(L)_{k} \geq \frac{m_v}{d}\alpha_{x}(L)_{K}> 0$.

To prove claim (e), notice that the height function with respect to $L_X\squareplus L_Y$ is the product of the height functions of $L_X$ and $L_Y$.  Since $\sepv((x_1,y_1),(x_2,y_2))=\sepv(x_1,x_2)+\sepv(y_1,y_2)$ is a distance function on $X\times Y$, we may take that as our distance function for $X\times Y$.

Let $\{(x_i,y_i)\}$ be a sequence of $k$-points approximating $(x,y)$.  If $\{x_i\}$ and $\{y_i\}$ are both eventually contained in $X-\{x\}$ and $Y-\{y\}$, respectively, then by the definition of $\alpha_x$ and $\alpha_y$, we must have
\[\alpha_{x\times y,X\times Y}(\{(x_i,y_i)\},L_X\squareplus L_Y)\geq\alpha_x(L_X)+\alpha_y(L_Y)\]
as desired.

If $\{x_i\}$ is eventually equal to $x$, we get
\[\alpha_{x\times y,X\times Y}(L_X\squareplus L_Y) =\alpha_{y,Y}(L_Y).\]
Similarly, if $\{y_i\}$ is eventually equal to $y$, we get
\[\alpha_{x\times y,X\times Y}(L_X\squareplus L_Y) =\alpha_{x,X}(L_X).\]

To finish the proof, it remains only to note that $\{x_i\}$ and $\{y_i\}$ are sequences of $k$-rational points, so that $\{x_i\}$ can only be eventually the constant sequence $\{x\}$ if $x$ is $k$-rational, and similarly for $y$.

Finally, statement (f) follows by the pigeonhole principle: if $\xseq$ is a sequence approximating $x$, then
infinitely many $x_i$ must lie on some component $X_j$, and by passing to a subsequence we may assume that all $x_i$
lie on $X_j$.  Thus $\alpha_{x,X}(L)$ is no more than the minimum in part (f).  The opposite inequality follows
from part (c).
\qed

\vspace{.1in}

\noindent
{\bf Remarks on extreme cases.} (a) If $\alpha_{x}(L)=\infty$ for one line bundle then $\alpha_{x}(A)=\infty$ for all ample line bundles $A$.  Indeed, for any sequence $\xseq$, if $\alpha_{x}(\xseq,L)=\infty$ then $\alpha_{x}(\xseq,A)=\infty$ for all ample line bundles $A$.  This follows immediately from the fact that there is some positive integer $n$ such that $nA-L$ is effective, giving $H_{nA}(x_i)\geq H_L(x_i)+O(1)$ for all $i$ and so $\alpha_{x}(\xseq,nA)\geq\alpha_{x}(\xseq,L)=\infty$.  Thus, by Proposition~\ref{prop:alpha}, part (a), $\alpha_{x}(\xseq,A)=\infty$.  

\noindent
(b) Assume that there is no nef line bundle $L$ so that $\alpha_x(L)=\infty$.
The concavity condition shows that $\alpha_x$ is a continuous function on
the ample cone.

\noindent
(c) If $X$ is smooth and $L$ is ample, then any sequence $\xseq$ such that $\alpha_{x}(\xseq,L)$ is finite must eventually lie in a fibre of the Albanese map $\pi\colon X\to A$.  This follows from the fact that $\alpha$ is infinite on Abelian varieties (see Example (c) in the introduction).  More precisely, let $D$ be an ample divisor on the Albanese variety $A$. Then there is some positive integer $n$ such that $L_n=nL+\pi^*D$ is ample.  If $\alpha_{x}(\xseq,\pi^*D)$ is finite, then clearly $\xseq$ is eventually contained in a fibre of the Albanese map, since $\pi$ does not increase distances by more than a bounded multiple.  Since $L_n-\pi^*D$ is effective, this means that $\alpha_{x}(\xseq,L_n)$ is also infinite unless $\xseq$ is eventually contained in a fibre of the Albanese map.  By Remark (a), this means that for any ample divisor $L$, $\alpha_{x}(\xseq, L)=\infty$, unless $\xseq$ is eventually contained in a fibre of the Albanese map.

\vspace{.1in}

\begin{lemma} \label{lem:Roth-on-P1}
Let $d$ be a positive integer, $L=\Osh_{\PP^1}(d)$, and $x\in \PP^1(\Qbar)$.  Then
\[
\alpha_{x}(L) =
\begin{cases}
\infty & \text{if $\kappa(x)\not\subseteq k_v$} \\
 d     & \text{if $\kappa(x) = k$} \\
 \frac{d}{2} & \text{otherwise.}
\end{cases}
\]
\end{lemma}

\bpf
If $\kappa(x)\not\subset k_v$ then there is no sequence of $k$-points converging (with respect to $\sepv(\cdot,\cdot)$)
to $x$ (see the Remark on page \pageref{rem:same-local-degree}), and hence $\alpha_{x}(L)=\infty$.
If $x\in \PP^1(k)$ then this is Lemma \ref{lem:projective} and Proposition \ref{prop:alpha}(a).
If $\kappa(x)\subseteq k_v$ but $\kappa(x)\neq k$ then $\alpha_{x}(\Osh_{\PP^1}(1))\geq \frac{1}{2}$ by Roth's
theorem for $\PP^1$, while $\alpha_{x}(\Osh_{\PP^1}(1))\leq \frac{1}{2}$ by a Dirichlet-type argument.  (This follows, for example, from Theorem~\ref{thm:general-dirichlet}.)  Thus $\alpha_{x}(\Osh_{\PP^1}(1))=\frac{1}{2}$, and
so $\alpha_{x}(L) = \frac{d}{2}$ by Proposition \ref{prop:alpha}(a) again. \epf

\noindent
{\bf Remark:} \/ Note that the cases in Lemma~\ref{lem:Roth-on-P1} depend sensitively upon the choice of extension of
$v_0$ to $\Qbar$.  For example, if $\kappa(x)$ is not a Galois extension of $k$, then it is possible that for some
choices of $v$ on $\Qbar$, $k_v$ contains $\kappa(x)$, and for others it does not.
This leads to radically different values of $\alpha_{x}(L)$.

\medskip

\begin{theorem}\label{thm:curve}
Let $C$ be any singular $k$-rational curve and $\varphi\colon\PP^1\rightarrow C$ the normalization map.
Then for any ample line bundle $L$ on $C$, and any $x\in C(\Qbar)$ we have the equality:
\[\alpha_{x,C}(L)=\min_{q\in \varphi^{-1}(x)} d/r_{q} m_{q}\]
where $d=\deg(L)$, $m_{q}$ is the multiplicity of the branch of $C$ through $x$ corresponding to $q$, and
\[r_{q}=
\begin{cases}
0 & \text{if $\kappa(q)\not\subseteq k_v$} \\
1 &\text{if $\kappa(q)=k$} \\
2 &\text{otherwise.}
\end{cases}
\]
\end{theorem}

\np
Here we use $r_q=0$ as a shorthand for $d/r_{q}m_{q}=\infty$.

\noindent
{\it Proof:}
Given any sequence $\xseqtox$ then by passing to a subsequence we can assume that all $x_i$ are on a single
branch through $x$.  More precisely, we can assume that none of the $x_i$ are the finitely many points where
$\varphi$ is not an isomorphism, and that $\{\varphi^{-1}(x_i)\}$ converges (with respect to $\sepv(\cdot,\cdot)$)
to one of the points $q\in \varphi^{-1}(x)$.
Conversely, given a sequence $\qseq$ of points of $\PP^1(k)$ converging to some $q$,
then $\{\varphi(q_i)\}$ converges to $x$.
Thus it suffices to study only sequences of this kind to compute $\alpha_{x}(L)$.

Given a sequence $\qseq\to q$ we have $H_{\varphi^{*}L}(q_i)=H_{L}(\varphi(q_i))$ for all $i$.  Furthermore since
the branch corresponding to $q$ has multiplicity $m_q$, locally $\varphi$ is described by functions in the $m_q$-th
power of the maximal ideal of $q$, and thus
$\sepv(x,\varphi(q_i))$ is equivalent to $\sepv(q,q_i)^{m_q}$ as $i\to\infty$.  Therefore, as in
Proposition \ref{prop:change-of-field} we have $\alpha_{x}(\{\varphi(q_i)\},L) =
\frac{1}{m_q}\alpha_{q}(\qseq,\varphi^{*}L)$,
and the theorem then follows from Lemma \ref{lem:Roth-on-P1}. \epf

\vspace{.1in}

\noindent
{\bf Remark:} \/ This is similar to Theorem~2.8 of \cite{McK}, except
that it is actually correct.  (The conclusion of Theorem~2.8 of
\cite{McK} neglects the possibility that the $r_q$ defined in
Theorem~\ref{thm:curve} might not be one.)  Theorem~\ref{thm:curve} also uses
the definition of $\alpha$ from this paper, rather than that of
\cite{McK}, and generalises the results to points defined over $\overline{k}$.

\newpage
\noindent
{\bf Examples}\label{ex:alpha}

\begin{itemize}
\item[(a)] If $X=\PP^N$, $L=\Osh_{\PP^N}(d)$ for some $d>0$, then $\alpha_x(L)=d$ for all
  points $x$ in $\PP^N(k)$.  This follows from Lemma~\ref{lem:projective}
  and Proposition~\ref{prop:alpha}(a).
\item[(b)] If $X=\PP^1\times\PP^1$, $L=\Osh_{\PP^1\times\PP^1}(a,b)$, with
  $a,b\geq 0$ then $\alpha_x=\min(a,b)$ for all $x\in X(k)$. This
  follows immediately from Proposition~\ref{prop:alpha}(e).
\item[(c)] Similarly if $X=\PP^{N_1}\times\cdots\times\PP^{N_r}$,
$L=\Osh_X(d_1,\ldots,d_r)$ with $d_i\geq 0$ then
$\alpha_x(L)=\min(d_1,\ldots,d_r)$.
\item[(d)] Taking $X=\PP^1\times\PP^1$, $L_1=L_2=\Osh_X(2,1)$, $L_3=\Osh_X(1,2)$
in example (b) gives $\alpha_x(L_i)=1$ for $i=1,2,3$, but
$\alpha_x(L_1+L_2)=2$ and $\alpha_x(L_1+L_3)=3$.
\end{itemize}

Part (d) shows that there can be no formula for determining
$\alpha_x(L_i+L_j)$ in terms of $\alpha_x(L_i)$ and $\alpha_x(L_j)$
alone, and that Proposition~\ref{prop:alpha}(b) is the best possible
general relation of this type.

\medskip
The following lemma, which we will use several times in the paper, allows us to
reduce to the case of geometrically irreducible varieties when studying $\alpha$.

\begin{lemma}\label{lem:irreducible-Z}
Let $Z$ be a variety over $\Spec(k)$, and set $Y$ to be the Zariski closure
of the points of $Z(k)$.  Then each irreducible component of $Y$ is
geometrically irreducible and for any line bundle $L$ on $Z$ and any $x\in Y(\Qbar)$ we have
$\alpha_{x,Z}(L)  = \alpha_{x,Y}(L|_{Y})$.
\end{lemma}

\bpf
Let $Y_1$, \ldots, $Y_r$  be the irreducible
components of $\Ybar:=Y\times_k\Qbar$;  we will show that each $Y_i$ is actually defined over $k$.
Let $Y_i$ be one such component.  Since $Y$ is a variety over $\Spec(k)$, all $\Gal(\Qbar/k)$ conjugates of $Y_i$ are
also components of $\Ybar$.
Let $I\subseteq\{1,\ldots, r\}$ be the subset of indices such that each $Y_j$, $j\in I$, is
a Galois conjugate of $Y_i$, and set $I'=\{1,\ldots, r\}\setminus I$.
Any point $y\in Z(k)$ contained in $Y_i$ is also contained in $Y_j$ for $j\in I$.
Therefore all points of $Z(k)$ are contained in $Y':=(\bigcap_{j\in I} Y_j )\bigcup (\bigcup_{j'\in I'}Y_{j'})$.
By construction $Y'$ is closed and defined over $k$.
If $I\neq \{i\}$ then $Y'$ is a proper subset of $\Ybar$.  This contradicts the construction
of $Y$ as the Zariski closure of $Z(k)$.  Thus $I=\{i\}$ and so $Y_i$ is defined over $k$.
Finally since $Y(k)=Z(k)$, it is clear that
$\alpha_{x,Z}(L)  = \alpha_{x,Y}(L|_{Y})$ for all line bundles $L$ and $x\in Y(\Qbar)$.
\epf

\section{Seshadri constants}\label{sec:seshadri}

In this section, we review some basic properties of Seshadri constants,  first introduced and studied
in \cite{Dem}.  Many foundational results on Seshadri constants are given in \cite[chap.\ 5]{PAG}.
The Seshadri constant is purely geometric in the sense that it only depends on the base change of the
variety to the algebraic closure.

\begin{definition}\label{def:seshadri}
Let $X$ be a  projective variety over $\Spec(k)$, $x$ a point of $X(\Qbar)$, and $L$ a
nef line bundle on $X$.  The {\em Seshadri constant},
$\ep_{x,X}(L)$, is defined to be

\[\ep_{x,X}(L) :=
\sup\left\{{\gamma\geq 0 \mid \pi^{*}L -\gamma E\,\,\,
\mbox{is nef}\,}\right\}\]
where  $\pi:\Xtil\longrightarrow X_{\Qbar}$ is the blowup of $X_{\Qbar}:=X\times_{k} \Qbar$ at $x$ with
exceptional divisor $E$.  Here, by abuse of notation, we also use $L$ for the base change of $L$ to $X_{\Qbar}$.
\end{definition}

The Seshadri constant is defined on the level of $\QQ$- or
$\RR$-divisors, and in the above definition $\gamma\geq 0$ is an
element of $\QQ$.  If $\gamma$ is allowed to be a real number, then
the $\sup$ in the definition can be replaced by a $\max$.

The idea behind the Seshadri constant is that it measures the local
positivity of $L$ at $x$.  From the definition, the Seshadri constant
only depends on the numerical equivalence class of $L$.
We will often just use $\ep_x(L)$ or $\ep_x$ for
$\ep_{x,X}(L)$ if $X$ or $L$ are clear from the context.

Since the Seshadri constant only depends on $X_{\Qbar}$,
for the rest of this section we assume our varieties are defined over a fixed algebraically closed field.
From Definition \ref{def:seshadri}, all of the properties of the Seshadri constant established below will hold for
varieties over $\Spec(k)$.

Another characterization of the Seshadri constant is given by the following.

\begin{proposition}\label{prop:seshcurvemult}
Let $X$ be a projective variety, $x\in X$, and $L$ a nef line bundle on $X$, then
\[
\ep_{x,X}(L) =
\inf_{x\in C\subseteq\, X}
\left\{
\frac{(L\cdot C)}{\mult_x(C)}
\right\}
\]
where the infimum is taken over all reduced irreducible curves $C$ passing
through $x$.
\end{proposition}

This alternate description of the Seshadri constant follows
immediately from the definition that a bundle $L'$ on
a variety $\Xtil$ is nef if and only if $L'\cdot C'\geq 0$ for all
reduced irreducible curves $C'$ in $\Xtil$, and the straightforward
observation that if $C'$ is the proper transform of $C$ in the blowup,
then $E\cdot C'=\mult_{x}(C)$, and $(\pi^{*}L)\cdot C'=L\cdot C$.

\noindent
{\bf Basic properties of $\ep$.} We start by computing $\ep$ when $X=\PP^n$.

\begin{lemma}\label{lem:ep-projective}
For any point $x\in \PP^n$, $\ep_{x}(\Osh_{\PP^n}(1))=1$.
\end{lemma}

\bpf
Let $\pi\colon\Ptil^n\longrightarrow \PP^n$ be the blowup of $\PP^n$ at $x$.  For any $\gamma>0$ set
$L_{\gamma}:=\pi^{*}(\Osh_{\PP^n}(1))-\gamma E$.  Then $L_{1}$
is base point free and defines the projection morphism $\Ptil^n\longrightarrow \PP^{n-1}$ with fibres the proper
transforms of lines in $\PP^n$ passing through $x$.  Thus $L_1$ is nef on $\Ptil^n$.  For any such fibre the
degree of $L_{\gamma}$ on the fibre is $1-\gamma$, hence $L_1$ is the boundary of the nef cone, and
$\ep_x(\Osh_{\PP^n}(1))=1$.  \epf

\vspace{.1in}

Note that Lemma~\ref{lem:ep-projective} shows that if $x\in\PP^n(k)$, then $\ep_x=\alpha_x$.  The following proposition extends the list of similarities between $\ep$ and $\alpha$ much further.

\begin{proposition}\label{prop:ex}
Let $X$ be a projective variety, $x\in X(\Qbar)$, and $L$ a nef line bundle on $X$.
\begin{enumerate}
\item For any positive integer $m$, $\ep_{x,X}(m\cdot L)=m\cdot \ep_{x,X}(L)$.  This allows an extension of the definition of $\ep_{x,X}(L)$ to $\QQ$-divisors $L$.
\item  $\ep_{x}$ is a concave function of $L$:
for any positive rational numbers $a$ and $b$, and any nef $\QQ$-divisors $L_1$ and $L_2$
$$\ep_x(a L_1+ b L_2) \geq a\ep_x(L_1)+b\ep_x(L_2).$$
\item If $Z$ is a subvariety of $X$ then for any point $z\in Z$
we have $\ep_{z,Z}(L|_Z)\geq \ep_{z,X}(L)$.
\item If $L$ is very ample then $\ep_{x}(L)\geq 1$, if $L$ is ample then
  $\ep_{x,X}(L)> 0$.
\item If $x$ and $y$ are points of varieties $X$ and $Y$, with nef
line bundles $L_X$ and $L_Y$ then
\[\ep_{x\times y,X\times Y}(L_X\squareplus L_Y)
= \min(\ep_{x,X}(L_X), \ep_{y,Y}(L_Y)).\]
\item Suppose that $X$ is reducible and let $X_1$,\ldots, $X_r$ be the irreducible components containing $x$.
Then $\ep_{x,X}(L) = \min(\ep_{x,X_1}(L|_{X_1}),\ldots, \ep_{x,X_r}(L|_{X_r}))$.
\end{enumerate}
\end{proposition}

\noindent
\bpf
The definition implies (a) immediately.
Part (b) is also clear from the definition: if
$\pi^{*}L_1-\ep_1\cdot E$ and $\pi^{*}L_2-\ep_2\cdot E$ are
nef on $\Xtil$, then so is
$\pi^{*}(aL_1+bL_2)-(a\ep_1+b\ep_2)\cdot E = a(\pi^{*}(L_1)-\ep_1\cdot E)+b(\pi^{*}(L_2)-\ep_2\cdot E)$.

To prove (c), it is enough to remark that the proper transform of $Z$ in the
blow up $\Xtil$ of $X$ at $z$ is the blow up $\Ztil$ of $Z$ at $z$,
and that the restriction of a nef bundle on $\Xtil$ will be a nef
bundle on $\Ztil$.

For (d), if $L$ is very ample then $L$ induces an embedding $X\hookrightarrow \PP^n$ in some projective space.
By part (c) and Lemma~\ref{lem:ep-projective} we conclude that
$\ep_{x,X}(L)\geq \ep_{x,\PP^n}(\Osh_{\PP^n}(1))=1$.
If $L$ is ample then some positive multiple $mL$ is very ample and so $\alpha_{x}(L)\geq\frac{1}{m}$ by the
first part of this statement and homogeneity.

The proper transforms of $X\times y$ and $x\times Y$ in the blow-up of
$X\times Y$ at $x\times y$ are the blowups $\Xtil$ and $\Ytil$ of $X$
at $x$ and $Y$ at $y$.  This and the observation that the restriction
of a nef bundle must be nef give
\[\ep_{x\times y,X\times Y}(L_X\squareplus L_Y) \leq
\min(\ep_{x,X}(L_X), \ep_{y,Y}(L_Y)).\]

To prove the other direction, we will use the description of
$\ep_{x\times y}$ from Proposition~\ref{prop:seshcurvemult}.
Let $\pi_X$ and $\pi_Y$ be the projections from $X\times Y$ to $X$ and
$Y$ and let $C$ be any irreducible curve in $X\times Y$ passing
through $x\times y$.

Let $\pi_X(C)$ be the reduced image of $C$.  Suppose that $C$ is not
contained in a fibre of $\pi_X$.  Then $\pi_X(C)$ is not equal to a
point, and if $d$ is the generic degree of the map $C\longrightarrow
\pi_X(C)$ we have $\pi^{*}L_X\cdot C = d(L_X\cdot\pi_X(C))$, and
$\mult_{x\times y}(C) \leq d \cdot \mult_x(\pi_X(C))$.

Since $\ep_x$ is the Seshadri constant for $L_X$ at $x$, we have

$$
\ep_x \leq \frac{L_X\cdot \pi_X(C)}{\mult_x(\pi_X(C))}
\leq\frac{d(L_X \cdot \pi_X(C))}{\mult_{x\times y}(C)} =
\frac{\pi^{*}_X L_X \cdot C}{\mult_{x\times y}(C)} \leq
\frac{(\pi^{*}_X L_X+\pi_Y^{*}L_Y)\cdot C}{\mult_{x\times y}(C)}
$$

\noindent
where the first inequality follows from
Proposition~\ref{prop:seshcurvemult} applied to $\ep_x$, the
second from the inequality on the multiplicities, and
the third from the fact that $\pi^{*}L_Y$ is nef.

Similarly, if $C$ is not contained in a fibre of $\pi_Y$ we have the
corresponding inequality with $\ep_y$ in place of $\ep_x$.
Since for any given curve $C$ one of these must be true we have

\[\min(\ep_x,\ep_y)\leq
\inf_{x\times y\in C\subseteq\, X\times Y}
\left\{
\frac{(\pi_X^{*}L_X+\pi_Y^{*}L_Y)\cdot C}{\mult_{x\times y}(C)}
\right\}
\stackrel{\scriptsize\ref{prop:seshcurvemult}}{=}
\ep_{x\times y}
\]
finishing the proof of (e).

For part (f) we use the fact that a line bundle is ample if and only if it is ample restricted to each component,
and that the blow up of each $X_i$ at $x$ is a component of $\Xtil$.
\epf

\vspace{.1in}

{\bf Examples}\label{examples:epsilon}

\begin{itemize}
\item[(a)] If $X=\PP^n$, $L=\Osh_{\PP^n}(d)$ then
$\ep_x(L)=d$ for all points $x$ in $\PP^n$.  This follows from
the computation for $\PP^n$ and $\Osh_{\PP^n}(1)$ in Lemma \ref{lem:ep-projective} along with
Proposition \ref{prop:ex}(a).
\item[(b)] If $X=\PP^1\times\PP^1$, $L=\Osh_{\PP^1\times\PP^1}(a,b)$,
with $a,b\geq 0$ then $\ep_x=\min(a,b)$ for all $x\in X$.  This follows immediately from
Proposition \ref{prop:ex}(e) and part (a) of the examples, but we can also prove this as follows.  Let $\Xtil$ be the blow up of $X=\PP^1\times\PP^1$
at a point $x$, $E$ the exceptional divisor and $F_1$ and $F_2$ the
pullback of the class of fibres from $X$.  The effective cone of
$\Xtil$ is generated by $F_1-E$, $F_2-E$, and $E$.  Dually, the nef
cone of $\Xtil$ is generated by $F_1$, $F_2$ and $F_1+F_2-E$.

Therefore for $aF_1+bF_2-\gamma E$ to be in the nef cone, the
condition is exactly that $\gamma\leq\min(a,b)$, i.e.,
$\ep_x(aF_1+bF_2)=\min(a,b)$.

\item[(c)] Similarly if $X=\PP^{n_1}\times\cdots\times\PP^{n_r}$,
$L=\Osh(d_1,\ldots,d_r)$ with $d_i\geq 0$, for
$i=1,\ldots, r$ then $\ep_x(L)=\min(d_1,\ldots,d_r)$.
\end{itemize}

As evidenced by our parallel statements in
Proposition~\ref{prop:alpha} and Proposition~\ref{prop:ex}
(and Lemmas \ref{lem:projective} and \ref{lem:ep-projective}, and the examples)
there is a great deal of formal similarity between
$\alpha_x$ and $\ep_x$.
See the discussion below on the Arakelov point of view for some heuristic reasons for this similarity.

For future reference we record the exact conditions on a curve $C$ and point $x\in C(\Qbar)$ so that
$\alpha_{x}(L)=\frac{1}{2}\ep_{x}(L)$.

\begin{lemma}\label{lem:roth-bound-for-curves}
Let $C$ be an irreducible curve over $\Spec(k)$, $x\in C(\Qbar)$ and $L$ any ample line bundle on $C$.
Then $\alpha_{x}(L) = \frac{1}{2}\ep_{x}(L)$ if and only if $C$ is a $k$-rational curve, $C$ is unibranch at $x$,
$\kappa(x)\neq k$, and $\kappa(x)\subseteq k_v$.
\end{lemma}

\bpf
Since $\ep_{x}(L)$ is always finite, the equality implies that $\alpha_{x}(L)$ is finite, and hence that $C$
is a $k$-rational curve.   Let $\varphi\colon \PP^1\longrightarrow C$ be the normalization map, and for any
$q\in \varphi^{-1}(x)$ let $m_q$ be the multiplicity at $x$ of the branch corresponding to $q$, and define $r_q$
as in Theorem \ref{thm:curve}.  By that theorem we have
$\alpha_{x}(L) = \min_{q\in \varphi^{-1}(q)}\{\frac{d}{r_q m_q}\}$ where $d=\deg(L)$.  By the definition of the
Seshadri constant we have $\ep_{x}(L) = \frac{d}{\mult_{x}C}=\frac{d}{\sum_{q\in \varphi^{-1}(x)} m_q}$.
Thus the equality
$\alpha_{x}(L) = \frac{1}{2} \ep_{x}(L)$ amounts to the equality
$$\max_{q\in \varphi^{-1}(x)} \{r_qm_q \} = 2\sum_{q\in \varphi^{-1}(x)} m_q.$$
Since $r_q\in\{0,1,2\}$ for each $q$, the only possible way to have equality above is if $\varphi^{-1}(x)$ consists
of a single point $q$ with $r_q=2$.  Given the definition of $r_q$ in Theorem \ref{thm:curve} this proves
the lemma.  \epf

\noindent
{\bf Arakelov point of view.}
For the rest of this section we discusses some parallels between $\alpha$ and $\ep$ from the point of
of Arakelov theory.   Although it does not explain those parallels, we think that this heuristic discussion
is useful.

Let $X$ be a projective variety over $\Spec(k)$ and $x$ a point of $X(k)$.
Let $\Xtil$ be the blow up of $X$ at $x$ with exceptional divisor $E$.
By Kleiman's characterization of the ample cone,
the definition of the Seshadri constant $\ep_{x}=\ep_{x}(L)$ is that for any $0<\gamma< \ep_x$ the set
$$
\left\{{ B\subseteq \Xtil_{\Qbar} \st \mbox{$B$ an irreducible curve},\,\, (L-\gamma E)\cdot B < 0 }\right\}
$$
is empty, and $\ep_x$ is the largest number with this property.

Let $\imXt$ be a proper integral model for $\Xtil$ over $\Spec(\Osh_{k})$.
We consider each point $y\in X(k)$, $y\neq x$, to also be a point of $\Xtil(k)$, and hence each $y$ gives rise to
a section $\sigma_y$ of $\imXt$ over $\Spec(\Osh_{k})$.  Choose suitable metrizations of $L$ and $E$ on the archimedean
places of $k$.  By the Arakelov construction of the intersection product on $\imXt$,
for any $\gamma>0$ we have $$h_{L-\gamma E}(y) = (L-\gamma E)\cdot \sigma_y.$$
(For details on Arakelov intersection theory, see for example \cite[\S III.2]{So}.)

Choose an embedding $\varphi\colon X\hookrightarrow \PP^r$ so that $x\mapsto [1\colon 0\colon\cdots\colon 0]$.
Let $Z_0$,\ldots, $Z_r$ be the coordinates on $\PP^r$ and define functions $u_i$, $i=1,\ldots, r$ on the
open subset $U$ where $Z_0\neq 0$ by $u_i=Z_i/Z_0$.
For each place $w$ of $k$, define a function $e_w\colon X(k)\rightarrow \RR_{\geq 0}$ by

$$e_w(y) = \left\{{\begin{array}{cl}
1 & \mbox{if $y\not\in U(k)$,} \\
\min\left({1, \max(\nrm{u_1(y)}_{w},\ldots, \nrm{u_r(y)}_{w})}\right) & \mbox{if $y\in U(k)$.} \\
\end{array}}\right.
$$

A short local calculation (see \cite[Lemma 3.1]{McKR}) shows that $-h_{E}(y) = \sum_{w} \log(e_{w}(y))$. By
Proposition \ref{prop:equiv-dist} and Lemma \ref{lem:local-distance-Zariski-open},
$e_w(\cdot)$ is equivalent to $\sepw(x,\cdot)$ on $U(k)$ for each place $w$.
Thus, up to a bounded constant which we ignore, we have

\begin{equation}\label{eqn:Arakelov-intersection}
\rule{0.5cm}{0cm} (L-\gamma E)\cdot\sigma_y =
h_{L-\gamma E}(y) =
h_L(y) + \gamma \log(\sepv(x,y)) + \gamma \left({\sum_{w\neq v} \log(\sepw(x,y)})\right) .
\end{equation}

By Proposition~\ref{prop:equiv-alpha} for any ample line bundle $L$ an equivalent description
of $\alpha_{x}=\alpha_{x}(L)$ is that for any $\gamma< \alpha_{x}$ the set

\begin{equation}\label{eq:finite-cond-exp}
\left\{{ y\in X(k) \st \sepv(x,y)^{\gamma}H_{L}(y) < 1}\right\}
\end{equation}

\noindent
is finite, and $\alpha_x$ is the largest number with this property.\footnote{The extra quantifier ``$C$'' in
Definition \ref{def:B-Northcott} can be absorbed by the condition that the finiteness is supposed to hold for
all $\gamma<\alpha_x$.  The purpose of this quantifier in Definition \ref{def:B-Northcott}
is to simplify arguments.}
Taking $\log$, the finiteness of \eqref{eq:finite-cond-exp} is equivalent to the finiteness of

\begin{equation}\label{eqn:finite-cond}
\left\{{ y\in X(k) \st h_{L}(y)+\gamma \log(\sepv(x,y))  < 0}\right\}
\end{equation}

\noindent
where $h_L$ is the logarithmic height.
Since the logarithmic height is only defined up to a bounded constant, ``finitely many'' is the best substitute
for ``none'', and this makes the definition of $\alpha_x$ look very much like the definition of $\epsilon_x$.
Equation \eqref{eqn:Arakelov-intersection} suggests an even closer parallel: that we interpret
$h_{L}(y) + \gamma\log(\sepv(x,y))$ as the intersection ``$(L-\gamma E_v)\cdot \sigma_y$'', where $E_v$ is meant
to indicate that we only count the local contribution of $E$ at the place $v$ when computing the intersection with
$\sigma_y$.

From this point of view the statements in Propositions \ref{prop:alpha} and \ref{prop:ex} have almost identical
proofs.  For instance, here are the arguments for the superadditivity of $\alpha_x$ and $\ep_x$ (part ({\em b}) of the
respective propositions).  For $\ep_x$ the argument is:  if there are no curves $B$ such that
$(L_1-\gamma_1 E)\cdot B<0$ and none such that $(L_2-\gamma_2 E)\cdot B <0$ then there are no curves $B$ such that
$\left((L_1-\gamma_1E)+(L_2-\gamma_2E)\right)\cdot B <0$.
For $\alpha_x$ the argument is: if there are only finitely many $y\in X(k)$ such that
$(L_1-\gamma_1 E_v)\cdot \sigma_y<0$, and only finitely many such that
$(L_2-\gamma_2 E_v)\cdot \sigma_y<0$, then there are only finitely many $y\in X(k)$ such that
$\left((L_1-\gamma_1E_v)+(L_2-\gamma_2E_v)\right)\cdot \sigma_y <0$.

We finish the discussion with another connection between the two invariants.  Since $\sepw(x,y)\in (0,1]$, we
have $\log(\sepw(x,y))\leq 0$ and thus by \eqref{eqn:Arakelov-intersection} if $\gamma >0$
the set \eqref{eqn:finite-cond} is contained in the set

\begin{equation}\label{eqn:finite-cond-all-places}
\left\{{ y\in X(k) \st (L-\gamma E)\cdot\sigma_y < 0 }\right\}.
\end{equation}

\noindent
If $0<\gamma < \ep_x(L)$ then the line bundle $L-\gamma E$ is ample, and therefore satisfies the Northcott property.
For such $\gamma$, the set \eqref{eqn:finite-cond-all-places} and hence the subset
\eqref{eqn:finite-cond} is finite.   This proves the inequality $\alpha_x(L)\geq \ep_x(L)$, a stronger inequality
than $\alpha_x(L)\geq \frac{1}{2}\ep_x(L)$, one of the main goals of this paper.
However, in the discussion above we have assumed that $x\in X(k)$.  For an arbitrary point $x\in X(\Qbar)$ with field of
definition $K$ this argument, along with the change of field formula of Proposition \ref{prop:change-of-field},
yields the inequality $\alpha_x(L)\geq \frac{1}{[K:k]}\,\ep_x(L)$.
This is the Liouville theorem mentioned in the introduction.
(See also \cite[Theorem 3.3]{McKR} for an extension involving $\gamma$ in the big cone and the asymptotic base locus.)

\pdfstringdef{\pdfbeta}{\337}   

\section{\texorpdfstring{The constant $\Ar_x(L)$}{The constant \pdfbeta(L)}}
\label{sec:f-def}

The proof of the general version of Roth's theorem will involve a third invariant of a point and an ample line bundle.
In this section we define this invariant and establish some basic facts to be used in the proof. 
As with the Seshadri constant this invariant only depends on the base change of the variety to an algebraically
closed field.    We start by describing the invariant in this case, and then give the general definition below.

First suppose that $X$ is an irreducible $n$-dimensional variety defined over an algebraically closed field.
For any ample line bundle $L$ on $X$ and $x\in X$,
let $\pi\colon\Xtil\longrightarrow X$ be the blow up at $x$ with exceptional divisor $E$,
and for any $\gamma\in \RR_{\geq 0}$ set $L_{\gamma}:=\pi^{*}L-\gamma\,E$.

\noindent
Let $\NS(\Xtil)_{\RR}$ be the real N\'eron-Severi group of $\Xtil$ and
let $\Vol(\cdot)$ be the volume function on $\NS(\Xtil)_{\RR}$.
Recall that the {\em volume}, $\Vol(M)$, of a line bundle $M$ on an $n$-dimensional variety
measures the asymptotic growth of the global sections of $M$.
Specifically $\Vol(M)$ is the unique real number so that
$h^0(mM) = \frac{\Vol(M)}{n!}m^n + O(m^{n-1})$ for $m\gg0$.
From the definition it follows that $\Vol(mM)=m^n\Vol(M)$ for
$m\geq 0$, so that $\Vol(\cdot)$ may be extended to $\QQ$-bundles.    By \cite[Corollary 2.2.45]{PAG} $\Vol(\cdot)$
depends only on the numerical class of $M$ and extends uniquely to a continuous function on the
real N\'eron-Severi group.
A line bundle $M$ is called {\em big} if $\Vol(M)\neq 0$.

Let $\geff=\geffx(L)=\sup \{\gamma\in \RR_{\geq 0} \st \mbox{$L_{\gamma}$ is effective}\}$.
We note that $\geff$ is always finite: if $A$ is an ample bundle on $\Xtil$ and $L_{\gamma}$ effective then
$L_{\gamma}\cdot A^{n-1}=(\varphi^{*}L)\cdot A^{n-1}-\gamma\,(E\cdot A^{n-1})>0$ and hence
$\geff \leq \frac{(\varphi^{*}L)\cdot A^{n-1}}{E\cdot A^{n-1}}$.
The big cone is the interior of the effective cone, and therefore by definition of $\geff$ we have
$\Vol(L_{\gamma})>0$ for all $\gamma\in [0,\geff)$, $\Vol(L_{\gamma})=0$ for all $\gamma>\geff$, and so
also $\Vol(L_{\geff})=0$ by continuity of the volume function.
We define a decreasing function (the ``asymptotic relative volume function'')
$f:[0,\infty)\longrightarrow[0,1]$ by $$f(\gamma) = \frac{\Vol(L_{\gamma})}{\Vol(L)},$$
and note that $f$ is supported on $[0,\geff]$.  Finally, define\footnote{To the best of our knowledge, the
number $\Ar_x(L)$ was first defined by Per Salberger in unpublished work dating from 2006, where it was used
to improve results of R.\ Heath-Brown on uniform upper bounds for the number of rational points of bounded height.
Salberger also proved Corollary~\ref{cor:beta-bound-irred} as a key step in this work.}
$$\Ar_{x}(L) = \int_{0}^{\infty} f(\gamma) \,d\gamma =
\int_{0}^{\geff} f(\gamma) \,d\gamma$$
to be the area under $f$.

\noindent
{\bf Example.} \label{ex:Ar-Pn}
Let $X=\PP^n$, and $L=\Osh_{\PP^n}(1)$. We will check below that for any $x\in \PP^n$ we have $\geffx(L)=1$,
$f(\gamma)= 1-\gamma^n$ on $[0,1]$, and hence $\Ar_{x}(L)=\frac{n}{n+1}$.
This will turn out (via Theorem \ref{thm:RothIII} or \ref{thm:RothI}) to
explain the approximation constants of $\frac{1}{2}$ for $\PP^1$ (from the classical Roth's theorem) or
$\frac{n}{n+1}$ for $\PP^n$ (from the Schmidt subspace theorem).

We now verify the claims above.
As in the proof of Lemma \ref{lem:ep-projective},
let $\pi\colon\Ptil^n\longrightarrow \PP^n$ be the blowup of $\PP^n$ at $x$ and recall that $L_1$ is
base point free and defines a projection morphism  $\varphi\colon\Ptil^n\longrightarrow\PP^{n-1}$. 
The degree of $L_{\gamma}$ on the fibres of $\varphi$ is $1-\gamma$.  Hence for rational $\gamma>1$,
and $m>0$ such that $m\gamma$ is an integer,
any global section of $mL_{\gamma}$ vanishes on all fibres of $\varphi$, and is therefore zero.
Thus $\geff\leq 1$.
When $0<\gamma<1$, $L_{\gamma}$ is ample, and thus effective, and we conclude that $\geff=1$.

If $M$ is a big and nef line bundle then $\Vol(M)=c_1(M)^n$ (see \cite[p.\ 148]{PAG}), and therefore
$\Vol(L_{\gamma}) = c_1(L_{\gamma})^n = c_1(L)^n + (-\gamma)^n E^n = 1-\gamma^n$ on $[0,1]$.  By this formula,
$\Vol(L)=\Vol(L_0)=1$, and so $f(\gamma)=1-\gamma^n$ as claimed.

We may also compute the volume directly.  Choosing
an affine chart $\AA^n$ centered at $x$ we can identify the global sections of $\Osh_{\PP^n}(m)$ with polynomials
in $z_1$, \ldots, $z_n$ of degree $\leq m$.  For $\gamma$ rational and $m$ such that $m\gamma$ is an integer,
the global sections of $mL_{\gamma}$ can be identified with the subspace of those polynomials whose lowest degree
term has degree at least $m\gamma$.  This subspace therefore has dimension $\binom{m+n}{n}-\binom{m\gamma-1+n}{n}$.
From the definition of volume we then compute that $\Vol(L_{\gamma}) = 1-\gamma^n$ as before.

\noindent
{\bf Example.} \label{ex:Ar-P1xP1}
Let $X=\PP^1\times \PP^1$, $L=\Osh_{X}(d_1,d_2)$ with $d_1,d_2\geq 1$, and let $x\in X$ be any point.
Choosing an affine chart $\AA^2$ centered at $x$, global sections of $mL$ may be identified
with the polynomials in two variables $z_1$, $z_2$ on $\AA^2$ spanned by the monomials $z_1^az_2^b$ with
$0\leq a\leq md_1$ and $0\leq b\leq md_2$.    For rational $\gamma$ and $m$ such that $m\gamma$ is integral,
the global sections of $mL_{\gamma}$ may be identified with the subspace of these polynomials vanishing to
order $\geq m\gamma$ at $x$, or equivalently, with the subspace spanned by the monomials $z_1^az_2^b$
with $a+b\geq m\gamma$.
For $\gamma> d_1+d_2$ we therefore have $H^0(mL_{\gamma})=0$ for all $m>0$, and for $\gamma< d_1+d_1$ we
have $H^0(mL_{\gamma})\neq 0$ for all $m\gg 0$ (and sufficiently divisible so that $m\gamma$ is an integer).
Thus $\geffx(L)=d_1+d_2$.

The exponent vectors $(a,b)$ of the monomials in $H^0(mL)$ are the lattice points in the rectangle
$[0,md_1]\times[0,md_2]$, while those of the monomials in the subspace $H^{0}(mL_{\gamma})$ are the subset of those
lattice points satisfying $a+b\geq m\gamma$.
Scaling the rectangle by $1/m$ and letting $m\to\infty$, we conclude that for
$\gamma\in [0,d_1+d_2]$ the ratio $f(\gamma) = \Vol(L_{\gamma})/\Vol(L)$ is the fraction of
the area of the rectangle $[0,d_1]\times [0,d_2]$ satisfying $a+b\geq\gamma$ (the shaded region shown below):

\bigskip

\newgray{lgray}{0.95}

\begin{centering}
\begin{pspicture}(0,0)(7,4.7)
\pspolygon[fillstyle=solid,fillcolor=lgray,linecolor=white](1,4)(7,4)(7,0)(5,0)
\pspolygon(0,0)(0,4)(7,4)(7,0)
\psline[linecolor=gray](0,5)(6,-1)
\rput(6.3,-0.5){\tiny\color{gray}$a+b=\gamma$}
\rput(-0.5,2){$d_1$}
\rput(3.5,4.5){$d_2$}
\end{pspicture}\\
\fig\label{fig:rectangle} \\
\end{centering}

\newcommand{\fcurve}[3]{
\parametricplot{0}{#1}{t 1 t t mul #1 #2 2 mul mul div sub}
\psline(!#1 1 #1 #2 2 mul div sub)(!#2 #1 #2 2 mul div)
\parametricplot{#2}{#3}{t #3 t sub dup mul 2 #1 #2 mul mul div}
\rput(#1,-0.150){\tiny $d_1$}
\rput(#2,-0.150){\tiny $d_2$}
\psline[linecolor=gray](#1,0.05)(#1,-0.05)
\psline[linecolor=gray](#2,0.05)(#2,-0.05)
}

\newcommand{\lcurve}[3]{
\parametricplot[linecolor=gray]{0}{#3}{t 1 t t mul #1 #2 2 mul mul div sub}
}

\newcommand{\filllcurve}[3]{
\pscustom[linecolor=lgray,fillstyle=solid,fillcolor=lgray]{%
\parametricplot[linecolor=gray]{0}{#3}{t 1 t t mul #1 #2 2 mul mul div sub}
\psline(#3,0)(0,0)(0,1)
}
}

\bigskip
\bigskip
\bigskip
So that for any $x\in X$ (and assuming that $d_1\leq d_2$ for the purposes of this formula)

\bigskip
\bigskip
\bigskip

\begin{tabular}{ccc}
$f(\gamma) = \left\{{
\begin{array}{cl}
1-\frac{\gamma^2}{2d_1d_2} & \mbox{if $0\leq \gamma\leq d_1$} \\
1+\frac{d_1}{2d_2} - \frac{\gamma}{d_2} & \mbox{if $d_1\leq \gamma \leq d_2$} \rule{0cm}{0.6cm}\\
\frac{(d_1+d_2-\gamma)^2}{2d_1d_2} & \mbox{if $d_2\leq \gamma\leq d_1+d_2$} \rule{0cm}{0.6cm}\\
\end{array}}\right.
$
&
\rule{0.5cm}{0cm}
&
\begin{tabular}{c}
\begin{pspicture}(0,-0.5)(7,2)
\psset{yunit=2}
\SpecialCoor
\filllcurve{3}{5}{5.477225575}
\lcurve{3}{5}{5.477225575}
\fcurve{3}{5}{8}
\psline[linecolor=gray,arrows=->](0,0)(0,1.2)
\psline[linecolor=gray,arrows=->](0,0)(8,0)
\rput(4,0.8){\tiny $y=f(\gamma)$}
\rput(8.2,0){\tiny $\gamma$}
\rput(0,1.35){\tiny $y$}
\end{pspicture}
\end{tabular}\\
& & \fig\label{fig:P1xP1-graph}
\end{tabular}

\bigskip

\noindent
with area $\Ar_x(L)=\int_{0}^{d_1+d_2} f(\gamma)\,d\gamma = \frac{d_1+d_2}{2}$.
(The shaded region in Figure \ref{fig:P1xP1-graph} is not connected with the shaded region in
Figure \ref{fig:rectangle} and will be explained below.)

\begin{lemma}\label{lem:beta-bound}
For any ample $L$, $x\in X$, and $\gamma\geq 0$ we have
$\Vol(L_{\gamma})\geq \Vol(L)-(\mult_{x}X)\cdot \gamma^n.$
\end{lemma}

\bpf
Since $\Vol(\cdot)$ is a continuous function, it suffices to prove the formula for rational $\gamma$.
For $m$ large and sufficiently divisible (i.e., so that $m\gamma$ is an integer) we have the exact sequence of sheaves

\begin{equation}\label{eqn:exact-gammaE-sequence}
0\longrightarrow mL_{\gamma} \stackrel{\cdot m\gamma E}{\longrightarrow} mL_0
\longrightarrow mL_0|_{m\gamma E}\longrightarrow 0
\end{equation}

\noindent
on $\Xtil$ where $m\gamma E$ is the subscheme defined by the $(m\gamma)^{\mbox{\scriptsize th}}$ power
of the ideal sheaf of the Cartier divisor $E$, and where $L_0=\pi^{*}L$.
This yields an exact sequence on global sections:
$$0\longrightarrow\Gamma(\Xtil,mL_{\gamma})\longrightarrow\Gamma(\Xtil,mL_0) \longrightarrow
\Gamma(m\gamma E, mL_0|_{m\gamma E}).  $$

Since $h^0(mL_0)=h^0(mL) = \frac{\Vol(L)}{n!}m^n  + O(m^{n-1})$ the lemma will follow if we show that
$h^0(mL_0|_{m\gamma E}) \leq \frac{\mult_{x}X}{n!}(\gamma m)^n  + O(m^{n-1}).$

Because $L$ can be trivialized in a neighbourhood of $x$, $L_0=\pi^{*}L$ is trivial in a neighbourhood of $E$,
and hence $L_0|_{m\gamma E} = \Osh_{m\gamma E}$.  Let $\IshE$ be the ideal sheaf of $E$ on $\Xtil$.
For any $\ell\geq 1$ we have $\IshE^{\ell}/\IshE^{\ell+1}=\Osh_{E}(-\ell E)$, and thus the exact sequence of sheaves

\begin{equation}\label{eqn:O-gamma-filtration}
0\longrightarrow \Osh_{E}(-\ell E) \longrightarrow \Osh_{(\ell+1)E} \longrightarrow \Osh_{\ell E} \longrightarrow 0.
\end{equation}

\noindent
This gives the inductive estimate

\begin{equation}\label{eqn:inductive-estimate}
h^0(L_0|_{m\gamma E}) = h^0(\Osh_{m\gamma E}) \leq \sum_{\ell=0}^{m\gamma-1} h^0(\Osh_{E}(-\ell E)).
\end{equation}

Choose an embedding $X\hookrightarrow \PP^m$ and let $\Ptil^m$ be the blow up of $\PP^m$ at the image of $x$,
with exceptional divisor $E'\cong \PP^{m-1}$.  Then $\Xtil$ is the proper transform of $X$ in $\Ptil^{m}$,
and $E=\Xtil\cap E'$.  Furthermore, $E$ has degree $\mult_{x}X$ as an $(n-1)$-dimensional subvariety of $\PP^{m-1}\cong E'$.     We thus have

\begin{equation}\label{eqn:O(kE)-estimate}
h^0(\Osh_{E}(-\ell E)) = \frac{\mult_{x}X}{(n-1)!} \ell^{n-1} + O(\ell^{n-2})\,\,\,\mbox{for $\ell\gg 0$}
\end{equation}

\noindent
since $\Osh_{E'}(-E')=\Osh_{\PP^{m-1}}(1)$, and so $h^0(\Osh_{E}(-\ell E)) = h^0(\Osh_{\PP^{m-1}}(\ell )|_{E})$ is
simply given by the Hilbert polynomial of $E$ for large $\ell$.
Summing \eqref{eqn:O(kE)-estimate} and using \eqref{eqn:inductive-estimate} we obtain the estimate
$h^0(mL_0|_{m\gamma E}) \leq \frac{\mult_{x}X}{n!}(\gamma m)^n  + O(m^{n-1})$, proving the lemma.
\epf

\medskip
\noindent
{\bf Remark.} \label{rem:volume-of-big-and-nef}
If $M$ is a big and nef line bundle, then $\Vol(M) = c_1(M)^n$.  In particular, for $\gamma\in [0,\ep_{x}(L)]$,
$\Vol(L_{\gamma}) = c_1(L_{\gamma})^n = c_1(L)^n + (-\gamma)^n E^n = \Vol(L) -(\mult_{x}{X})\cdot \gamma^n$, i.e.,
the lower bound from Lemma \ref{lem:beta-bound} is an equality on $[0,\ep_{x}(L)]\subseteq [0,\geff]$.  In general the
inequality in Lemma \ref{lem:beta-bound} is strict on $(\ep_{x}(L),\geff]$ (i.e., $H^0(mL_0|_{m\gamma E})$ fails
to impose independent conditions on $H^0(mL_0)$ for $\gamma$ in that range).  As an example,
the shaded region in Figure \ref{fig:P1xP1-graph} shows the (normalized) lower bound
$\frac{1}{\Vol(L)}(\Vol(L)-\gamma^2)$ in the
case $X=\PP^1\times\PP^1$. The lower bound is equal to $f(\gamma)$ up until $d_1=\ep_{x}(L)$, but drops away
from $f(\gamma)$ immediately after.

\begin{corollary}\label{cor:beta-bound-irred}
For any ample $L$ and $x\in X$ we have
$\Ar_{x}(L) \geq \frac{n}{n+1} \sqrt[n]{\frac{\Vol(L)}{\mult_{x}X}}\geq
\frac{n}{n+1}\ep_{x}(L)$.  In general, both these inequalities are strict.
\end{corollary}

\bpf
Let $g(\gamma)= 1 - \frac{\mult_{x}X}{\Vol(L)}\gamma^n$ and set $\omega=\sqrt[n]{\frac{\Vol(L)}{\mult_{x}X}}$
(i.e, the solution to $g(\omega)=0$).  By Lemma \ref{lem:beta-bound} we have
$f(\gamma)\geq g(\gamma)$ for all $\gamma\geq 0$.
For any $\gamma\in[0,\omega)$, $g(\gamma)>0$ and hence $f(\gamma)>0$, so we conclude that
$\omega\leq \geff$.  Therefore
$$\Ar_{x}(L) = \int_{0}^{\geff} f(\gamma)\, d\gamma \geq \int_{0}^{\omega} g(\gamma)\, d\gamma =
\textstyle \frac{n}{n+1} \sqrt[n]{\frac{\Vol(L)}{\mult_{x}X}}.$$
The inequality
$ \sqrt[n]{\frac{\Vol(L)}{\mult_{x}X}}= \sqrt[n]{\frac{c_1(L)^n}{\mult_{x}X}}\geq \ep_{x}(L)$ is
\cite[Proposition 5.1.9]{PAG}.
In the example of $X=\PP^1\times \PP^1$, $L=\Osh_{X}(d_1,d_2)$ (with $d_1\leq d_2$) the inequalities are
$\frac{d_1+d_2}{2} > \frac{2}{3}\sqrt{2d_1d_2} > \frac{2}{3}d_1$, i.e, all are strict.  \epf

We now give the definition of $\Ar$ in general.

\begin{definition}\label{def:Ar}
Let $X$ be a variety over $\Spec(k)$, $x\in X(\Qbar)$, and $L$ an ample line bundle on $X$.
Then we define
$$\Ar_{x}(L) = \min(\Ar_{x,\Xbar_1}(L|_{\Xbar_1}),\ldots, \Ar_{x,\Xbar_{\ell}}(L|_{\Xbar_{\ell}})),$$
where $\Xbar_1$,\ldots, $\Xbar_{\ell}$ are the irreducible components of $\Xbar = X\times_{k}\Qbar$ containing $x$.
\end{definition}

It will be important for us that part of Corollary \ref{cor:beta-bound-irred} holds in the general case.

\begin{corollary}\label{cor:beta-bound}
Let $X$ be an irreducible $n$-dimensional variety over $\Spec(k)$.  Then for any $x\in X(\Qbar)$ and any
ample $L$ we have
$\Ar_{x}(L) \geq \frac{n}{n+1}\ep_{x}(L)$.  In general, this inequality is strict.
\end{corollary}

\bpf
Let $\Xbar=X\times_{k}\Qbar$ with irreducible components $\Xbar_1$,\ldots, $\Xbar_\ell$.  Then each component
is $n$-dimensional, hence applying Corollary \ref{cor:beta-bound-irred} we have
$\Ar_{x,\Xbar_i}(L|_{\Xbar_i})\geq \frac{n}{n+1}\ep_{x,\Xbar_i}(L|_{\Xbar_i})$ for each $i=1$, \dots, $\ell$.
By Definition \ref{def:Ar} and Proposition \ref{prop:ex}(f) we then conclude that
$\Ar_{x}(L) \geq \frac{n}{n+1}\ep_{x}(L)$.
\epf

\medskip
\noindent
{\bf Remark.}
\label{rem:f-computed-over-K}
Let $X$ be absolutely irreducible, $x\in X(\Qbar)$ be any point and $K$ its field of definition.
Set $X_{K} = X\times_{k}K$, $\pi_K\colon \Xtil_{K}\longrightarrow X_{K}$ to be the blow
up of $X_K$ at the closed point corresponding to $x$,  and $E_{K}$ to be the exceptional divisor.
For any $\gamma\geq 0$ set $L_{\gamma,K}=\pi_{K}^{*}L_K-\gamma E_K$, where  $L_K$ is the base change of $L$ to $X_{K}$.
We similarly define $\Xtil_{\Qbar}$, $E_{\Qbar}$, and $L_{\gamma,\Qbar}$.
Since $x$ is defined over $K$ it follows that $\Xtil_{K}\times_{K}\Qbar = \Xtil_{\Qbar}$
and hence that $\dim_{\Qbar} H^0(\Xtil_{\Qbar},mL_{\gamma,\Qbar}) = \dim_{K} H^{0}(\Xtil_{K},mL_{\gamma,K})$ for all
$m>0$ and $\gamma\geq 0$ with $m\gamma$ an integer.  Thus the dimension of $mL_{\gamma,\Qbar}$, and hence the
asymptotic growth (i.e., the volume) of $L_{\gamma,\Qbar}$ may be computed ``over $K$''.  In particular,
$\Vol(L_{\gamma,K})/\Vol(L)= \Vol(L_{\gamma,\Qbar})/\Vol(L)$  for all $\gamma\geq 0$.

\medskip
We will investigate $\Ar_{x}(L)$ further in \S\ref{sec:more-about-Ar}.   The facts above
are all we need for our application to the general versions of Roth's theorem.

\medskip

\section{An approximation theorem}
\label{sec:central-thm}

This section is devoted to proving Theorem \ref{thm:simul-approx-I} below.  This theorem is the
central theorem of the paper in the sense that, together with lines of reasoning common in diophantine
approximation%
\footnote{\ldots as well as Propositions \ref{prop:alpha}(f) and \ref{prop:ex}(c), and
Corollary \ref{cor:beta-bound} \ldots}
this theorem implies most of the results in \S\ref{sec:Roth-theorems}--\S\ref{sec:unram-bounds}.

We fix the following notation for the rest of the section.
Let $S$ be a finite set of places of $k$, each extended in some way to $\Qbar$.
Let $X$ be an irreducible $n$-dimensional variety over $\Spec(k)$.
For each $v\in S$ choose a point $x_v\in X(\Qbar)$,
and let $\sepv(\cdot,\cdot)$ be a distance function (as in \S\ref{sec:rat-approx}) computed with respect to $v\in S$.
We are interested in simultaneously approximating each $x_v$, where the distance to $x_v$ is computed with $\sepv$.
To simplify notation, let $\alpha_v = \alpha_{x_v}$ computed with respect to $\sepv$.

A large part of this article is concerned with the approximation constant $\alpha$, and we will state
the results of this section in terms of $\alpha$ and in terms of the usual finiteness conditions;
both versions are equivalent.

\noindent
Fix an ample $\QQ$-bundle $L$.
For a sequence of positive real numbers $\{R_v\}_{v\in S}$ we consider the following two equivalent conditions:

\begin{equation}\label{eqn:alpha-condition}
\rule{1.5cm}{0cm}
\left\{
\begin{tabular}{c}
\mbox{\begin{minipage}{0.85\textwidth}
There is a proper subvariety $Z$ of $X$ such that for all infinite sequences $\xseq$ of distinct points of
$X(k)\setminus Z(k)$, there is at least one $v\in S$ so that
$$\alpha_v(\xseq,L) \geq \frac{1}{R_v}.$$
\end{minipage}} \\
\end{tabular}\right\}\phantom{.}
\end{equation}

\noindent
and

\begin{equation}\label{eqn:finite-condition}
\rule{1.5cm}{0cm}
\left\{
\begin{tabular}{c}
\mbox{\begin{minipage}{0.85\textwidth}
There is a proper subvariety $Z$ of $X$ such that for any collection $\{\delta_v\}_{v\in S}$
with each $\delta_v>0$,
there are only finitely many solutions $y\in X(k)\setminus Z(k)$ to
$$\sepv(x_v,y) \leq  H_{L}(y)^{-(R_v+\delta_v)}\,\,\,\mbox{for all $v\in S$}.$$
\end{minipage}} \\
\end{tabular}\right\}.
\end{equation}

\medskip
We think of the constants $R_v$ as ``Roth constants''
for this approximation problem, generalizing $R=2$ in the case $X=\PP^1$.
Although indexed by the place $v\in S$, it is the local geometry around $x_v$,
also indexed by $v$, which influences the constants $R_v$ for which \eqref{eqn:alpha-condition} and
\eqref{eqn:finite-condition} hold.

\medskip

\begin{theorem}\label{thm:simul-approx-I}
Given a collection $\{R_v\}_{v\in S}$ of positive real numbers, if

\begin{equation}\label{eqn:R-condition}
\sum_{v\in S} \Ar_{x_v}(L)R_v > 1
\end{equation}

\noindent
then \eqref{eqn:alpha-condition} and \eqref{eqn:finite-condition} hold with respect to the collection
$\{R_v\}_{v\in S}$.
\end{theorem}

It is sometimes common (e.g., as in the Schmidt subspace theorem) to switch the order of quantifiers
in condition \eqref{eqn:finite-condition} and specify $\{\delta_{v}\}_{v\in S}$ before having to specify $Z$.
In this case one may relax the strict inequality in \eqref{eqn:R-condition} to allow equality.
We record this variation for future use.

\begin{corollary}\label{cor:switch-quantifiers}
If $\{R_v\}_{v\in S}$ is a sequence of positive real numbers such that
$\sum_{v\in S} \Ar_{x_v}(L)R_v \geq  1$, then given any sequence $\{\delta_v\}_{v\in S}$ of positive real
numbers there is a proper subvariety $Z$ so that there are only finitely many
solutions $y\in X(k)\setminus Z(k)$ to
$$\sepv(x_v,y) \leq  H_{L}(y)^{-(R_v+\delta_v)}\,\,\,\mbox{for all $v\in S$}.$$
\end{corollary}

\noindent
{\em Proof of Corollary \ref{cor:switch-quantifiers}. }
Given such collections $\{R_v\}_{v\in S}$ and $\{\delta_v\}_{v\in S}$ set $\delta'_v=\frac{\delta_v}{2}$ and
$R'_v=R_{v}+\delta'_v$ for each $v\in S$.  Since each $\Ar_{x_v}(L)>0$, and since each $R'_v>R_v$, we conclude
that $\sum_{v} \Ar_{x_v}(L)R'_v>1$, and thus we may apply Theorem \ref{thm:simul-approx-I} to the collection
$\{R'_v\}_{v\in S}$.  By the theorem, there exists a $Z\subset X$ such that there are only finitely many $y\in X(k)\setminus Z(k)$ satisfying
$$\sepv(x_v,y) \leq  H_{L}(y)^{-(R_v'+\delta_v')}\,\,\,\mbox{for all $v\in S$}.$$
Since $R'_v+\delta'_v = R_v+\delta_v$ for all $v\in S$, this establishes the corollary.
\epf

The following slight improvement in Theorem \ref{thm:simul-approx-I} is useful as a first step in induction.
The content is that if $\dim X=1$ one may again relax the strict inequality in \eqref{eqn:R-condition} to allow
equality and also take $Z=\emptyset$.

\begin{corollary}\label{cor:equality-in-dim-one}
Assume that $\dim X=1$.  If $\{R_v\}_{v\in S}$ is a sequence of positive real numbers such that
$\sum_{v\in S} \Ar_{x_v}(L)R_v \geq  1$, then given any sequence $\{\delta_v\}_{v\in S}$ of positive real numbers,
there are only finitely many solutions $y\in X(k)$ to

\begin{equation}\label{eqn:curve-finite}
\sepv(x_v,y) \leq  H_{L}(y)^{-(R_v+\delta_v)}\,\,\,\mbox{for all $v\in S$}.
\end{equation}

\noindent
In particular, there is no need for an exceptional subset $Z$.
\end{corollary}

\bpf
By Corollary \ref{cor:switch-quantifiers} there is a $Z$ (depending on our choice of $\{\delta_{v}\}_{v\in S}$)
so that there are only finitely many solutions $y\in X(k)\setminus Z(k)$ to \eqref{eqn:curve-finite}.
Since $Z$ is of dimension zero, $Z(k)$ is finite, and so there are only finitely
many $y\in X(k)$ satisfying \eqref{eqn:curve-finite}.  \epf

We will prove Theorem \ref{thm:simul-approx-I} at the end of this section,
after dealing with some preliminary material.
The key input in the proof of the theorem is the powerful and flexible
approximation theorem of Faltings-W\"{u}stholz, which we now outline in the form we will use.

For each $v\in S$ let $\Kv$ be a finite extension of $k$ (we use this notation so that there is no confusion with
$K_v$, the completion of a field $K$ at $v$).  Let $L$ be a very ample line bundle on $X$ and set
$V=\Gamma(X,L)$.  For each $v\in S$, set $V_{\Kv} = V\otimes_{k} \Kv$.  We suppose that for each $v$ we're
given a decreasing filtration

$$V_{\Kv} = V_{\Kv}^{0} \supseteq V_{\Kv}^{1} \supseteq V_{\Kv}^{2} \supseteq \cdots \supseteq V_{\Kv}^{r_v}
\supseteq V_{\Kv}^{r_v+1}=\{0\}$$

\noindent
of $\Kv$-vector spaces, and an increasing sequence $0< c_{v,1} < c_{v,2} < \cdots < c_{v,r_v}$ of positive
real numbers.  For any $k$-subspace $W\subseteq V$ we set $W_{\Kv} = W\otimes_{k} \Kv$ and
$W_{\Kv}^{j} = V_{\Kv}^j \cap W_{\Kv}$ for $j=1,\ldots, r_v+1$.
We define the $v$-th piece of the slope, $\mu_{v}(W)$, by
$$\mu_{v}(W) = \frac{1}{\dim W} \sum_{j=1}^{r_v} c_{v,j} \dim\left(W_{\Kv}^{j}/W_{\Kv}^{j+1}\right)
=\frac{1}{\dim W}\sum_{j=1}^{r_v} c_{v,j}\left(\dim W_{\Kv}^{j} - \dim W_{\Kv}^{j+1}\right).$$
Finally, we define the slope $\mu(W)$ of $W$ to be
$\mu(W) = \sum_{v\in S} \mu_{v}(W).$

Although there are an infinite number of possible subspaces $W$, once the data of the filtration is fixed,
there are only finitely
many possible values for the slope.  Let $\muo$ be the largest slope appearing, and among the subspaces of slope $\muo$,
let $\Whn$ be one of the largest dimension.   A short calculation shows that if $W'$ is a subspace with slope $\muo$,
then $W'\subseteq \Whn$, so $\Whn$ is the largest subspace of slope $\muo$ both in dimension and in the
partial ordering induced by inclusion.
The bundle $\Whn$ is often called the ``maximal destabilizing bundle'', or the ``first step in the
Harder-Narasimhan filtration''.    We now fix $\Whn$ to be this subspace (rather than an arbitrary variable subspace).
Note that $\Whn\neq \{0\}$.

\np
Given the destabilizing bundle $\Whn$, set $Z=\left\{z\in X \st  \mbox{$s(z)=0$ for all $s\in \Whn$}\right\}$.  Since
$\Whn$ is a nonzero subspace of $H^0(X,L)$,  $Z$ is a proper subvariety of $X$.

\np
Next, for each $v\in S$ we fix a $v$-adic norm on $L$ extending our chosen valuation $v$.   Given a global
section $s$ of $L$ and a point $y\in X(\Qbar)$ we denote the $v$-adic norm of $s$ in the fibre at $y$ by
$\Cnrm{s(y)}_{v}$.

Choosing an affine open set $U$ where $L$ is trivial, each global section $s$ may be identified with a function
$g_s$ via the trivialization.    The only fact about the norm which we will need is that
for any $x\in U(\Qbar)$, locally  (with respect to $\sepv$) near $x$ the
functions $\Cnrm{s(\cdot)}_{v}$ and $\nrm{g_s(\cdot)}_{v}$ are equivalent.
In particular, if $\sepv(x_v,y_i)\to 0$ as $i\to \infty$
then the asymptotics of $\Cnrm{s(y_i)}_{v}$ and $\nrm{g_s(y_i)}_v$ are the same.

Finally, for each $v\in S$ and $j\in\{1,\ldots, r_v\}$
we choose a $\Kv$-basis $\{s_{v,j,\ell}\}_{\ell\in I_{v,j}}$ for $W_{\Kv}^{j}$.
With this notation, the theorem \cite[Theorem 9.1]{FW} of Faltings-W\"{u}stholz is:

\begin{theorem} \label{thm:FW-simul}
(Faltings-W\"{u}stholz) If $\mu(\Whn)>1$ then there are only finitely many solutions
$y\in X(k)\setminus Z(k)$ such that
$$\Cnrm{s_{v,j,\ell}(y)}_v < H_{L}(y)^{-c_{v,j}}\,\,\,\mbox{for all $v\in S$,
$j\in \{1,\ldots, r_v\}, \ell\in I_{v,j}$}.$$
\end{theorem}

By definition of $\Whn$ we have the elementary estimate $\mu(\Whn)\geq \mu(V)$
and we will ensure the hypothesis $\mu(\Whn)> 1$ by simply checking that $\mu(V)>1$.
The next lemma allows us deduce $\mu(V)>1$ from condition \eqref{eqn:R-condition}.

\begin{lemma}\label{lem:approx-BR}
Suppose that $f$ is a continuous function defined on an interval $[0,\geff]$ with $f(\geff)=0$,
and set $\Ar=\int_{0}^{\geff} f(\gamma)\,d\gamma$.
Given any positive real number $R$ and any $\delta'>0$ it is possible to choose a non-negative integer $r$ and
rational numbers
$0=\gamma_0<\gamma_1 < \gamma_2 < \cdots <\gamma_r < \geff$ so that, if we define $c_j$ by $c_j = \gamma_j R$
and set $\gamma_{r+1}=\geff$, we have
$$\sum_{j=1}^{r} c_j\left(f(\gamma_j)-f(\gamma_{j+1})\right)> \Ar R-\delta'.$$
\end{lemma}

\bpf
Substituting $c_j=R\gamma_j$ we have
$$
\sum_{j=1}^{r} c_j\left(f(\gamma_j)-f(\gamma_{j+1})\right) =
\sum_{j=1}^{r} R\gamma_j\left(f(\gamma_j)-f(\gamma_{j+1})\right) =
R\left({ \sum_{j=1}^{r} (\gamma_j-\gamma_{j-1}) f(\gamma_j) }\right),
$$
and we recognize the final term as $R$ times the right-hand-sum approximation to the integral of $f$.
By choosing $r$ and rational $\gamma_1$,\ldots, $\gamma_r\in (0,\geff)$
we can clearly arrange for this approximation to be as close as we want to $\Ar$.  \epf

\label{sec:proof-of-approx-I}
\np
{\em Proof of Theorem \ref{thm:simul-approx-I}: \/}
The idea of the proof is simple.
For each $v\in S$ we filter the space of global sections of $mL$ (with $m\gg 0$) by the order
of vanishing at $x_v$.  (Using sections of $mL$ instead of $L$ allows us to get the better estimate
on the resulting slope.)  Writing out what the Faltings-W\"ustholz theorem gives us with respect to the resulting
filtration yields Theorem \ref{thm:simul-approx-I}.  We now explicitly carry out these steps.

If $X(k)$ is not Zariski-dense, then \eqref{eqn:alpha-condition} and \eqref{eqn:finite-condition} hold with
$Z=\overline{X(k)}$.  We may therefore assume that $X(k)$ is Zariski dense and hence by
Lemma \ref{lem:irreducible-Z} that $X$ is geometrically irreducible.

For each $v\in S$ we let $\Kv$ be the field of definition of $x_v$, $\pi_{v}\colon\Xvtil\longrightarrow \Xv$
the blow up of $\Xv=X\times_{k}\Kv$ at the closed point corresponding to $x_v$. Let $\Ev$ denote the exceptional
divisor  and for $\gamma\geq 0$ we put
$\Lvg=\pi_{v}^{*}L-\gamma\Ev$.  Then $\pi_v$, $\Ev$ and $\Lvg$ are all varieties over $\Spec(\Kv)$.
For $\gamma\geq 0$ set $f_{v}(\gamma) = \frac{\Vol(\Lvg)}{\Vol(L)}$.
Since $X$ is geometrically irreducible (and by the remark on page \pageref{rem:f-computed-over-K}) the integral
of $f_v$ is $\Ar_{x_v}(L)$.
By Lemma \ref{lem:approx-BR}
and the hypothesis \eqref{eqn:R-condition},
for each $v\in S$ we may choose $r_v$ and rational
$0<\gamma_{v,1} < \cdots < \gamma_{v,r_v} < \geffxv(L)$ so that

$$\sum_{v\in S} \left(\sum_{j=1}^{r_v} c_{v,j}\left(f_v(\gamma_{v,j})-f_v(\gamma_{v,j+1})\rule{0cm}{0.4cm}\right)\right) > 1,$$
with $c_{v,j} = R_v \gamma_{v,j}$ for $v\in S$, $j=1,\ldots, r_v$.  By taking $m$ sufficiently divisible we may ensure
that $mL$ is an integral line bundle and that each $m\gamma_{v,j}$ is an integer.

For any $\gamma\geq 0$, ${\dim_{\Kv}\Gamma(m \Lvg)}/{\dim_{\Kv}\Gamma(m L)}\to f_{v}(\gamma)$ as $m\to\infty$,
and so by taking $m$ sufficiently large we may also ensure that
each $\dim_{\Kv}{\Gamma(m \Lvgj)}/\dim_{\Kv}{\Gamma(m L)}$ is sufficiently close to $f_v(\gamma_{j,v})$ so that

\begin{equation}\label{eqn:simul-slope-condition}
\sum_{v\in S} \frac{1}{\dim_{\Kv}\Gamma(mL)} \left( \sum_{j=1}^{r_v} c_{v,j} \left( \dim_{\Kv}{\Gamma(m \Lvgj)}-
\dim_{\Kv}{\Gamma(m \Lvgjpo)}\right)\right) > 1.
\end{equation}

Set $V=\Gamma(X,mL)$ and we identify $V_{\Kv}$ with $\Gamma(\Xvtil,m \Lvz)$ as vector spaces\footnote{If $X$ is not normal,
$\Gamma(X,mL)\otimes_{k}\Kv$ may only be a proper subspace of $\Gamma(\Xvtil,m\Lvz)$. However, since the volume
is a birational invariant, the asymptotic calculations go through without change and we omit further mention of this
detail.}.
We give a decreasing filtration on each
$V_{\Kv}$ by setting $V_{\Kv}^{j} = \Gamma(m L_{\gamma_{j,v},v})$ for $j=1,\ldots, r_v$, and
choose a $\Kv$-basis $\{s_{v,j,\ell}\}_{\ell\in I_{v,j}}$ for each $V_{\Kv}^{j}$.
As above we let $\Whn$ be the maximal destabilizing subspace and $Z$ the base locus of the
sections in $\Whn$.
Equation \eqref{eqn:simul-slope-condition} is the statement that $\mu(V)>1$, and so we conclude that $\mu(\Whn)>1$ too.
We may therefore apply Theorem \ref{thm:FW-simul} and conclude that there are only finitely many solutions
$y\in X(k)\setminus Z(k)$ to

\begin{equation}\label{eqn:contradict-me}
\Cnrm{s_{v,j,\ell}(y)}_{v}^{\frac{1}{m R_v \gamma_{v,j}}}  \, H_{L}(y) < 1\,\,\,
\mbox{for all $v\in S$, $j=1,\ldots, r_v,\ell\in I_{v,j}$}.
\end{equation}

Now suppose that \eqref{eqn:alpha-condition} is false for this choice of $Z$.
Then there exists a sequence $\yseq$ of $k$-points of $X$,
with no $y_i$ contained in $Z$ such that $\alpha_{v}(\yseq,L) < \frac{1}{R_v}$ for each $v\in S$.
This means that for all sufficiently small $\delta'>0$, and each $v\in S$,
$\sepv(x_v,y_i)^{\frac{1}{R_v}-\delta'} H_{L}(y_i) \to 0\,\,\,\mbox{as $i\to \infty$}$.

Since each $s_{v,j,\ell}$ is in $V_{\Kv}^{j}$, each $s_{v,j,\ell}$ is in the $(m\gamma_{v,j})^{\mbox{\scriptsize th}}$
power of the maximal ideal of $x_v$, and so for any $\delta>0$ and for large enough $i$ (depending on $\delta$) we have
$\Cnrm{s_{v,j,\ell}(y_i)}_v \leq \sepv(x_v,y_i)^{m\gamma_{v,j}-\delta}$.
But then for large enough $i$

\begin{equation}\label{eqn:Growth-bound}
\Cnrm{s_{v,j,\ell}(y_i)}_{v}^{\frac{1}{m R_v \gamma_{v,j}}}  \, H_{L}(y_i)  \leq
\sepv(x_v,y_i)^{\frac{1}{R_v} - \left(\frac{\delta}{m R_v \gamma_{v,j}}\right)} H_{L}(y_i)
\end{equation}

\noindent
{for all $v\in S$, $j=1,\ldots, r_v,\ell\in I_{j,v}$}.
For small enough $\delta>0$ the right hand side of
\eqref{eqn:Growth-bound} tends to $0$ as $i\to \infty$.  This contradicts \eqref{eqn:contradict-me} and therefore
assertion \eqref{eqn:alpha-condition} holds.  This proves Theorem \ref{thm:simul-approx-I}. \epf

\section{Roth theorems}
\label{sec:Roth-theorems}

Let $X$ be an irreducible $n$-dimensional variety over $\Spec(k)$.
In this section we present theorems giving lower bounds for $\alpha_{x}(L)$ independent of the field
of definition of $x\in X(\Qbar)$, and in particular lower bounds in terms of $\ep_{x}(L)$.
In the remaining sections of the paper we will deal with simultaneous approximation but in order to clarify the ideas
we start by approximating with respect to a single place, either archimedean or non-archimedean.
As in the beginning of the paper, we fix a place $v_0$ of $k$, an extension $v$ of $v_0$ to $\Qbar$,
and compute $\alpha_{x}$ with respect to $\sepv(\cdot,\cdot)$.

\begin{theorem}\label{thm:RothI}
For any ample $\QQ$-bundle $L$ and any $x\in X(\Qbar)$ either

\begin{enumerate}
\item $\alpha_x(L) \geq \Ar_x(L)$

\medskip
\hfill or \hfill\rule{0.1cm}{0.0cm}

\medskip
\medskip
\item There exists a proper subvariety $Z\subset X$,  irreducible over $\Qbar$,  with $x\in Z(\Qbar)$
so that $\alpha_{x,X}(L) = \alpha_{x,Z}(L|_{Z})$, i.e., ``$\alpha_x(L)$ is computed on a proper subvariety of $X$''.
\end{enumerate}
\end{theorem}

\bpf
If $\alpha_{x}(L) < \Ar_{x}(L)$ then choose any $R>0$ such that $\alpha_{x}(L) < \frac{1}{R} < \Ar_{x}(L)$.
Then $\Ar_{x}(L)R>1$ so by Theorem \ref{thm:simul-approx-I} in the case of a single place we conclude that there
is a proper subvariety $Z$ such that for all sequences $\xseq$ of
$k$-points with $\alpha_{x}(\xseq,L)\leq \frac{1}{R}$, all but finitely
many of the points lie in $Z$. We conclude that $\alpha_{x,Z}(L|_{Z}) = \alpha_{x}(L)$.
To see that we may assume that $Z$ is geometrically irreducible, apply
Lemma \ref{lem:irreducible-Z} to $Z$, and use Proposition \ref{prop:alpha}(f) to replace $Z$ by a component
of the resulting variety $Y$.
\epf

\np
By Corollary \ref{cor:beta-bound} we have $\Ar_{x}(L)\geq \frac{n}{n+1}\ep_x(L)$.
Thus Theorem \ref{thm:RothI} implies the weaker theorem:

\begin{theorem}\label{thm:RothII} {\rm (}Schmidt type{\rm )}
Under the same hypothesis as Theorem \ref{thm:RothI}, either

\begin{enumerate}
\item $\alpha_x(L) \geq \frac{n}{n+1} \ep_{x}(L)$

\medskip
\hfill or \hfill\rule{0.1cm}{0.0cm}

\medskip
\medskip
\item $\alpha_{x}(L)$ is computed on a proper subvariety $Z$ of $X$ (irreducible over $\Qbar$, as above).
\end{enumerate}
\end{theorem}

This immediately yields

\begin{theorem}\label{thm:RothIII} {\rm (}Roth type{\rm )}
With the same hypotheses as above, $\alpha_{x}(L)\geq \frac{1}{2}\ep_{x}(L)$, with equality if and only if
both $\alpha$ and $\ep$ are computed on a $k$-rational curve $C$ such that
(1) $C$ is unibranch at $x$, (2) $\kappa(x)\neq k$,
(3) $\kappa(x)\subset k_v$, and (4) $\ep_{x,C}(L|_C) = \ep_{x_v,X}(L)$.
\end{theorem}

\noindent
{\em Proof of Theorem \ref{thm:RothIII}: \/}
If $\alpha_{x}(L)\geq \frac{n}{n+1} \ep_{x}(L)$ and $n>1$ then this is stronger than
$\alpha_{x}(L)\geq \frac{1}{2}\ep_x(L)$ so we are done.  If not, then by Theorem \ref{thm:RothII} we pass to a
smaller irreducible subvariety.  Since the Seshadri constant can only go up when restricting to a subvariety (Proposition \ref{prop:alpha}(c)),
we are done by induction.  Finally, in the case of equality we conclude that we must have gone all
the way down to a curve $C$ (irreducible over $\Qbar$),
and $\ep_x$ must also be computed on $C$, or the inequality would be strict (i.e., (4) above holds).
Conditions (1), (2), and (3) then follow from Lemma \ref{lem:roth-bound-for-curves}.

Conversely, if $C$ is a $k$-rational curve passing through $x$ and satisfying (1), (2), and (3) above then
Lemma \ref{lem:roth-bound-for-curves}  gives $\alpha_{x,C}(L|_C) = \frac{1}{2}\ep_{x,C}(L|_C)$.
If in addition (4) holds then we have
$$\textstyle\frac{1}{2}\ep_{x,X}(L) = \frac{1}{2}\ep_{x,C}(L|_C)=\alpha_{x,C}(L|_C) \geq \alpha_{x,X}(L),$$
where the last inequality is Proposition \ref{prop:alpha}(c).  By the first part of the theorem we always
have $\alpha_{x,X}(L)\geq \frac{1}{2} \ep_{x,X}(L)$, and thus equality must hold.
\epf

Here is a form of Theorem \ref{thm:RothIII} expressed in language closer to the usual statement of Roth's theorem.

\begin{corollary}\label{cor:one-place-Roth}
With the same hypothesis as above, for any $\delta>0$ there are only finitely many solutions $y\in X(k)$ to
$$\sepv(x,y) < {H_{L}(y)^{-\left(\frac{2}{\ep_{x}(L)}+\delta\right)}}.$$
\end{corollary}

\bpf
This is immediate from Theorem \ref{thm:RothIII} and Proposition \ref{prop:equiv-alpha}. \epf

Other variations on the deduction of Theorem \ref{thm:RothIII} from Theorem \ref{thm:RothII} are possible;
here are two examples.

\begin{corollary}\label{cor:control-dim-Z}
If $\alpha_{x}(L) < \frac{m}{m+1}\ep_{x}(L)$ for some $m<n$ then $\alpha_{x}(L)$ is computed on a subvariety $Z$
of dimension $<m$.
\end{corollary}

\begin{corollary}\label{cor:no-rat-curve}
If $x\in X(\Qbar)$, and no rational curve passes through $x$
then $\alpha_{x}(L)\geq \frac{2}{3}\ep_{x}(L)$.
Equivalently, for any $\delta>0$ there are only finitely many solutions $y\in X(k)$ to
$$\sepv(x,y) < {H_{L}(y)^{-\left(\frac{3}{2\ep_{x}(L)}+\delta\right)}}.$$
\end{corollary}

\vspace{0.1cm}

\noindent
{\bf Remark.} Theorems \ref{thm:RothI} and \ref{thm:RothII} were stated for a variety $X$ irreducible over $k$
since if $X$ were reducible, and $x\in X(\Qbar)$ not on a component of maximal dimension $n$, the estimate
$\Ar_{x}(L)\geq \frac{n}{n+1}\ep_{x}(L)$ would not hold (the volume only measures top-dimensional asymptotics).
However by using Propositions \ref{prop:alpha}(f) and \ref{prop:ex}(f) to reduce to the irreducible
components of $X$ it follows that Theorem \ref{thm:RothIII}
and Corollaries \ref{cor:one-place-Roth}, \ref{cor:control-dim-Z},
and \ref{cor:no-rat-curve} above still hold when $X$ is reducible.

\section{Simultaneous approximation}
\label{sec:simul-approx}

In this section we apply Theorem \ref{thm:simul-approx-I} to study simultaneous approximation.
As in \S\ref{sec:central-thm} we let $S$ be a finite set of places of $k$, each extended in some way to $\Qbar$  and
$X$ be an irreducible $n$-dimensional variety over $\Spec(k)$.
For each $v\in S$ we choose a point $x_v\in X(\Qbar)$,
and let $\sepv(\cdot,\cdot)$ be the distance function (as in \S\ref{sec:rat-approx}) computed with respect to $v\in S$.
Again, to simplify notation, we set $\alpha_v$ to be $\alpha_{x_v}$ computed with respect to $\sepv$.

We are interested in understanding how well sequences of $k$-points can simultaneously approximate each $x_v$.
An example of this, showing how Theorem \ref{thm:RothIII} and Corollary \ref{cor:one-place-Roth} generalize
to simultaneous approximation, is given in the introduction.
We will also consider the case of sequences $\xseq$ not contained in a subvariety $Z$, and obtain results
along the lines of Theorem \ref{thm:RothII} or Corollary \ref{cor:control-dim-Z}.

There is a general mechanism for proving such simultaneous approximation results due to Mahler.
The basic idea is that these generalizations are equivalent to studying simultaneous approximations with
weights.  We next review these ideas, and then use Theorem \ref{thm:simul-approx-I} to deduce the appropriate
weighted versions.

\begin{definition}
A weighting function $\xi$ is a function $\xi\colon S\longrightarrow [0,1]$ such that $\sum_{v\in S} \xi_v = 1$.
\end{definition}

\noindent
Here, and in the rest of the paper, we use $\xi_v$ for the value of $\xi$ at $v$.

It will be useful to be able to reduce verifying a statement for infinitely many weighting functions to
verifying a slightly stronger statement for only finitely many weighting functions.  This is the purpose
of the following lemma.

\begin{lemma}\label{lem:finite-xi}
Let $S$ be a finite set, and $\{\Delta'_v\}_{v\in S}$ and $\{\Delta_v\}_{v\in S}$ collections of positive
real numbers with $\Delta'_v< \Delta_v$ for all $v\in S$.
Then there exists a finite set $\Xi$ of weighting functions $\xi'\colon S\longrightarrow[0,1]$ so that
given any function $\xi\colon S\longrightarrow \RR_{\geq 0}$ satisfying $\sum_{v\in S} \xi_v\geq 1$
there is a weighting function $\xi'\in \Xi$ satisfying $\xi'_v\Delta'_v \leq \xi_v \Delta_v$ for all $v\in S$.
\end{lemma}

\bpf
Let $N$ be any positive integer so that
$\min_{v\in S}\{\Delta_v/\Delta_v'\}-{\#S}/{N}\geq 1$,
and $\Xi$ the finite set of weighting functions $\xi'\colon S\longrightarrow [0,1]$
such that $N\xi_v'$ is an integer for all $v\in S$ (i.e., all $\xi_v'$ are rational with denominator dividing $N$).
Given a function $\xi\colon S\longrightarrow \RR_{\geq 0}$ with $\sum_{v} \xi_v\geq 1$ set
$$\xi''_{v} = \frac{\left\lfloor{\frac{N\cdot \Delta_v\cdot \xi_{v}}{\Delta'_v}}\right\rfloor}{N}\,\,\,\,
\mbox{for each $v\in S$}.$$

Then $\xi_v'' \leq \frac{\Delta_v}{\Delta'_v} \xi_v$, and so $\xi_v'' \Delta'_v \leq \xi_v \Delta_v$ for each $v\in S$.
Furthermore, each $\xi_v''$ is rational and nonnegative with $N\xi_{v}''$ an integer.  Since
$$\xi''_{v} \geq  \frac{\frac{N\cdot \Delta_v\cdot \xi_{v}}{\Delta_v'}-1}{N} =
\frac{\Delta_v}{\Delta'_v}\xi_v - \frac{1}{N},$$
for each $v\in S$ we conclude that
$$\sum_{v\in S} \xi''_{v} \geq
\sum_{v\in S} \left({\frac{\Delta_v}{\Delta'_v}\xi_v-\frac{1}{N}}\right) \geq
\left(\sum_{v\in S} \min_{v\in S}\left\{\frac{\Delta_v}{\Delta'_v}\right\}\xi_v\right)-\#S/N \geq
\min_{v\in S}\left\{\frac{\Delta_v}{\Delta'_v}\right\}-\#S/N
\geq 1.$$

\noindent
Therefore there exists a weighting function
$\xi'\in \Xi$ with $\xi'_{v} \leq \xi''_v$ for all $v\in S$.
\epf

The following proposition shows the equivalence between statements on simultaneous approximation as in the
introduction, and versions of simultaneous approximation with weights.

\begin{proposition}\label{prop:weighted-equiv}
Let $Z$ be a proper subvariety of $X$,  and $L$ an ample $\QQ$-bundle.
Then for any collection $\{R_v\}_{v\in S}$ of positive real numbers the following conditions are equivalent.

\begin{equation}\label{eqn:A}
\rule{1.5cm}{0cm}
\left\{
\begin{tabular}{c}
\mbox{\begin{minipage}{0.85\textwidth}
For all weighting functions $\xi\colon S\longrightarrow [0,1]$ and all sequences $\xseq$
of $X(k)\setminus Z(k)$ there is at least one $v\in S$ with $\xi_v\neq 0$ such that
$\alpha_{v}(\xseq,L) \geq \frac{1}{R_v\xi_v}$.
\end{minipage}} \\
\end{tabular}\right\}\phantom{.}
\end{equation}

\medskip
\begin{equation}\label{eqn:B}
\rule{1.5cm}{0cm}
\begin{tabular}{c}
\mbox{\begin{minipage}{0.85\textwidth}
For all sequences $\xseq$ of $X(k)\setminus Z(k)$,
$\sum_{v\in S} \frac{1}{R_v\alpha_{v}(\xseq,L)} \leq 1$.
\end{minipage}} \\
\end{tabular}\phantom{.}
\end{equation}

\medskip
\begin{equation}\label{eqn:C}
\rule{1.5cm}{0cm}
\left\{
\begin{tabular}{c}
\mbox{\begin{minipage}{0.85\textwidth}
For all weighting functions $\xi\colon S\longrightarrow [0,1]$ and any collection
$\{\delta_v\}_{v\in S}$ of positive real numbers, there are only finitely many solutions $y\in X(k)\setminus Z(k)$ to
$$ \sepv(x_v,y)^{\frac{1}{R_v}} < H_{L}(y)^{-\xi_v(1+\delta_v)}\,\,\mbox{for all $v\in S$}. $$
\end{minipage}} \\
\end{tabular}\right\}\phantom{.}
\end{equation}

\begin{equation}\label{eqn:D}
\rule{1.5cm}{0cm}
\left\{
\begin{tabular}{c}
\mbox{\begin{minipage}{0.85\textwidth}
For all $\delta>0$ there are only finitely many solutions $y\in X(k)\setminus Z(k)$ to
$$\prod_{v\in S} \sepv(x_v,y)^{\frac{1}{R_v}} < H_{L}(y)^{-(1+\delta)}.$$
\end{minipage}} \\
\end{tabular}\right\}{.}
\end{equation}
\end{proposition}

\newcommand{\labA}{\eqref{eqn:A}}
\newcommand{\labB}{\eqref{eqn:B}}
\newcommand{\labC}{\eqref{eqn:C}}
\newcommand{\labD}{\eqref{eqn:D}}

\noindent
\bpf
$\labA\implies\labB$: Given a sequence $\xseq$ in $X(k)\setminus Z(k)$, set
$D=\sum_{v\in S} \frac{1}{R_v\alpha_{v}(\xseq,L)}$.
If all $\alpha_{v}(\xseq,L)=\infty$ then $D=0$ and so the inequality in $\labB$ holds.  We may therefore
assume that $D\neq 0$, i.e., that there is some $v\in S$
so that $\alpha_{v}(\xseq,L) < \infty$.
Define a weighting function by $\xi_v=\frac{1}{R_v\alpha_{v}(\xseq,L)D}$ for each $v\in S$.  By
$\labA$ there is a $v\in S$ with $\xi_v\neq 0$ so that the inequality in $\labA$ holds.
Writing out the definition of $\xi_v$ and clearing denominators gives $\labB$ (recall that $\alpha_{v}(\xseq,L)>0$
by Proposition \ref{prop:alpha}(d)).

\noindent
$\labB\implies\labA$: If $\labA$ is false then there is a weighting function $\xi$ and
a sequence $\xseq$ in $X(k)\setminus Z(k)$ such that $\frac{1}{R_v\alpha_{v}(\xseq,L)} > \xi_v$
for all $v\in S$ such that $\xi_v\neq 0$.
Summing gives a contradiction to $\labB$.

\noindent
$\labC\implies\labD$: Assume $\labD$ is false and fix any $\delta>0$.  For each of the infinitely many
solutions $y_i$ in $X(k)\setminus Z(k)$ to inequality $\labD$, define $\xi_{v,i}$ so that
$$ \sepv(x_v,y_i)^{\frac{1}{R_v}} = H_{L}(y_i)^{-\xi_{v,i}(1+\delta)}$$
for each $v\in S$.  Taking the product and using the fact the $y_i$ are solutions to the inequality in
$\labD$ we conclude that $\sum_{v\in S} \xi_{v,i} > 1$.  Fix any positive $\delta'$ less than $\delta$. Applying Lemma \ref{lem:finite-xi},
with $\Delta'_v=1+\delta'$ and $\Delta_v=1+\delta$ for all $v\in S$ we obtain
a finite set $\Xi$ of weighting functions so that for any $\xi\colon S\longrightarrow \RR_{\geq 0}$
satisfying $\sum_{v} \xi_v\geq 1$, there is a $\xi'\in \Xi$ satisfying $\xi'_v(1+\delta') \leq \xi_v(1+\delta)$ for
all $v\in S$.  In particular, there is a $\xi'_i\in \Xi$ for each function $\xi_i$ as above.  Since $\Xi$ is a finite
set, by passing to a subsequence of $\yseq$ there is a $\xi'\in \Xi$ which works for all $i$.
Note that since $L$ is ample, we may assume that $H_L(y_i)>1$ for all $i$ by omitting finitely many $y_i$.  Choosing $\delta_v=\delta'$ for each $v$, we have infinitely many solutions to
$ \sepv(x_v,y)^{\frac{1}{R_v}} < H_{L}(y)^{-\xi_v'(1+\delta_v)},\,\,\mbox{for all $v\in S$}$,
contradicting $\labC$.

\noindent
$\labD\implies\labA$:  Assume that $\labA$ is false, so that there is a sequence $\xseq$ in $X(k)\setminus Z(k)$
and a weighting function $\xi$ such that $\alpha_{v}(\xseq,L) < \frac{1}{R_v\xi_v}$ for each $v\in S'$,  where
$S'=\{v\in S\st \xi_v\neq 0\}$.
For $\delta>0$
small enough we will still have $\alpha_{v}(\xseq,L) < \frac{1}{R_v\xi_v(1+\delta)}$ for each $v\in S'$,
and so by definition of $\alpha_v$,
$\sepv(x_v,x_i)^{\frac{1}{R_v\xi_v(1+\delta)}}H_{L}(x_i)\to0$ or equivalently
$\sepv(x_v,x_i)^{\frac{1}{R_v}}H_{L}(x_i)^{\xi_v(1+\delta)}\to0$,
as $i\to\infty$ for all $v\in S'$.
Thus by omitting finitely many of the initial $x_i$ we can make the product
$$\prod_{v\in S'} \left(\sepv(x_v,x_i)^{\frac{1}{R_v}}H_{L}^{\xi_v(1+\delta)}\right)
=
\left(\prod_{v\in S'} \sepv(x_v,x_i)^{\frac{1}{R_v}}\right)H_{L}(x_i)^{(1+\delta)}$$
as small as desired.
The product $\prod_{v\in S\setminus S'} \sepv(x_v,x_i)^{\frac{1}{R_v}}$ is bounded
since each distance function $\sepv(\cdot,\cdot)$ is bounded.
Hence after omitting finitely many of the initial $x_i$ the rest satisfy
$$\prod_{v\in S} \sepv(x_v,x_i)^{\frac{1}{R_v}} < H_{L}(x_i)^{-(1+\delta)}$$
contradicting $\labD$.

\noindent
$\labA\implies\labC$:  Assume that $\labC$ is false.  Then there is a weighting function $\xi$ and a
collection $\{\delta_v\}_{v\in S}$ so that the inequalities in $\labC$ have infinitely many solutions.
Let $S'=\{v\in S\st \xi_v\neq 0\}$ and
let $\yseq$ be a sequence of these solutions ordered by height.  Then
$\sepv(x_v,y_i)^{\frac{1}{R_v\xi_v(1+\delta_v)}}H_L(y_i)<1\,\,\, \mbox{for all $v\in S'$},$
so we conclude that $\frac{1}{R_v\xi_v(1+\delta_v)}\in A_{x_v}(\yseq,L)$.  Thus
$\alpha_{v}(\yseq,L)\leq \frac{1}{R_v\xi_v(1+\delta_v)} < \frac{1}{R_v\xi_v}$ for $v\in S'$, contradicting $\labA$.
\epf

\bigskip
We now use Theorem \ref{thm:simul-approx-I} to establish cases where the equivalent conditions in
Proposition \ref{prop:weighted-equiv} hold.

\begin{theorem}\label{thm:simul-approx-II}
In each of the following two cases there is a proper subvariety $Z\subset X$
so that the equivalent conditions in Proposition \ref{prop:weighted-equiv} hold with respect to the
given collection $\{R_v\}_{v\in S}$.

\medskip

\begin{enumerate}
\item Any choice of $\{R_v\}_{v\in S}$  such that $R_v> \frac{1}{\Ar_{x_v}(L)}$ for each $v\in S$.

\medskip
\item Any choice of $\{R_v\}_{v\in S}$ such that $R_v> \frac{n+1}{n\,\ep_{x_v}(L)}$ for each $v\in S$.
\end{enumerate}

\noindent
In the case $n=\dim X=1$  equality in (a) and (b) is sufficient, and one may take $Z=\emptyset$.
\end{theorem}

\bpf
By Corollary \ref{cor:beta-bound}, $\Ar_{x_v}(L)\geq \frac{n}{n+1} \ep_{x_v}(L)$,
so the condition in (b) implies the condition in (a), and it therefore suffices to prove (a).
Given such a collection $\{R_v\}_{v\in S}$ choose $\{R'_v\}_{v\in S}$ so that
$R_v > R'_v > \frac{1}{\Ar_{x_v}(L)}$ for each $v\in S$.  Applying Lemma \ref{lem:finite-xi} with $\Delta'_v=R'_v$
and $\Delta_v=R_v$ for each $v\in S$, we obtain a finite set of weighting functions $\Xi$
so that for any weighting function $\xi$ there is $\xi'\in \Xi$ satisfying
$\xi'_vR'_v \leq \xi_v R_v$ for all $v\in S$.

Temporarily fix $\xi'\in \Xi$ and set $S'=\{v\in S\st \xi'_v\neq 0\}$.
By our choice of $R_v'$ we have
$\sum_{v\in S'} \Ar_{x_v}(L) \xi'_v R'_v > \sum_{v\in S'} \xi'_v =1$.
Applying Theorem \ref{thm:simul-approx-I} to the collection $\{\xi_v'R'_v\}_{v\in S'}$ we obtain
a proper subvariety $Z_{\xi'}$ such that for any sequence $\xseq$ in $X(k)\setminus Z_{\xi'}(k)$
there is at least one $v\in S'$ with $\alpha_{v}(\xseq,L) \geq \frac{1}{\xi'_v R'_v}$.

Set $Z$ to be the union of the finitely many $Z_{\xi'}$ over all $\xi'\in \Xi$.  Given an arbitrary weighting function $\xi$ and a sequence $\xseq$ in $X(k)\setminus Z(k)$, let $\xi'\in \Xi$
be a weighting function such that $\xi'_v R'_v \leq \xi_v R_v$ for all $v\in S$.
Then since $X(k)\setminus Z(k) \subseteq X(k)\setminus Z_{\xi'}(k)$
we conclude that there is some $v\in S$ with $\xi'_v\neq 0$ so that
$\alpha_{v}(\xseq,L) \geq \frac{1}{\xi'_v R'_v} \geq \frac{1}{\xi_v R_v}$.

Finally the statements about equality in the case $\dim X = 1$ follow as in the proof of Corollary
\ref{cor:equality-in-dim-one}.  (After proving the equivalent version of Corollary \ref{cor:switch-quantifiers}.)
\epf

\medskip
As in Theorem \ref{thm:RothIII} inducting on dimension yields a version with $Z=\emptyset$.

\begin{theorem}\label{thm:simul-approx-III}
Set $R_v=\frac{2}{\ep_{x_v}(L)}$ for each $v\in S$.  Then
the conditions in Proposition \ref{prop:weighted-equiv} hold with respect to the collection
$\{R_v\}_{v\in S}$ and $Z=\emptyset$.
\end{theorem}

\bpf
We will show condition \eqref{eqn:C} holds for the collection $\{R_v\}_{v\in S}$ and with $Z=\emptyset$, i.e,
given any weighting function $\xi$ and any $\delta>0$ we will show that there are only finitely many
solutions $y\in X(k)$ to

\begin{equation}\label{eqn:simul-induct}
\sepv(x_v,y)^{\frac{\ep_{x_v}(L)}{2}} \leq  H_{L}(y)^{-\xi_v\left(1+\delta\right)}\,\,\,\mbox{for all $v\in S$}.
\end{equation}

Suppose a weighting function $\xi$ is given.
When $\dim X =1$ the result we want to prove is Theorem \ref{thm:simul-approx-II}(b).
If $\dim X = n>1$ then $2>\frac{n+1}{n}$ so by Theorem \ref{thm:simul-approx-II}(b) again there is a
proper subvariety $Z'\subset X$ such that there are only finitely many $y\in X(k)\setminus Z'(k)$ satisfying
\eqref{eqn:simul-induct}.
Let $Z_j$ be an irreducible component of $Z'$.
By induction there are only finitely many solutions $y\in Z_j(k)$ to the equations
$$\sepv(x_v,y)^{\frac{\ep_{x_v,Z_j}(L)}{2}} < H_{L}(y)^{-\xi_v\left(1+\delta\right)}\,\,\,\mbox{for all $v\in S$}.$$
Since $\ep_{x_v,Z_j}(L|_{Z_j})\geq \ep_{x_v,X}(L)$ this is a stronger statement than the one we are claiming,
i.e., this implies that there are only finitely many solutions $y\in Z_j(k)$ to \eqref{eqn:simul-induct}.
Thus there are only finitely many solutions $y\in X(k)$ to \eqref{eqn:simul-induct}.  \epf

\medskip
\begin{corollary} \label{cor:alpha-hypersurface}
For any sequence $\xseq$ in $X(k)$

\begin{equation}\label{eqn:alpha-hypersurface}
\sum_{v\in S} \frac{\ep_{x_v}(L)}{\alpha_{v}(\xseq,L)} \leq 2.
\end{equation}

\noindent
Equivalently, for any $\delta>0$ there are only finitely many solutions $y\in X(k)$ to

$$\prod_{v\in S} \sepv(x_v,y)^{\ep_{x_v}(L)} < H_{L}(y)^{-(2+\delta)}.$$

\end{corollary}

\bpf
These are conditions \eqref{eqn:B} and \eqref{eqn:D} respectively when $Z=\emptyset$
and with the choice of $R_v=\frac{2}{\ep_{x_v}(L)}$ for all $v\in S$.
These conditions hold by Theorem \ref{thm:simul-approx-III}.
\epf

\np
{\bf Equality.}
As in Theorem \ref{thm:RothIII} it is useful to study the case of ``equality'' in Theorem
\ref{thm:simul-approx-III}.
By ``equality'' we mean that there is a sequence $\xseq$ so that \eqref{eqn:alpha-hypersurface} is an equality.
Equivalently, in terms of condition \eqref{eqn:A}, equality means that for the given sequence $\xseq$
there a weighting function $\xi$ such that

\begin{equation}\label{eqn:equality-in-A}
\alpha_v(\xseq,L) = \frac{\ep_{x_v}(L)}{2\,\xi_v}\,\,\,\mbox{for all $v\in S'$.}
\end{equation}

\noindent
where $S'=\{v\in S\st \xi_v\neq 0\}$.

\begin{theorem}\label{thm:equality-of-hypersurface}
Suppose that $\xseq$ is a sequence so that we have equality in \eqref{eqn:alpha-hypersurface}.
Let $S'=\{v\in S\st \alpha_{v}(\xseq,L)< \infty\}$ (note that $S'$ is nonempty --- otherwise equality in
\eqref{eqn:alpha-hypersurface} is impossible).  Then there is a $k$-rational curve $C$ containing infinitely
many $x_i$ such that for all $v\in S'$:
(1) $C$ is unibranch at $x_v$ (in particular, $C$ contains $x_v$) (2) $\kappa(x_v)\neq k$,
(3) $\kappa(x_v)\subset k_v$, and (4) $\ep_{x_v,C}(L|_C) = \ep_{x_v,X}(L)$.

Conversely, given a $k$-rational curve $C$ satisfying these
conditions with respect to a non-empty subset $S'\subseteq S$, then for any weighting function
$\xi\colon S'\longrightarrow(0,1]$ (extended by $0$ to a weighting function on $S$)
there is a sequence $\xseq$ of points of $C(k)$ such that \eqref{eqn:equality-in-A} holds.
\end{theorem}

For the converse direction of Theorem \ref{thm:equality-of-hypersurface} we require a
``simultaneous weighted Dirichlet'' result on $\PP^1$, which seems to be generally known,
but for which we could not find a reference.   We first prove this result, which is slightly involved, below.
The proof of Theorem \ref{thm:equality-of-hypersurface} appears after Corollary \ref{cor:simultaneous-dirichlet}.

We are indebted to Damien Roy for the following argument.

\begin{theorem}\label{thm:general-dirichlet}
Let $S$ be a finite set of places of $k$ containing all the archimedean places.
For each place $v$ of $S$, let $e_v\in[0,2)$ be a real number between $0$ and $2$, satisfying $e=\sum_{v\in S} e_v<2$.
For each $v$ in $S$, let $x_v$ be an algebraic element not in $k$ of the completion $k_v$ of $k$ at $v$.
Then there exist infinitely many elements $y\in k$ such that $\nrm{y-x_v}_v<H(y)^{-e_v}$ for all $v$ in $S$.
\end{theorem}

\bpf
Let $R$ be the ring of $S$-integers of $k$, and embed $R$ in $V=\prod_{v\in S} k_v$ via the diagonal embedding.
This embedding also induces an embedding of $R^2$ in $V^2$.
Let $B$ be a large, positive real number, and for each $v\in S$
set $f_v=e_v/e$.

There is a convex subset $D$ of $V^2$ of finite volume (i.e., Haar measure) with the property that $D$ contains
a complete set of representatives for the abelian group $V^2/R^2$.
For any positive real number $N$ let $A_N$ be the set of
vectors $(\avec,\bvec)\in V^2$ such that $\nrm{a_v-x_vb_v}_v<B^{-f_v}$ and $\nrm{b_v}_v<NB^{f_v}$ for all places $v$
in $S$.
Choose $N$ large enough so that the volume of $\frac{1}{2}A_N$ is greater than the volume of $D$ and set $A=A_N$
(note that the choice of $N$ does not depend on $B$).
We will show that $A$ contains a nonzero element of $R^2$ by generalizing the proof of Minkowski's famous
result in the geometry of numbers, as found in \cite[\S1.4]{Ne}.

To see this, consider the sets $\frac{1}{2}A\cap(D+u)$ as $u$ varies over elements of $R^2$.
They clearly cover the set $\frac{1}{2}A$, and for each $u$, we have
$\frac{1}{2}A\cap(D+u)=((\frac{1}{2}A-u)\cap D)+u$.  Therefore, the volume
of $\frac{1}{2}A\cap(D+u)$ is equal to that of $(\frac{1}{2}A-u)\cap D$.  If the sets $(\frac{1}{2}A-u)\cap D$ were
pairwise disjoint, then by summing over $u$, we would find that the volume of $\frac{1}{2}A$ is at
most the volume of $D$, in contradiction to our choice of $A$.  We conclude
that the sets $(\frac{1}{2}A-u)\cap D$ are not disjoint.

We may therefore find elements $u,v\in R^2$ and $a_1,a_2\in A$ such that $\frac{1}{2}a_1-u=\frac{1}{2}a_2-v$.
Since $A$ is convex and closed under multiplication by $-1$, it follows that $u-v$ is a nonzero
element of $R^2\cap A$, as desired.
Let $(a_B,b_B)\in R^2$ be such an element (so $a_B,b_B\in R$, and for all $v\in S$
$\nrm{a_B-x_vb_B}_v<B^{-f_v}$ and $\nrm{b_B}_v<NB^{f_v}$).
Since at least one of the $x_v$ is not in $k$, at least one $\nrm{a_B-x_vb_B}_v\neq 0$, and
as $B$ goes to infinity we obtain infinitely many such pairs.
Now, the height of $a_B/b_B$ is at most $\prod \nrm{b_B}_v$ (since the $v$-adic valuation of $a_B$ is essentially
determined by those of $b_B$ and $x_v$), so we deduce that $H([a_B\colon b_B])\leq B\sum f_v = B$.

Since $e<2$ we may choose $\delta>0$ small enough so that $2e_v/e-\delta > e_v$ for each $v\in S$.  Fix one
such $\delta$.
By the Schmidt Subspace Theorem (see \cite[Corollary 7.2.5]{BG}) applied to the linear forms $a-bx_v$ and $b$
over the places $v$ of $S$, it follows that there is a finite set of lines in $k^2$ which contain all
pairs $(a,b)\in k^2$ satisfying $\nrm{a-x_vb}_{v} < H([a\colon b])^{-f_v}$ and
$\nrm{b}_v < H([a\colon b])^{f_v-\delta}$ for all $v\in S$.

If this finite set of lines contains infinitely many of the $(a,b)\in R^2$ satisfying
$\nrm{a-x_vb}_v< H([a\colon b])^{-f_v}$ and $\nrm{b}_v< H([a\colon b])^{f_v}$ constructed above,
then there is an infinite set of such pairs lying on one of the lines.
That is, there is a fixed $m\in k$ and an infinite set of pairs
$(a,ma)\in R^2$ so that $\nrm{a}_v\cdot\nrm{1-x_vm}_v=\nrm{a-x_vma}_v<H([a\colon ma])^{-f_v}=H([1\colon m])^{-f_v}$.
Since none of the $x_v$ are in $k$, none of the $1-x_vm$ are zero, and this implies
that $\nrm{a}_{v}< C$ for all $v\in S$ and some constant $C$.
Since $\nrm{a}_v\leq 1$ for all $v\not\in S$ this implies that $H([a\colon 1])$ is bounded,  contradicting
the fact that there are infinitely many different $a$.

Thus, there are an infinite
number of pairs $(a,b)\in R^2$ which satisfy $\nrm{a-x_vb}_v < H([a\colon b])^{-f_v}$ and
$\nrm{b}_v \geq H([a\colon b])^{f_v-\delta}$, and hence infinitely many $a/b\in k$ satisfying
$\nrm{a/b-x_v}_v\leq H([a\colon b])^{-2f_v+\delta}$ $=H([a\colon b])^{-2e_v/e+\delta} \leq H([a\colon b])^{-e_v}$.
\epf

Given a sequence $\yseq\subseteq \PP^1(k)$ set $\tau_{v}(\yseq)=1/\alpha_{v}(\yseq,\Osh_{\PP^1}(1))$.
Since $\tau$ is the reciprocal of $\alpha$, if $\tau'< \tau_{v}(\yseq)$ (respectively $\tau'> \tau_{v}(\yseq)$)
then $\lim_{i\to\infty} \sepv(x_v,y_i)^{1/\tau'}H(y_i)=0$ (respectively $=\infty$).
The content
of Theorem \ref{thm:general-dirichlet} is that given any finite set $S$ of places of $k$, and any collection
$\{e_v\}_{v\in S}$ of elements of $[0,2]$ with $\sum e_v < 2$, there is a sequence $\yseq$ such that
$e_v\leq \tau_{v}(\yseq)$ for all $v\in S$.   By a simple diagonal argument we now see that if we choose
the $e_v$ so that $\sum e_v=2$, we may achieve equality.

\begin{corollary}\label{cor:simultaneous-dirichlet}{\bf (Simultaneous weighted Dirichlet):}
Let $S$ be a finite set of places of $k$, and $\{e_v\}_{v\in S}$ a collection of elements of $(0,2]$ such
that $\sum e_v = 2$.  Then there is a sequence $\yseq$ of $k$-points of $\PP^1$ such that
$e_v=\tau_{v}(\yseq)$ for all $v\in S$.
\end{corollary}

\bpf
Let $n_0$ be large enough so that $e_v-\frac{1}{n}>0$ for all $n\geq n_0$ and
all $v\in S$.
By Theorem \ref{thm:general-dirichlet} for each $n\geq n_0$
there is a sequence $\yseqn_{i\geq 0}$ such that $e_v-\frac{1}{2n} \leq \tau_{v}(\yseqn)$.  Since
$e_v-\frac{1}{n} < e_v-\frac{1}{2n}$, we have
$\lim_{i\to\infty} \sepv(x_v,y_{i,n})^{\frac{1}{e_v-1/n}}H(y_{i,n})=0$ for all $v\in S$.
For each fixed $n$, by choosing $i$ large enough, we may pick $y_n=y_{i,n}$ so that
$\sepv(x_v,y_n)^{\frac{1}{e_v-1/n}}H(y_n)< \frac{1}{n}$ and $\sepv(x_v,y_n)< 1$ for all $v\in S$.
In this way we construct a sequence
$\{y_n\}_{n\geq n_0}$ which we simply call $\ynseq$.

\np
Fix $\delta>0$ small enough that $e_{v}-\delta>0$ for each $v\in S$.  For large $n$ we have
$e_{v}-\delta < e_{v}-\frac{1}{n}$ and hence
$$\sepv(x_v,y_{n})^{\frac{1}{e_v-\delta}}H(y_n)< \sepv(x_v,y_{n})^{\frac{1}{e_v-1/n}}H(y_n)< \frac{1}{n}.$$
Therefore
$\lim_{n\to\infty} \sepv(x_v,y_{n})^{\frac{1}{e_v-\delta}}H(y_n) =0$ and so $e_v-\delta\leq \tau_{v}(\ynseq)$.
Letting $\delta$ go to zero we conclude that $e_v\leq \tau_{v}(\ynseq)$ for each $v\in S$.
By Roth's theorem for $\PP^1$ (e.g., Corollary \ref{cor:alpha-hypersurface}) $\sum_v \tau_{v}(\ynseq)\leq 2$.
Since $\sum_{v} e_v=2$ we conclude that $e_v=\tau_{v}(\ynseq)$ for each $v\in S$.
\epf

\np
{\em Proof of Theorem \ref{thm:equality-of-hypersurface}:} \/
In the induction proving Theorem \ref{thm:simul-approx-III}, in order to arrive at equality we must have
gone all the way down to curve $C$, necessarily $k$-rational (since there are infinitely many
rational points, and the approximation constants are finite).  The first result then follows by Roth's theorem
for $\PP^1$ (with the appropriate modification for the singularity, as in Theorem \ref{thm:curve} for a single point).
The converse direction is Corollary \ref{cor:simultaneous-dirichlet} with the choice $e_v=2\xi_v$ for all $v\in S$,
combined with the appropriate modification for the singularity, again as in Theorem \ref{thm:curve}. \epf

As in the case of a single place there are other variations on the deduction of Theorem \ref{thm:simul-approx-III}
from Theorem \ref{thm:simul-approx-II}.

\begin{corollary}\label{cor:simul-var-I}
For any positive integer $m < n$, if we choose $R_v$ so that $R_v>\frac{m}{(m+1)\,\ep_{x_v}(L)}$
for each $v\in S$ then there is a subset $Z$ of $X$ with $\dim Z <m$ such that the equivalent
conditions in Proposition \ref{prop:weighted-equiv} hold with respect to $Z$ and $\{R_v\}_{v\in S}$.
\end{corollary}

\begin{corollary}\label{cor:simul-no-rat-curve}
Suppose that there is no $k$-rational curve passing through any of the $x_v$, $v\in S$.
Then the conditions of Proposition \ref{prop:weighted-equiv} hold with $Z=\emptyset$ and
$R_v = \frac{3}{2\,\ep_{x_v}(L)}$ for all $v\in S$.
\end{corollary}

\noindent
{\em Proof of Corollary \ref{cor:simul-no-rat-curve}: \/}
We prove that condition \eqref{eqn:C} holds with respect to this data, i.e., that given any weighting
function $\xi$ and any collection $\{\delta_v\}_{v\in S}$ of positive real numbers, there are only finitely
many solutions $y\in X(k)$ to

\begin{equation}\label{eqn:simul-no-rat-curve}
\sepv(x_v,y)^{\frac{2\ep_{x_v}(L)}{3}} < H_{L}(y)^{-\xi_v(1+\delta_v)}\,\,\,\mbox{for all $v\in S$}.
\end{equation}

Given the collection $\{\delta_v\}_{v\in S}$ set $\delta'_v=\frac{\delta_v}{2}$ and
$R'_v=\frac{3}{2\,\ep_{x_v}(L)}+\delta'_v$ for all $v\in S$.
By Corollary \ref{cor:simul-var-I} with $m=2$ there is a curve $Z'$, depending on $\{R_v'\}$, so that
there are only finitely many solutions $y\in X(k)\setminus Z'(k)$ to \eqref{eqn:simul-no-rat-curve}.
By hypothesis, there is no $k$-rational curve passing through any of the $x_v$, and so we conclude
that there are only finitely
many solutions $y\in Z'(k)$ to \eqref{eqn:simul-no-rat-curve}.  Thus there are only finitely many
solutions $y\in X(k)$ to \eqref{eqn:simul-no-rat-curve}. \epf

As in Theorem \ref{thm:simul-approx-III} it is probably simplest to express Corollary \ref{cor:simul-no-rat-curve}
in terms of condition \eqref{eqn:B},
i.e., as an inequality governing the position of the point $(\alpha_{v_1}(\xseq,L),\ldots, \alpha_{v_s}(\xseq,L))$
in $\RR^{s}$.  Assuming the hypotheses of the Corollary, for any sequence $\xseq$ of $k$ points,
$$
\sum_{v\in S} \frac{\ep_{x_v}(L)}{\alpha_{v}(\xseq,L)} \leq \frac{3}{2}.
$$

\noindent
{\bf Remark.} \label{rem:equality}
It is clear that it is possible to continue this type of argument if in each dimension $m$ we knew
the types of $m$-dimensional subvariety $Z$ where ``equality'' occurs,   i.e., where there is a sequence
$\xseq$  of points of $Z(k)$, with no subsequence contained in a proper subvariety of $Z$,  satisfying
$$
\sum_{v\in S} \frac{\ep_{x_v,Z}(L)}{\alpha_{v}(\xseq,L)} = \frac{m+1}{m}.
$$
One necessary condition on such a $Z$ is that $Z$ must be Seshadri exceptional (see \S\ref{sec:more-about-Ar})
with respect to each point $x_v$ where $\alpha_{v}(\xseq,L)<\infty$.
(Here Seshadri exceptional means as a subvariety of itself, not as a subvariety of $X$.)
It would already be interesting to work out the case of surfaces.
For instance
$\PP^2$ is such a surface if none of the points $x_v$ lie on $k$-rational lines.

\noindent
{\bf Remark.}
In this section we have used a different constant $R_v$ at each place when describing results on simultaneous
approximation.
By replacing each $R_v$ with the largest (i.e., the worst) of the $R_v$ we obtain a weaker statement, but with the
advantage of the same constant at each place.
Thus, for example, Theorem \ref{thm:simul-approx-III} implies the following product version.

\begin{corollary}\label{cor:simul-product}
Let $\ep=\min_{v\in S}(\ep_{x_v}(L))$.  Then for any $\delta>0$ there are only finitely many solutions $y\in X(k)$ to
$$\prod_{v\in S} \sepv(x_v,y) \leq H_{L}(y)^{-\left(\frac{2}{\ep}+\delta\right)}.$$
\end{corollary}

\section{Improvements via unramified covers}
\label{sec:unram-bounds}

Theorem \ref{thm:simul-approx-I} and an idea due to Robinson-Roquette \cite{RR} and Macintyre \cite{Mac}
(see also \cite[p.\ 100 and \S7.7]{Se}) allow us to give sharper versions of the theorems so far.

In this section by {\em unramified cover} we mean a finite surjective unramified map
$\varphi\colon Y_1\longrightarrow Y_2$ in the category of varieties over $\Spec(k)$, with
both $Y_1$ and $Y_2$ irreducible.

Let $\varphi \colon Y\longrightarrow X$ be an unramified cover and $x$ be any point of $X(\Qbar)$.
As we will see below, for any ample bundle $L$ on $X$, $\min_{y\in \varphi^{-1}(x)}(\Ar_{y}(\varphi^{*}L))$
and $\min_{y\in \varphi^{-1}(x)}(\ep_{y}(\varphi^{*}L))$
are at least as large as $\Ar_{x}(L)$ and $\ep_{x}(L)$ respectively.
We will define $\Arhat_x$ and $\ephat_x$ as suprema over such unramified covers.
The point of this section is that the theorems in \S\ref{sec:Roth-theorems} and \S\ref{sec:simul-approx}
hold with $\Ar$ and $\ep$ replaced by $\Arhat$ and $\ephat$. 
The basic idea is to lift a sequence $\xseq$ on $X$ to a sequence $\yseq$ on $Y$ and use the bounds there;
however the lift involves a change of field, and this introduces a factor which seems to make the result
strictly worse.
Fortunately, by using simultaneous approximation on $Y$ we can exactly cancel out this factor.
In particular, even to get such a result for a single place of $k$ we must use simultaneous approximation on
the cover $Y$.

We first check that $\Ar$ and $\ep$ are weakly increasing in unramified covers;
thus the theorems using $\Arhat$ and $\ephat$ are stronger than the original ones.

\begin{lemma}\label{lem:beta-in-covers}
Let $\varphi \colon Y \longrightarrow X$ be an unramified cover,
$L$ an ample line bundle on $X$, $x$ any point of $X(\Qbar)$, and $y\in \varphi^{-1}(x)$.  Then

\begin{enumerate}
\item $\Ar_{y}(\varphi^{*}L) \geq  \Ar_{x}(L)$, and
\item $\ep_{y}(\varphi^{*}L) \geq \ep_{x}(L)$.
\end{enumerate}
\end{lemma}

\bpf
Let $\Xbar_1$,\ldots, $\Xbar_{r}$ and $\Ybar_1$,\dots, $\Ybar_s$
be the irreducible components of $X\times_{k}\Qbar$
and $Y\times_{k}\Qbar$ containing $x$ and $y$ respectively.   Each $\Ybar_i$ maps to some $\Xbar_j$, and this map
expresses $\Ybar_i$ as an unramified cover of $\Xbar_j$.
Since $\Ar$ and $\ep$ are defined as minima over irreducible components,
establishing the conclusion of the lemma for each map $\Ybar_i\longrightarrow \Xbar_j$ establishes the lemma
for $Y\longrightarrow X$.  Thus we are reduced to the case of studying unramified covers over an algebraically
closed field.   To reduce notation we continue to use $X$ and $Y$ as the names of the varieties, rather than
$\Xbar_j$ and $\Ybar_i$, and $\varphi$ as the name of the map.

Let $\pi_{X}\colon \Xtil\longrightarrow X$ be the blow up of $X$ at $x$ with exceptional divisor $E_x$, and for
any $\gamma\geq 0$ set $L_{X,\gamma}=\pi_{X}^{*}L-\gamma E_x$ and $f_X(\gamma)=\frac{\Vol(L_{X,\gamma})}{\Vol(L)}$.
We similarly let $\pi_{Y}\colon\Ytil\longrightarrow Y$ be the blow up of $Y$ at $y$ with exceptional divisor $E_y$,
and for any $\gamma\geq 0$ we set $L_{Y,\gamma} = \pi_{Y}^{*}\varphi^{*}L - \gamma E_y$ and
$f_Y(\gamma) = \frac{\Vol(L_{Y,\gamma})}{\Vol(\varphi^{*}L)}.$

We first prove (b).
Let $\varphi^{-1}(x) = \{y_1,\ldots, y_\ell\}$ with $y_1=y$.
Since $\varphi$ is unramified, the fibre product
$Y\times_{X}\Xtil$ is the blow up of $Y$ at the points $y_1,\ldots, y_\ell$.
Let $\psi_X$ and $\psi_Y$ be the maps from $Y\times_{X}\Xtil$ to $\Xtil$ and $\Ytil$ respectively (the map
to $\Ytil$ being the blow down at the points of $\varphi^{-1}(x)$ different from $y$).  Thus we have the following
commutative diagram of maps
$$
\begin{array}{c}
\begin{pspicture}(-2,-2)(2,2.5)
\rput(-2,2){$Y\times_X\times\Xtil$}
\rput(-2,0){$\Ytil$}
\rput(-2,-2){$Y$}
\rput(2,-2){$X$}
\rput(2,2){$\Xtil$}
\psset{arrows=->}
\psline(-2,1.6)(-2,0.3)
\rput(-1.6,0.95){\small$\psi_{Y}$}
\psline(-2,-0.3)(-2,-1.7)
\rput(-1.6,-1){\small$\pi_{Y}$}
\psline(-1.7,-2)(1.7,-2)
\rput(0,-1.7){\small $\varphi$}
\psline(2,1.6)(2,-1.7)
\rput(2.4,0){{\small $\pi_{X}$}.}
\psline(-1.0,2)(1.7,2)
\rput(0,2.3){\small $\psi_{X}$}
\end{pspicture}
\end{array}
$$

For $i=2,\ldots, \ell$ let $E_i$ be the exceptional divisor of $\psi_{Y}$ lying over $y_i$.
The description of the fibre product as a further blowup of $\Ytil$ shows that for any
$\gamma$ we have the equality of line bundles

\begin{equation}\label{eqn:fibre-prod-pullback}
\psi_{X}^{*} L_{X,\gamma} = \psi_{Y}^{*} L_{Y,\gamma} - \gamma\left(\sum_{i=2}^{d} E_i\right).
\end{equation}

If $0\leq\gamma\leq \ep_{x}(L)$ then
$L_{X,\gamma}$ is nef on $\Xtil$ and so $\psi_{X}^{*}L_{X,\gamma}$ is nef on $Y\times_{X}\Xtil$.
Equation \eqref{eqn:fibre-prod-pullback} then implies that
$L_{Y,\gamma}$ is nef on $\Ytil$.   This proves (b).

We will prove (a) by showing the inequality $f_{Y}(\gamma)\geq f_{X}(\gamma)$ for all $\gamma\geq 0$.
Since both $f_X$ and $f_Y$ are continuous functions, it suffices to prove the inequality
for rational $\gamma$.

Set $\Esh=\varphi_{*}\Osh_{Y}$ and let $d$ be the generic rank of $\Esh$. By the projection formula,
for any $m>0$ we have
$\varphi_{*}(\varphi^{*}mL) = mL \otimes_{\Osh_{x}} \Esh$, and so $H^0(Y,\varphi^{*}mL) = H^0(X,mL\otimes \Esh)$.
The volume measures the leading term in the asymptotic growth of global sections, and for this purpose
tensoring with the (generic) rank $d$ sheaf $\Esh$ has the same effect as tensoring with $d$ copies of $\Osh_{X}$.
Therefore $\Vol(\varphi^{*}L) = d \Vol(L)$.
Similarly, for any $\gamma\geq 0$ we have
$\Vol(\psi_{X}^{*}L_{X,\gamma}) = d \Vol(L_{X,\gamma})$.

For any rational $\gamma\geq 0$, and any $m\geq 0$ sufficiently divisible so that $m\gamma$ is integral,
multiplying \eqref{eqn:fibre-prod-pullback} by $m$ shows that global sections of
$m\psi_{X}^{*}L_{X,\gamma}$ are a subspace of the global sections of $m\psi_{Y}^{*}L_{Y,\gamma}$
and so $\Vol(\psi_{Y}^{*}L_{Y,\gamma})\geq \Vol(\psi_{X}^{*}L_{X,\gamma})$.
Finally, since $\psi_{Y}$ is birational, $\Vol(L_{Y,\gamma})=\Vol(\psi_{Y}^{*}L_{Y,\gamma})$.
We thus have
$$f_{Y}(\gamma) = \frac{\Vol(\psi_{Y}^{*}L_{Y,\gamma})}{\Vol(\varphi^{*}L)}
\geq \frac{\Vol(\psi_{X}^{*}L_{X,\gamma})}{\Vol(\varphi^{*}L)} = \frac{d\Vol(L_{X,\gamma})}{d\Vol(L)}=f_X(\gamma),$$
and integrating gives $\Ar_{y}(\varphi^{*}L)\geq \Ar_{x}(L)$.
\epf

In the category of schemes over $X$,
consider the full subcategory whose objects are the unramified covers $\varphi\colon Y\longrightarrow X$ as above.
If $(Y_1,\varphi_1)$ and $(Y_2,\varphi_2)$ are objects
and $\psi\colon Y_1\longrightarrow Y_2$ a morphism in this category, then $\psi$ expresses $Y_1$ as an unramified
cover of $Y_2$, and thus Lemma \ref{lem:beta-in-covers} applies. In particular, for any $y_2\in Y_2(\Qbar)$,
$\min_{y_1\in \psi^{-1}(y_2)}(\Ar_{y_1}(\varphi_1^{*}L)) \geq \Ar_{y_2}(\varphi_2^{*}L)$
and similarly for $\ep$.

\begin{definition}\label{def:Arhat}
Let $X$ be an irreducible variety over $\Spec(k)$, $L$ an ample line bundle on $X$ and $x\in X(\Qbar)$.
We define
$$
\begin{array}{rclcrcl}
\Arhat_{x}(L) &=&
\displaystyle \sup_{\varphi\colon Y\longrightarrow X} \min_{y\in \varphi^{-1}(x)}\Ar_{y}(\varphi^{*}L)
&\mbox{and} &
\ephat_{x}(L) &=&
\displaystyle \sup_{\varphi\colon Y\longrightarrow X } \min_{y\in \varphi^{-1}(x)}(\ep_{y}(\varphi^{*}L)),\\
\end{array}
$$
where the suprema are over the set of unramified covers $\varphi\colon Y \longrightarrow X$.
\end{definition}

In the arguments below it will be important to know we can find a single unramified cover which approximates
finitely many of the $\Arhat_{x}(L)$.

\begin{lemma}\label{lem:single-unram-cover}
Let $X$ be an irreducible variety over $\Spec(k)$, $L$ an ample line bundle on $X$, and $x_1$,\ldots, $x_\ell$
finitely many points of $X(\Qbar)$.  Suppose that $\Ar_1$,\ldots, $\Ar_{\ell}$ are positive real numbers
with $\Ar_i < \Arhat_{x_i}(L)$ for $i=1$,\ldots, $\ell$.  Then there exists an unramified cover
$\varphi\colon Y \longrightarrow X$ such that $\min_{y\in \varphi^{-1}(x_i)}(\Ar_{y}(\varphi^{*}(L)) > \Ar_i$
for $i=1$,\ldots, $\ell$.
\end{lemma}

\bpf
By the definition of $\Arhat$, for each $i$ there is an unramified cover $\varphi_i\colon Y_i\longrightarrow X$
such that $\min_{y\in \varphi_i^{-1}(x_i)}(\Ar_{y}(\varphi_i^{*}(L)) >  \Ar_i$.  Let $Y$ be any irreducible component
of the fibre product $Y_1\times_{X} \cdots \times_{X} Y_{\ell}$, and $\varphi\colon Y\longrightarrow X$ the
induced map.  The natural projection maps of the fibre product induce maps $\psi_i\colon Y\longrightarrow Y_i$
for each $i$, and $\psi_i$ expresses $Y$ as an unramified cover of $Y_i$.   For any $y\in \varphi^{-1}(x_i)$,
$\psi_i(y)\in \varphi_i^{-1}(x_i)$, and hence an application of
Lemma \ref{lem:beta-in-covers} to the unramified cover $\psi_i$ shows that $Y$ has the desired property.
\epf

\medskip\noindent
\hypertarget{rem:Galois-symmetries-a}{{\bf Remarks on Galois symmetries.}} 
(a) Suppose that $x\in X(\Qbar)$. Given any algebraic conjugate $x'$ of $x$
let $\sigma\in \Gal(\Qbar/k)$ be an element such that $\sigma(x)=x'$.
Given any unramified cover $\varphi\colon Y\longrightarrow X$, the action of $\sigma$ on $Y(\Qbar)$ then
takes points of $Y(\Qbar)$ lying over $x$ to points lying over $x'$.
We conclude that for any ample line bundle $L$ on $X$, $\min_{y\in \varphi^{-1}(x)} \Ar_{y}(\varphi^{*}L)=
\min_{y'\in\varphi^{-1}(x')}\Ar_{y'}(\varphi^{*}L)$. 
This Galois symmetry argument also shows that
$\ephat_{x}(L)=\ephat_{x'}(L)$ and $\Arhat_{x}(L)=\Arhat_{x'}(L)$ for any ample line bundle $L$ on $X$.

\noindent
\hypertarget{rem:Galois-symmetries-b}{(b)}
Let $v_0$ be a place of $k$, $v$ and $v'$ two extensions of $v_0$ to $\Qbar$, and
$\sigma\in \Gal(\Qbar/k)$ such that $v'=v\circ\sigma$. 
Fix a point $x\in X(\Qbar)$ and set $x'=\sigma(x)$.
Suppose that a sequence $\xseq\subseteq X(k)$ converges to $x$ with respect to a distance function $\sepv$.
If we define $\sepvp$ by
using the same embedding $X\hookrightarrow \PP^{r}_{k}$ used to define $\sepv$, then applying $\sigma$ to formula
\eqref{eqn:distance-formula-arch} or \eqref{eqn:distance-formula-nonarch} shows that
$\sepvp(x',x_i)=\sepv(x,x_i)$ for all $i\geq 0$.  More generally, if $\sepvp$ is defined by using a different
embedding of $X$ then this result combined with Proposition \ref{prop:equiv-dist}
shows that $\sepvp(x',x_i)$ and $\sepv(x,x_i)$ are equivalent as $i\to \infty$.
Summarizing, 
if $\xseq$
converges to $x\in X(\Qbar)$ with respect to $\sepv$, and if $v'$ is a different extension of $v_0$ to $\Qbar$,
then $\xseq$ will converge to an algebraic conjugate of $x$ with respect to $\sepvp$,
with the same essential speed of convergence.

\medskip\noindent
{\bf Lifting sequences.}  Let $\psi\colon Y'\longrightarrow X$ be an unramified cover.
By the theorem of Chevalley-Weil \cite[Theorem 10.3.11]{BG} there is a finite extension $F/k$ such that
all points $\{y\in Y'(\Qbar) \st \psi(y) \in X(k)\}$ are defined over $F$ (this field $F$ is not unique, since
any larger field will also work).
Fix such a field $F$.  It will be convenient for us that the covering variety is also a variety over $\Spec(F)$.
To do this we let $Y$ be an irreducible component of $Y'\times_{k} F$, and
$\varphi\colon Y\longrightarrow X$ the induced map.
Via the natural map $\Spec(F)\longrightarrow\Spec(k)$, $Y$ is a variety over $\Spec(k)$, and
$\varphi$ is an unramified cover.  Furthermore, all points of $Y(\Qbar)$ lying over points of $X(k)$ are again
defined over $F$. 

Given a sequence $\xseq$ of points of $X(k)$, for each $i$  we arbitrarily choose $y_i\in Y(F)$ with $\varphi(y_i)=x_i$.
We call such a sequence $\yseq$ a {\em lift} of $\xseq$.
This lift is somewhat haphazard, but by further passing to a subsequence we may obtain a lift with better properties.

Let $v_0$ be a place of $k$, extended to a place $v$ on $\Qbar$,
and suppose that there is $x_v\in X(\Qbar)$ such that $\sepv(x_v,x_i)\to 0$, i.e., that $\xseq$ approximates $x_v$
with respect to $\sepv(\cdot,\cdot)$.   Let $w_0$ be a place of $F$ lying over $v_0$, and $w$ an extension
of $w_0$ to $\Fbar$.  In this situation we define $x_w$, an algebraic conjugate of $x$, as follows.
On $\Qbar=\Fbar$, $w$ gives a place $v'$ of $\Qbar$ extending $v_0$, but which may not be equal to $v$.
We then apply \hyperlink{rem:Galois-symmetries-b}{(b)} of the 
\hyperlink{rem:Galois-symmetries-a}{`Remarks on Galois symmetries'} above to obtain
an algebraic conjugate $x_w$ of $x$.  With respect to $\sepvp$, $\xseq$ converges to $x_w$.

Returning to the problem of improving the lift,
since $Y(F_w)$ is compact by passing to a subsequence we may assume
that the sequence $\yseq$ has a limit $y_w\in Y(F_w)$.
Since $Y$ is a variety over $\Spec(F)$, the place $w$ gives a distance function $\sepw(\cdot,\cdot)_{F}$ on $Y$.
The topology on $Y(F_w)$ is that induced by $\sepw(\cdot,\cdot)_{F}$, and so this means that
$\sepw(y_w,y_i)_{F}\to 0$ as $i\to \infty$.
Furthermore, by continuity we have $\varphi(y_w)=x_w$, in particular, $y_w\in Y(\Qbar)$.

We will need a generalization obtained by repeating this procedure.
Let $T_{v}$ be the set of places of $F$ over $v_0$, each extended to a place of $\Fbar$.
As above, for each such $w\in T_v$ we obtain an algebraic conjugate $x_w$ of $x_v$.  (These conjugates are
not necessarily distinct.)
By applying the procedure above to
each $w\in T_{v}$ in turn, we may find $\Qbar$-points $y_w\in \varphi^{-1}(x_w)$ for each $w\in T_{v}$, and
a subsequence of $\yseq$ so that for each $w\in T_{v}$, $\sepw(y_w,y_i)_{F}\to 0$ as $i\to \infty$.

Finally, given a finite set $S$ of places of $k$ extended to $\Qbar$, we may repeat this process for each $v\in S$.
We record the conclusion below.

\begin{proposition}\label{prop:lift-approx}
Let $\psi\colon Y'\longrightarrow X$ be an unramified cover, $S$ a finite set of
places of $k$ each extended to $\Qbar$, and $F/k$ a finite extension so that for all $x\in X(k)$, 
all points of $\psi^{-1}(x)$ are defined over $F$.
We replace $Y'$ by a component $Y$ of $Y'\times_{k}F$, and let $\varphi\colon Y\longrightarrow X$ be the
induced map.
Suppose that $\{x_v\}_{v\in S}$ are a set of points of $X(\Qbar)$, and that
$\xseq$ is a sequence of $k$-points so that $\sepv(x_v,x_i)\to 0$ for each $v\in S$.
For each $v\in S$ let $T_{v}$ be the set of places of $F$
lying over $v_0=v|_{k}$ each extended to a place of $\Fbar$.
For each such $w$ we let $x_w$ be the corresponding algebraic conjugate of $x_v$ as defined above. 

Then by passing to a subsequence of $\xseq$ we may find a lift
$\yseq$ of $\xseq$ to $Y$,
and for each $v\in S$ and $w\in T_{v}$ a $\Qbar$-point $y_w\in \varphi^{-1}(x_w)$, such that
$\sepw(y_w,y_i)_{F}\to 0$ as $i\to\infty$.
\end{proposition}

We next compare the asymptotics of $\sepv(x_v,x_i)_{k}$ with $\sepw(y_w,y_i)_{F}$, and the resulting effect on
$\alpha$.

\begin{lemma}\label{lem:etale-distance}
Let $\psi\colon Y'\longrightarrow X$ be an unramified cover, $v_0$ a place of $k$, $v$ an extension of $v_0$
to $\Qbar$, and $\xseq$ a sequence of points of $X(k)$ converging to $x\in X(\Qbar)$
with respect to $\sepv(\cdot,\cdot)$.
Let $F/k$ be a finite extension so that for all $x\in X(k)$, all points of $\psi^{-1}(x)$ are defined over $F$.
Let $Y$ be a component of $Y'\times_{k} F$ and $\varphi\colon Y\longrightarrow X$ the induced map.
Let $\yseq$ be a lift of $\xseq$ to $Y(\Qbar)$, $w_0$ a place of $F$ lying over $v_0$,
and $w$ an extension of $w_0$ to $\overline{F}$.
Suppose that $\yseq$ converges to $y\in Y(\Qbar)$ with respect to $\sepw(\cdot,\cdot)_{F}$.
Finally, set $m_w=[F_w\colon k_v]$ {\rm (}$=[F_{w_0}\colon k_{v_0}]${\rm )} and $e=[F\colon k]$.

Then $\sepw(y,y_i)_{F}$ is asymptotically equivalent to $\sepv(x,x_i)^{m_w}_{k}$ as $i\to \infty$, and for any
line bundle $L$ on $X$, $\alpha_{y}(\yseq,\varphi^{*}L)_{F} =  \frac{e}{m_w}\alpha_{x}(\xseq,L)_{k}$.
\end{lemma}

Here, as in Proposition \ref{prop:change-of-field},
the subscript $k$ or $F$ on $\sepv(\cdot,\cdot)_k$ or $\sepw(\cdot,\cdot)_F$ respectively indicates the field 
being used to normalize the distance function.  The subscript on $\alpha(\cdot,\cdot)$
similarly indicates the field used to normalize the distance function and the height.

\bpf
For each $y_i$ in the sequence, we have $H_{\varphi^{*}L}(y_i)_{F} = H_{L}(x_i)_{F}= H_{L}(x_i)_{k}^{e}$.
Thus once we show that
$\sepw(y,y_i)_{F}$ is asymptotically equivalent to $\sepv(x,x_i)_{k}^{m_w}$ as $i\to \infty$ the equality
$\alpha_{y}(\yseq,\varphi^{*}L)_{F} =  \frac{e}{m_w}\alpha_{x}(\xseq,L)_{k}$ follows immediately as in the proof of
Proposition \ref{prop:change-of-field}.

The geometric point $x\in X(\Qbar)$ defines a morphism $\Spec(\Qbar)\longrightarrow X$, whose image
is a closed point $[x]\in X$.  Let $U'$ be an affine neighbourhood of $[x]$,
and let $u_1$, \ldots, $u_r$ be elements of
$\Gamma(U',\Osh_{X})$ which generate the maximal ideal of $[x]$.
Let $K$ be the Galois closure (over $k$) of the field of definition of $y$, and
set $x'=\varphi(y)$; as before $x'$ is an algebraic conjugate of $x$.
Since $y$ maps to $x'$, and since $x'$ and $x$ have isomorphic fields of definition,
$K$ also contains the field of definition of $x$.
Pulling back to $U_{K}':=U'\times_{k}K$,
the functions $u_1$,\ldots, $u_r$ cut out all the geometric points lying over $[x]$.
By passing to a smaller affine open $U\subset U_{K}'$ we may assume that $u_1$,\ldots, $u_r$
generate the maximal ideal of $x$.
By Lemma \ref{lem:local-distance}, $\sepv(x,x_i)_k$ is equivalent to
$\max(\nrm{u_1(x_i)}_{v},\ldots, \nrm{u_r(x_i)}_{v})$ as $i\to\infty$, and so
$\sepv(x,x_i)^{m_w}_{k}$ is equivalent to
$\max(\nrm{u_1(x_i)}_{v}^{m_w},\ldots, \nrm{u_r(x_i)}_{v}^{m_w})$ as $i\to\infty$.

Let $V'_K=\varphi^{-1}_{K}(U'_K)$.
Since $\varphi$ is unramified, $\varphi^{*}u_1$,\ldots, $\varphi^{*}u_r$ cut out the reduced points lying
over all algebraic conjugates of $x$.  Therefore
if we restrict to a small enough affine open neighbourhood $V$ of $y$ in $V_K'$
then $\varphi^{*}u_1$,\ldots, $\varphi^{*}u_r$ generate the maximal ideal of $y$.
Since the $\yseq$ converge to $y$ with respect to $\sepw(\cdot,\cdot)_F$, we may apply
Lemma \ref{lem:local-distance} again to conclude that $\sepw(y,y_i)_F$ is equivalent to
$\max(\nrm{(\varphi^{*}u_1)(y_i)}_{w},\ldots, \nrm{(\varphi^{*}u_r)(y_i)}_{w})$ as $i\to \infty$.
Since $(\varphi^{*}u_j)(y_i) = u_j(\varphi(y_i))=u_j(x_i)$ for each $j=1,\ldots, r$ and each $i$,
we will be done if we show that $\nrm{u_j(x_i)}_{w}=\nrm{u_j(x_i)}_{v}^{m_w}$ for all $i$, $j$.

The value $u_j(x_i)$ is the value of the residue class of $u_j$ in $\kappa(x_i)$. Since the $x_i$ are defined
over $k$, $\kappa(x_i)=k$ and so $u_j(x_i)\in k$ for all $i$, $j$.
Since $w_0$ lies over $v_0$ (and $w$ is normalized with respect to $F$),
the restriction of $\nrm{\cdot}_{w}$ to $k\subset F$ is $\nrm{\cdot }_{v}^{m_w}$. This proves the lemma.
\epf

Applying Lemma \ref{lem:etale-distance} to the lift produced in Proposition \ref{prop:lift-approx} yields
the following corollary.

\begin{corollary}\label{cor:lift-asymptotics}
Assume the setup and notation of Proposition \ref{prop:lift-approx}, and let $\yseq$  and
$\{y_w\}_{w\in T_{v}, v\in S}$ be the lift and set of points provided by its conclusion.
Then for any $\QQ$-bundle $L$ on $X$, every $v\in S$, and $w\in T_{v}$ we have
$\alpha_{y_w}(\yseq,\varphi^{*}L)_{F} = \frac{[F\colon k]}{[F_w\colon k_v]}\, \alpha_{x_v}(\xseq,L)_{k}$.

Here $\alpha_{y_w}$ is computed with respect to $\sepw(\cdot,\cdot)_{F}$ and $\alpha_{x_v}$ with respect to
$\sepv(\cdot,\cdot)_{k}$.
\end{corollary}

\noindent
We are now ready to establish the version of Theorem \ref{thm:simul-approx-I} with $\Arhat$ in place of $\Ar$.

\begin{theorem}\label{thm:unram-approx-I}
Let $X$ be an irreducible variety defined over $k$, $S$ a finite set of places of $k$, each extended to $\Qbar$,
and for each $v\in S$ choose an $x_v\in X(\Qbar)$.   Suppose that $L$ is an ample $\QQ$-bundle on $X$ defined over $k$,
and that $\{R_v\}_{v\in S}$ are a collection of positive real numbers such that

\begin{equation}\label{eqn:R-condition-asymp}
\sum_{v\in S} \Arhat_{x_v}(L)R_v > 1.
\end{equation}

\noindent
Then \eqref{eqn:alpha-condition} and \eqref{eqn:finite-condition} hold with respect to the collection
$\{R_v\}_{v\in S}$.
\end{theorem}

\bpf
By condition \eqref{eqn:R-condition-asymp} and Lemma \ref{lem:single-unram-cover} we may find
an unramified cover $\psi\colon Y'\longrightarrow X$ satisfying

\begin{equation}\label{eqn:R-condition-initial-cover}
\sum_{v\in S} \left(\min_{y_v\in \psi^{-1}(x_v)}(\Ar_{y_v}(\psi^{*}L))\right)R_v > 1.
\end{equation}

\noindent
Let $F/k$ be a finite extension so that
any sequence $\xseq$ of $k$-points of $X$ can be lifted to a sequence $\yseq$ of points of $Y'$ defined over
$F$, and fix $F$ for the rest of the proof.  Let $Y$ be an irreducible component of
$Y'\times_{k}F$, and $\varphi\colon Y\longrightarrow X$ the induced map.
Since $Y$ is a component of $Y'\times_{k}F$, by definition of $\Ar$ (Definition \ref{def:Ar})
we have $\Ar_{y,Y}(\varphi^{*}L)\geq \Ar_{y,Y'}{\psi^{*}L}$ for every point $y\in Y(\Qbar)$, and
so by \eqref{eqn:R-condition-initial-cover}

\begin{equation}\label{eqn:R-condition-cover}
\sum_{v\in S} \left(\min_{y_v\in \varphi^{-1}(x_v)}(\Ar_{y_v}(\varphi^{*}L))\right)R_v > 1.
\end{equation}

We will prove the theorem in the form of condition \eqref{eqn:alpha-condition}, that is, we will show that
there is a proper subvariety $Z \subset X$, such that for all sequences $\xseq$ of $k$-points
of $X\setminus Z$ there is at least one $v\in S$ so that $\alpha_{v}(\xseq,L)\geq \frac{1}{R_v}$.
Here as in \S\ref{sec:central-thm} for each $v\in S$ we use $\alpha_{v}$ to mean $\alpha_{x_v}$,
computed with respect to the distance $\sepv(\cdot,\cdot)$.

For each $v\in S$ we set $T_{v}$ to be the set of places of $F$ lying over $v_0:=v|_{k}$, each extended
to a place of $\Fbar$.
Each $v\in S$ and $w\in T_{v}$ determine an algebraic conjugate $x_w$ of $x_v$ as above.

Let $\xseq$ be a sequence of $k$-points of $X$.
If there is a $v\in S$ so that $\sepv(x_v,x_i)$ does not go to $0$ as $i\to \infty$, then $\alpha_{v}(\xseq,L)=\infty$,
and the statement to be proved is trivially satisfied.   We may therefore restrict ourselves to studying
sequences $\xseq$ so that $\sepv(x_v,x_i)\to 0$ as $i\to \infty$ for each $v\in S$,
and we do so for the rest of the proof.
We note again that passing to a subsequence can only possibly lower the value of $\alpha$, so we may
freely do so in proving the result.

Set $T=\bigsqcup_{v\in S} T_{v}$.
By Proposition \ref{prop:lift-approx}
if $\xseq$ converges to each $x_v$ with respect to $\sepv(\cdot,\cdot)$, then for any lift $\yseq$ of $\xseq$  
we may, after passing to a subsequence of $\yseq$, choose a $\Qbar$-point $y_w\in \varphi^{-1}(x_w)$ 
for each $w\in T$ so that $\yseq$ converges to $y_w$ with respect to $\sepw(\cdot,\cdot)_F$.

Thus, up to passing to a subsequence, for any sequence $\xseq$ of $k$-points of $X$ there is
a choice of $y_w\in \varphi^{-1}(x_w)$ for each $w\in T$ such that

\begin{equation}\label{eqn:lipstick}
\rule{1.5cm}{0cm}
\left\{
\begin{tabular}{c}
\mbox{\begin{minipage}{0.85\textwidth}
the sequence $\{x_i\}$ lifts to a sequence $\yseq$ in $Y$ (necessarily defined over $F$) which converges with
respect to $\sepw(\cdot,\cdot)_F$ to $y_w$ for each $w\in T$.
\end{minipage}} \\
\end{tabular}\right\}\phantom{.}
\end{equation}

We will show that for each of the finitely many elements $q$ of the product set
$\prod_{w\in T} \varphi^{-1}(x_w)$, i.e., each of the finitely many choices of a $y_w\in \varphi^{-1}(x_w)$
for each $w\in T$, there exists
a proper subvariety $Z_q\subset X$ (depending on these choices) so that for any
sequence $\xseq$ of $k$-points of $X\setminus Z_q$ satisfying \eqref{eqn:lipstick}
there is at least one $v\in S$ so that $\alpha_{v}(\xseq,L)\geq \frac{1}{R_v}$.
Taking $Z$ to be the union over the finitely many such $Z_q$ then yields the theorem.

We now assume that we have fixed $q=\{y_w\}_{w\in T} \in \prod_{w\in T} \varphi^{-1}(x_w)$
and prove the existence of such a $Z_q$.
For each $w\in T$ set $g_w=\frac{[F\colon k]}{[F_w\colon k_v]}$, where $v\in S$ is such that $w\in T_v$.
If $\xseq$ is a sequence of $k$-points satisfying \eqref{eqn:lipstick} above, and $\yseq$ such a lift,
Corollary \ref{cor:lift-asymptotics} gives
$\alpha_{w}(\yseq,\varphi^{*}L)_{F} =g_w\,\alpha_{v}(\xseq,L)_{k}$ for each $v\in S$ and $w\in T_{v}$.
Here, as before, we use
$\alpha_{w}$ to mean $\alpha_{y_w}$ computed with respect to $\sepw(\cdot,\cdot)_F$, and the subscripts $F$
and $k$ to indicate the field used to normalize the distance and the heights.

For each $w\in T$ set $R_{w}'=R_{v}/g_w$, where again $v\in S$ is such that $w\in T_v$.
Combining: (1)
for each $v\in S$ we have $\sum_{w\in T_v} \frac{1}{g_w}=1$;
(2) for each $v\in S$, $w\in T_v$,
$\min_{y\in \varphi^{-1}(x_v)}\Ar_{y}(\varphi^{*}L)=\min_{y'\in \varphi^{-1}(x_w)}\Ar_{y'}(\varphi^{*}L)$
(see \hyperlink{rem:Galois-symmetries-a}{`Remarks on Galois symmetries'} 
\hyperlink{rem:Galois-symmetries-a}{(a)});
and (3) inequality \eqref{eqn:R-condition-cover},
we conclude that
$$\sum_{w\in T} \Ar_{y_w}(\varphi^{*}L) R_w' > 1.$$

Working over $\Spec(F)$,
we now apply Theorem \ref{thm:simul-approx-I} to the collection $\{R_w'\}_{w\in T}$ and line bundle
$\varphi^{*}(L)$, and let $Z_q'$ be the resulting proper subvariety of $Y$.
Set $Z_q$ to be the image of $Z_q'$ in $X$.
Since $Y$ is irreducible, $Z_q'$ is of dimension strictly less than $Y$, and hence $Z_q$ is again a proper
subvariety of $X$.

Suppose that $\xseq$ is a sequence of $k$-points of $X\setminus Z_q$ satisfying \eqref{eqn:lipstick},
and let $\yseq$ be such a lift.
Then $\yseq$ is contained in the $F$-points of $Y\setminus Z_q'$ and thus by construction of $Z_q'$
there is at least one $v\in S$ and $w\in T_v$
so that $$\alpha_{w}(\yseq,\varphi^{*}L)_{F}\geq \frac{1}{R_w'} = \frac{g_w}{R_v}.$$
Since $\alpha_{w}(\yseq, \varphi^{*}L)_{F}=g_w\,\alpha_{v}(\xseq,L)_{k}$ we conclude that
$\alpha_{v}(\xseq,L)_{k} \geq \frac{1}{R_v}$.  Thus $Z_q$ has the required property, and
this completes the proof of Theorem \ref{thm:unram-approx-I}. \epf

\medskip
\noindent
{\bf Remarks.} (a) Besides the fact that $\beta$ is weakly increasing in unramified covers, the
keys to the proof are (1) the fact that $\alpha_{v}(\xseq,L)_k$ and $\alpha_{w}(\yseq,\varphi^{*}L)_F$
differ by a factor of $1/g_w = \frac{[F_w:k_v]}{[F:k]}$, whenever $w_0 := w|_{F}$ lies over
$v_0:=v|_{k}$, and (2) for each $v\in S$, $\sum_{w\in T_v} 
\frac{1}{g_w}=1$, which allows us to get rid of this factor by using simultaneous approximation.

\np
(b) If one of the $\Arhat_{x_v}(L)$ is infinite then condition \eqref{eqn:R-condition-asymp} holds for
any collection $\{R_v\}_{v\in S}$ of positive numbers.

\begin{proposition} \label{prop:hat-version}
Let $X$ be an irreducible variety of dimension $n$, $L$ an ample line bundle on $X$, and $x\in X(\Qbar)$.
Then

\begin{enumerate}
\item $\Arhat_{x}(L) \geq \frac{n}{n+1}\ephat_{x}(L)$, and
\item for any irreducible subvariety $Z\subseteq X$ and any $x\in Z(\Qbar)$, $\ephat_{x,Z}(L|_{Z})\geq \ephat_{x,X}(L)$.
\end{enumerate}
\end{proposition}

\bpf
Let $\varphi\colon Y\longrightarrow X$ be an unramified cover.
By Corollary \ref{cor:beta-bound} we have $\Ar_{y}(\varphi^{*}(L))\geq \frac{n}{n+1}\ep_{y}(\varphi^{*}L)$ for
each $y\in \varphi^{-1}(x)$.  Thus $
\min_{y\in \varphi^{-1}(x)}(\Ar_{y}(\varphi^{*}L)) \geq \frac{n}{n+1}
\min_{y\in \varphi^{-1}(x)}(\ep_{y}(\varphi^{*}L))$, and (a) follows after taking the supremum over such covers.

For part (b), let $Z'$ be any irreducible component of $\varphi^{-1}(Z)$, where $\varphi\colon Y \longrightarrow X$
is an unramified cover as above.
The induced map $\psi\colon Z'\longrightarrow Z$ expresses $Z'$ as an unramified cover over $Z$,
and for any $z\in \psi^{-1}(x)$ we have
$\ep_{z,Z'}(\psi^{*}(L|_{Z})) =\ep_{z,Z'}((\varphi^{*}L)|_{Z'}) \geq \ep_{z,Y}(\varphi^{*}L)$
by Proposition \ref{prop:ex}(c).  Since $z\in \varphi^{-1}(x)$, this implies
$$\min_{z\in\psi^{-1}(x)}(\ep_{z,Z'}(\psi^{*}(L|_{Z}))) \geq \min_{z\in \psi^{-1}(x)}(\ep_{z,Y}(\varphi^{*}L))
\geq \min_{y\in \varphi^{-1}(x)}(\ep_{y,Y}(\varphi^{*}L)).$$
Taking the suprema over unramified covers of $Z$ and $X$ we deduce (b).
\epf

Once Theorem \ref{thm:simul-approx-I} is established, the approximation results in
\S\ref{sec:central-thm}---\S\ref{sec:simul-approx}
follow from that theorem, Corollary \ref{cor:beta-bound}, Propositions \ref{prop:ex}(c) and \ref{prop:alpha}(f),
as well as arguments common in Diophantine approximation.
The necessary results about $\alpha$, $\Ar$, $\ep$, and their asymptotic versions $\Arhat$ and $\ephat$ 
needed to make these arguments are summarized in the following table.

\noindent
\hspace{-0.3cm}
\begin{tabular}{|lc|cl|}
\hline
\bf Results about $\alpha$, $\Ar$ and $\ep$
& & &
\bf Results about $\alpha$, $\Arhat $ and $\ephat$ \\
\hline
\underline{Theorem \ref{thm:simul-approx-I}} & & & \underline{Theorem \ref{thm:unram-approx-I}} \\
\begin{minipage}{0.42\textwidth}
$\sum_{v} \Ar_{x_v}(L)R_v > 1 \implies \eqref{eqn:alpha-condition} + \eqref{eqn:finite-condition}$
\end{minipage}
& & &
\begin{minipage}{0.42\textwidth}
$\sum_{v} \Arhat_{x_v}(L)R_v > 1 \implies \eqref{eqn:alpha-condition} + \eqref{eqn:finite-condition}$
\end{minipage}  \\
 & & & \\
\underline{Corollary \ref{cor:beta-bound}} & & & \underline{Proposition \ref{prop:hat-version}(a)} \\
\begin{minipage}{0.4\textwidth}
$\Ar_{x}(L) \geq \frac{n}{n+1}\ep_{x}(L)$
\end{minipage}
& & &
\begin{minipage}{0.4\textwidth}
$\Arhat_{x}(L) \geq \frac{n}{n+1}\ephat_{x}(L)$
\end{minipage} \\
 & & & \\
\underline{Proposition \ref{prop:ex}(c)} & & & \underline{Proposition \ref{prop:hat-version}(b)} \\
\begin{minipage}{0.4\textwidth}
$\ep_{x,Z}(L|_{Z})\geq \ep_{x,X}(L)$
\end{minipage}
& & &
\begin{minipage}{0.4\textwidth}
$\ephat_{x,Z}(L|_{Z})\geq \ephat_{x,X}(L)$
\end{minipage} \\
 & & & \\
\underline{Proposition \ref{prop:alpha}(f)} & & & \underline{Proposition \ref{prop:alpha}(f)} \\
\begin{minipage}{0.47\textwidth}
$\alpha_{x,X}(L) = \min(\alpha_{x,X_1}(L|_{X_1}),\ldots, \alpha_{x,X_r}(L|_{X_r}))$
\end{minipage}
& & &
\begin{minipage}{0.48\textwidth}
$\alpha_{x,X}(L) = \min(\alpha_{x,X_1}(L|_{X_1}),\ldots, \alpha_{x,X_r}(L|_{X_r}))$
\end{minipage} \\
\hline
\end{tabular}

\medskip
By using Theorem \ref{thm:unram-approx-I} in place of Theorem \ref{thm:simul-approx-I} and Proposition
\ref{prop:hat-version}(a,b) in place of Corollary \ref{cor:beta-bound} and Proposition \ref{prop:ex}(c)
respectively, the arguments in \S\ref{sec:central-thm}---\S\ref{sec:simul-approx}
hold with $\Arhat$ and $\ephat$ used in place of $\Ar$ and $\ep$.  Explicitly, we have the following synthesis
of the arguments in \S\ref{sec:central-thm}---\S\ref{sec:unram-bounds}.

\begin{corollary}\label{cor:everything}
Theorems \ref{thm:simul-approx-I}, \ref{thm:RothI}, \ref{thm:RothII}, \ref{thm:RothIII}, \ref{thm:simul-approx-II},
\ref{thm:simul-approx-III}, and \ref{thm:equality-of-hypersurface}, Corollaries \ref{cor:switch-quantifiers},
\ref{cor:equality-in-dim-one}, \ref{cor:one-place-Roth}, \ref{cor:control-dim-Z}, \ref{cor:no-rat-curve},
\ref{cor:alpha-hypersurface}, \ref{cor:simul-var-I}, \ref{cor:simul-no-rat-curve}, and
\ref{cor:simul-product} hold with $\Ar$ and $\ep$ replaced by $\Arhat$ and $\ephat$.
\end{corollary}

\medskip
\noindent
{\bf Remark.} The larger the values of $\Arhat$ and $\ephat$  the stronger these types of results are.
In particular, this means that given any lower bounds for $\Arhat$ and $\ephat$  the
results listed in Corollary \ref{cor:everything} hold with the lower bounds used in place of
$\Ar$ or $\ep$.  One method of getting lower bounds for $\Arhat$ and $\ephat$ which still takes into account
the asymptotic behaviour of covers is to consider only \etale Galois covers $\varphi\colon Y \longrightarrow X$
with $Y$ irreducible.
By the transitivity of the Galois action, for any ample line bundle $L$ on $X$, both $\ep_{y}(\varphi^{*}L)$
and $\Ar_{y}(\varphi^{*}L)$ are independent of $y\in \varphi^{-1}(x)$, and thus we avoid worrying which point
in the fibre achieves the minimum.

In particular, setting

$$
\begin{array}{rclcrcl}
\Arhatet_{x}(L) &=&
\displaystyle \sup_{\substack{\varphi\colon Y\longrightarrow X \\ y\in\varphi^{-1}(x)}} \Ar_{y}(\varphi^{*}L)
&\mbox{and} &
\ephatet_{x}(L) &=&
\displaystyle \sup_{\substack{\varphi\colon Y\longrightarrow X \\ y\in\varphi^{-1}(x)}} \ep_{y}(\varphi^{*}L).\\
\end{array}
$$
where the suprema are over irreducible \etale Galois covers $\varphi\colon Y \longrightarrow X$, we
obtain lower bounds $\Arhat_{x}(L) \geq \Arhatet_{x}(L)$ and $\ephat_{x}(L)\geq \ephatet_{x}(L)$
for all $x\in X(\Qbar)$ and ample $L$.

\medskip
\noindent
{\bf Example.}
\label{ex:ab-var-etale}
Let $X$ be an abelian variety and let
$[m]\colon X\longrightarrow X$ denote the multiplication by $m$ map.
For any ample line bundle $L$, $[m]^{*}L$ has the same numerical class as $m^2 L$,
and so $\ep_{x}([m]^{*}L)=m^2\ep_{x}(L)$ for any $x\in X(\Qbar)$.
In particular, $\ephatet_{x}(L)=\infty$ and thus $\ephat_{x}(L)=\infty$.
Therefore for any $x\in X(\Qbar)$, $\alpha_{x}(L)\geq \frac{1}{2}\ephat_{x}(L) = \infty$
by the unramified cover version of Theorem \ref{thm:RothIII}.
(This gives another proof of example (c) on page \pageref{ex:alpha-intro} of the introduction.)

\medskip
\noindent
{\bf Remark.}  If $X$ is normal then any unramified cover of $X$ is \'etale, and any such cover
can be dominated by a Galois \etale cover.
Thus if $X$ is normal $\Arhatet$ and $\ephatet$ agree with $\Arhat$ and $\ephat$.

One of the themes of this article is the comparison of $\alpha$ and $\ep$.  In light of Corollary \ref{cor:everything}
it is natural to ask if $\ephat$ has same formal properties shared by $\alpha$ and $\ep$ (i.e. perhaps we have
been writing the wrong article).
We have not defined $\ephat$ when $X$ is reducible, and we do it now by simply adopting one of the desired properties
of $\ephat$ as the definition.  If $X$ is reducible over $k$, $x\in X(\Qbar)$ and $X_1$,\ldots, $X_r$ the irreducible
components passing through $x$ then we set
$\ephat_{x,X}(L) = \min(\ephat_{x,X_1}(L|_{X_1}),\ldots, \ephat_{x,X_r}(L|_{X_r}))$.

\begin{proposition}\label{prop:exhat}
Let $X$ be a projective variety defined over $k$, $x\in X(\Qbar)$, and $L$ a nef $\QQ$-divisor on $X$.
Consider the following assertions:

\begin{enumerate}
\item For any positive integer $m$, $\ephat_{x}(m\cdot L)=m\cdot \ephat_{x}(L)$.
\item  $\ephat_{x}$ is a concave function of $L$:
for any positive rational numbers $a$ and $b$, and any nef $\QQ$-divisors $L_1$ and $L_2$
$$\ephat_x(a L_1+ b L_2) \geq a\ephat_x(L_1)+b\ephat_x(L_2).$$
\item If $Z$ is a subvariety of $X$ then for any point $z\in Z(\Qbar)$
we have $\ephat_{z,Z}(L|_Z)\geq \ephat_{z,X}(L)$.
\item If $L$ is very ample then $\ephat_{x}(L)\geq 1$, if $L$ is ample then
  $\ephat_{x,X}(L)> 0$.
\item If $x$ and $y$ are points of varieties $X$ and $Y$, with nef
line bundles $L_X$ and $L_Y$ then
\[\ephat_{x\times y,X\times Y}(L_X\squareplus L_Y)
= \min(\ephat_{x,X}(L_X), \ephat_{y,Y}(L_Y)).\]
\item Suppose that $X$ is reducible and let $X_1$,\ldots, $X_r$ be the irreducible components containing $x$.
Then $\ephat_{x,X}(L) = \min(\ephat_{x,X_1}(L|_{X_1}),\ldots, \ephat_{x,X_r}(L|_{X_r}))$.
\end{enumerate}

Then (a), (c), (d), and (f) hold.  We do not know if (b) and (e) hold in general, but they do hold
when $X$ is normal (respectively $X$ and $Y$ are normal).
\end{proposition}

\bpf
Part (f) holds by definition of $\ephat$.   It follows from the definition that establishing any of (a)---(e) for
irreducible $X$ implies the corresponding result for reducible $X$, so from now on we assume that $X$ (or $Y$)
is irreducible over $k$.
Then parts (a) and (d) follow immediately from Proposition \ref{prop:ex}(a,d) and the definition of $\ephat$,
and part (c) is Proposition \ref{prop:hat-version}(b).

The difficulty with (b) is that the definition of $\ephat$ involves the minimum over covers, and it is not clear
that the minimum of all three of $\ephat(aL_1+bL_2)$, $\ephat(L_1)$ , and $\ephat(L_2)$ happen at the same point
and can be compared.  However for \etale Galois covers, since we do not have to worry about the minimum, we can
compare at any point and then it is clear the inequality holds by Proposition \ref{prop:ex}(b).
Thus, in particular, (b) holds when $X$ is normal.

Similarly, if $X$ and $Y$ are normal, so that again we may just consider \etale Galois covers, (e) follows
from Proposition \ref{prop:ex}(e) and the fact that such any such cover is a product of an \etale Galois
cover of $X$ with an \etale Galois cover of $Y$.  (Specifically, let
$\Xbar_1$,\ldots, $\Xbar_r$ and $\Ybar_1$, \ldots, $\Ybar_s$ be the irreducible components of
$X\times_{k}\Qbar$ and $Y\times_{k}\Qbar$ respectively.
Note that all $\Xbar_i$ and $\Ybar_j$ are isomorphic over $\Qbar$, and that $r=1$ and $s=1$ if
$X$ and $Y$ are geometrically connected.
For any \etale Galois cover $\varphi\colon V\longrightarrow X\times Y$,
after passing to the algebraic closure, which we do when computing
$\ep$, each connected component of $V\times_{k}\Qbar$ is an \etale Galois cover of some $\Xbar_i\times \Ybar_j$,
and hence is a product of an \etale Galois covers of $\Xbar_i$ and $\Ybar_j$.    These \etale Galois covers of
$\Xbar_i$ and $\Ybar_j$ may be descended to Galois covers of $X$ and $Y$ respectively.)
\epf

\comment{Check this last statement}

\noindent
{\bf Remark.}
From the arguments for (b) and (e) above, it may seem that $\ephatet$ is a better substitute for $\ephat$, since
for $\ephatet$ properties (b) and (e) hold for any variety, even non-normal ones.  However, if $X$ is not normal,
it is not clear that property (c) holds for $\ephatet$.  In the argument of Proposition \ref{prop:hat-version}(b)
it was necessary to pass to a component of a cover of $Z$, and a component of an \etale cover is not necessarily
\'etale.  This is one of the reasons for the definition of $\ephat$ as a supremum over unramified covers.

\section{\texorpdfstring{More about $\Ar_x(L)$}{More about \pdfbeta(L)}}
\label{sec:more-about-Ar}

In this section we discuss interpretations of and further results and remarks about $\Ar_{x}(L)$.
For simplicity we assume that $X$ is irreducible and defined over an algebraically closed field.

\noindent
{\bf Heuristic Interpretation of $\mathbf{\Ar}$.}
Let $L$ be an ample $\QQ$-bundle on $X$ and $x\in X$. As in \S\ref{sec:f-def}
we define a function $f(\gamma)=\Vol(L_{\gamma})/\Vol(L)$ for $\gamma\geq 0$, and set $\geff=\geffx(L)$.
The function $f$ is decreasing with $f(0)=1$ and $f(\geff)=0$ (Figure \ref{fig:P1xP1-graph} is a good illustration).
By \cite[Corollary C]{LM} or \cite[Theorem A]{BFJ} the volume function is first-differentiable and hence so is $f$.
The function $1-f$ therefore satisfies the criteria to be a cumulative distribution function.

It is straightforward to say what the associated probability distribution is measuring.
Suppose for the sake of discussion that $L$ is an integral line bundle and base point free.
For a fixed $\gamma> 0$, what is the probability that
a randomly chosen section of $V=\Gamma(X,L)$ vanishes to order $\geq \gamma$ at $x$?   Since the set of sections
vanishing to order $\geq \gamma$ at $x$ forms a proper subspace $W_\gamma$ of $V$, under the usual
probability measure the chance is zero.  However if we instead decide the ratio
$\dim W_{\gamma}/\dim V$ is a good measure of the chance that a section of $V$ lies in $W_{\gamma}$,
and further decide that we should really ask the question asymptotically,
that is, assign the limit $\dim W_{m\gamma}/\dim \Gamma(X,mL)$ as $m\to \infty$
as the probability of the event, then we arrive exactly at $f(\gamma)$.  Therefore (under this strange distribution)
$1-f(\gamma)$ is the probability that a section vanishes to order $\leq \gamma$, and $-f'(\gamma)$ the probability
density function for vanishing to order exactly $\gamma$.

The first computation one usually does when given a probability measure is to compute the expected value.
Since $-f'$ is supported on $[0,\geff]$, and since $f(\geff)=0$, integration by parts gives
$$\mathbb{E}(\gamma) = -\int_{0}^{\geff} \gamma f'(\gamma)\,d\gamma  =
\left.{-\gamma f(\gamma)\rule{0cm}{0.6cm}}\right|_{\gamma=0}^{\gamma=\geff} +
\int_{0}^{\geff} f(\gamma)\,d\gamma = -  0 + 0 + \Ar_{x}(L) = \Ar_{x}(L).
$$
This gives an interpretation of $\Ar_{x}(L)$: under the probability distribution above
$\Ar_{x}(L)$ is the expected order of vanishing at $x$ of a section of $L$.

The idea that the probability an element of a vector space $V$ lies in a subspace $W$ should be
$\dim W/\dim V$ is counter to our intuition under the uniform measure,
however it is exactly this type of probability measure which is used by Faltings-W\"ustholz in the proof of their
approximation theorem (see \cite[\S4]{FW}).
Thus, with the exception of the passage to the limiting distribution, which is simply to get
better control over the behaviour of the line bundle, $-f'$ is the probability measure used in the proof of the
Faltings-W\"{u}stholz approximation theorem.
It is therefore completely natural that the expected order of vanishing at $x$
governs approximation results as in Theorem \ref{thm:simul-approx-I}.

\medskip
\noindent
{\bf Other results.}
In Corollary \ref{cor:beta-bound} we showed the inequalities
$\Ar_{x}(L) \geq \frac{n}{n+1} \sqrt[n]{\frac{\Vol(L)}{\mult_{x}X}}\geq \frac{n}{n+1}\ep_{x}(L)$, and
we have used this to deduce approximation theorems involving $\ep$ from those involving $\Ar$.
If the inequalities are strict then replacing $\Ar$ by $\frac{n}{n+1}\ep$ produces a weaker result.
It is therefore natural to ask when these inequalities are equalities.

\begin{theorem}\label{thm:inequal-is-equal}
Let $X$ be an $n$-dimensional irreducible variety, $x\in X$ and $L$ an ample $\QQ$-bundle on $X$.
Then the following conditions are equivalent.

\begin{enumerate}
\item $\Ar_{x}(L) =  \frac{n}{n+1}\sqrt[n]{\frac{\Vol(L)}{\mult_{x}X}}$

\medskip
\item $\sqrt[n]{\frac{\Vol(L)}{\mult_{x}X}} = \ep_{x}(L)$

\medskip
\item $\Ar_{x}(L) = \frac{n}{n+1}\ep_{x}(L)$

\medskip
\item $\ep_{x}(L) = \geffx(L)$.

\end{enumerate}

\end{theorem}

\bpf
To simplify the notation somewhat, set
$\Ar_{x}=\Ar_{x}(L)$,
$\omega_{x} = \sqrt[n]{\frac{\Vol(L)}{\mult_{x}X}}$,
$\ep_{x}=\ep_{x}(L)$,
and $\geff=\geffx(L)$.

\noindent
(a) $\implies$ (b): The estimate $\Ar_{x}\geq \frac{n}{n+1}\omega_{x}$ resulted from integrating the lower bound
$\Vol(L_{\gamma})/\Vol(L)\geq 1 - \frac{\mult_{x}(X)}{\Vol(L)}\gamma^n$ over $[0,\omega_{x}]$.
The equality in (a) is therefore equivalent to the two statements:

\noindent
\hspace{-0.6cm}
\begin{tabular}{rl}
\bf \forceandshowlabel{eqn:9.1.a.1}{(9.1.a.1)} & $\Vol(L_{\gamma})=\Vol(L) - (\mult_{x}X)\gamma^n$ for $\gamma\in [0,\omega_{x}]$,
and \\
\bf \forceandshowlabel{eqn:9.1.a.2}{(9.1.a.2)} &  $\geff=\omega_{x}$. \rule{0cm}{0.6cm}\\
\end{tabular}

Here (as usual) $L_{\gamma} = \pi^{*}L -\gamma E$ and $\pi\colon \Xtil\longrightarrow X$ is the blow up of $X$ at $x$
with exceptional divisor $E$.
We will see that \ref{eqn:9.1.a.1} implies (b).  We first recall an extension of the idea of volume
to arbitrary cohomology groups.  For any line bundle $M$ on an $n$-dimensional variety  $Y$,  and any $0\leq i\leq n$
we set
$$\hhat^{i}(M) = \lim_{m\to\infty} \frac{\dim H^{i}(Y,mM)}{m^n/n!}$$
so that $\hhat^{0}(M) = \Vol(M)$.  As in the case of the volume, the groups $\hhat^{i}$ depend only on the
numerical class of $M$, make sense for $\QQ$-divisors, and for fixed $i$ extend to continuous functions on
$\NS(Y)_{\RR}$ (see \cite[p.\ 1477]{K}).  We will also need a slight variation of this idea.
As in \S\ref{sec:f-def} for any rational $\gamma>0$
and $m$ such that $m\gamma$ is an integer we denote by $m\gamma E$ the subscheme defined by the $(m\gamma)$-th power
of the defining equation for $E$.  For any $0\leq i\leq n$ we set

$$\hhat^{i}(\Osh_{\gamma E}) = \lim_{m\to \infty} \frac{\dim H^{i}(\Xtil, \Osh_{m\gamma E})}{m^n/n!}$$
where the limit runs over all $m$ such that $m\gamma$ is an integer.  Note that ``$\hhat^{i}(\Osh_{\gamma E})$''
is being defined as an atomic symbol --- we are not giving any meaning to $\Osh_{\gamma E}$ as a scheme.
Since $\Osh_{E}(-E)$ is ample on ${E}$, it follows from Serre vanishing and \eqref{eqn:O-gamma-filtration}
that $\hhat^{i}(\Osh_{\gamma E})=0$ for all $i>0$.  Combined with this, the
argument in the proof of Lemma \ref{lem:beta-bound} actually shows that
$\hhat^{0}(\Osh_{\gamma E}) = (\mult_{x}X)\gamma^{n}$.

The asymptotic cohomology groups are birational invariants.  Since $L$ is ample, $h^{i}(X,mL)=0$ for all $m\gg 0$,
and hence (pulling back to $\Xtil$) $\hhat^{i}(L_0)=0$ for all $i>0$.  The long exact sequence associated to
\eqref{eqn:exact-gammaE-sequence} then implies that for any rational $\gamma \geq 0$,
$\hhat^{i}(L_{\gamma})=0$ for all $i\geq 2$ and that
$$\Vol(L_{\gamma}) - \hhat^{1}(L_{\gamma}) = \Vol(L) - \hhat^{0}(\Osh_{\gamma E}) = \Vol(L) - (\mult_{x} X )\gamma^n.$$
Thus \ref{eqn:9.1.a.1} is equivalent to the statement that
$\hhat^{1}(L_{\gamma})=0$ for all $0\leq \gamma \leq \omega_{x}$.

Let $A$ be any ample bundle on $\Xtil$.  By \cite[Theorem A]{dFKL} $L_{\gamma}$ is ample if and only if
$\hhat^{i}(L_{\gamma}-tA)=0$ for all $i>0$ and all sufficiently small $t$.  Let $s$ be any number $0< s < \ep$
so that $A=L_{s}$ is ample on $\Xtil$. Then
$L_{\gamma}-tA = (1-t)\pi^{*}L - (\gamma-ts)E  = (1-t)L_{\frac{\gamma-ts}{1-t}}$. The asymptotic cohomology groups
are homogeneous of degree $n$, so
$$\hhat^{i}(L_{\gamma} - t L_{s}) = \hhat^{i}\left({(1-t)L_{\frac{\gamma-ts}{1-t}}}\right) =
(1-t)^{n}\hhat^{i}\left( L_{\frac{\gamma-ts}{1-t}}\right)$$
for all $i\geq 0$.  If $0<\gamma < \omega_{x}$, then for small enough $t$ we have
$0\leq\frac{\gamma-ts}{1-t}<\omega_{x}$ too,
and hence by \ref{eqn:9.1.a.1} and the equation above $\hhat^{i}(L_{\gamma}-tA)=0$ for all $i>0$.

Summarizing, condition \ref{eqn:9.1.a.1} and Theorem A of \cite{dFKL} imply that
$L_{\gamma}$ is ample for all $0<\gamma < \omega_{x}$.  Thus $\omega_{x}\leq \ep_{x}$.
The opposite inequality, $\ep_{x}\leq \omega_{x}$, is \cite[Proposition 5.1.9]{PAG}
(this already appeared in the proof of Corollary \ref{cor:beta-bound-irred}) and thus $\ep_{x}=\omega_{x}$,
i.e., (b) holds.

\noindent
(b) $\implies$ (c)+(d): Since $\Vol(L_{\gamma}) = \Vol(L) - (\mult_{x}X)\gamma^n$ for $\gamma\in [0,\ep_{x}]$, and since
condition (b) is that $\omega_{x}=\ep_{x}$, we have

\noindent
\hspace{-0.60cm}
\begin{tabular}{rl}
\bf \forceandshowlabel{eqn:9.1.b.1}{(9.1.b.1)}
&  $\Vol(L_{\gamma})=\Vol(L) - (\mult_{x}X)\gamma^n$ for
$\gamma\in [0,\omega_{x}]$.
\end{tabular}

\noindent
Condition \ref{eqn:9.1.b.1} shows that $\Vol(L_{\gamma})>0$ for $0\leq \gamma < \omega_{x}$, and
that $\Vol(L_{\omega_{x}})=0$, hence $\omega_{x}$ is the boundary of the effective cone, i.e.,

\vspace{-0.3cm}
\noindent
\hspace{-0.60cm}
\begin{tabular}{rl}
\bf \forceandshowlabel{eqn:9.1.b.2}{(9.1.b.2)} &  $\geff=\omega_{x}$. \rule{0cm}{0.6cm} \\
\end{tabular}

\noindent
Given these two conditions,

$$\Ar_{x} = \int_{0}^{\geff} \Vol(L_{\gamma})/\Vol(L)\, d\gamma =
\int_{0}^{\omega_x} 1 - \textstyle \frac{\mult_{x}X}{\Vol(L)} \gamma^n\,d\gamma = \frac{n}{n+1}\omega_{x} = \frac{n}{n+1}\ep_{x}.$$
Thus (c) holds.
Since (d) is condition \ref{eqn:9.1.b.2} it is also clear that (b) implies (d).

\noindent
(c) $\implies$ (a)+(b):
This is clear from the inequalities $\Ar_{x} \geq \frac{n}{n+1}\omega_{x}\geq \frac{n}{n+1}\ep_{x}$.

\noindent
(d) $\implies$ (b): This is immediate from the inequalities $\geff \geq \omega_{x} \geq \ep_{x}$.  \epf

\medskip
\noindent
{\bf Remark.}
Condition (b) of Theorem \ref{thm:inequal-is-equal} seems the easiest one to check in practice.
Condition (d) is also tractable;
it is the
statement that along the ray $\pi^{*}L-\gamma E$ ($\gamma\geq 0$), the point where the ray exits the nef cone
is the same point where the ray exits the effective cone.

\bigskip
\noindent
{\bf Seshadri Exceptional Subvarieties.}
Recall that by \cite[Proposition 5.1.9]{PAG} for any irreducible subvariety $V\subseteq X$ of positive dimension
passing through $x$ we have the
inequality

\refstepcounter{eqcount}
\begin{equation}\label{eqn:seshadri-except}
\ep_{x}(L) \leq \left({\frac{c_1(L)^{\dim V}\cdot V}{\mult_{x} V}}\right)^{\frac{1}{\dim V}},
\end{equation}

\noindent
and that there are irreducible subvarieties $V$ for which \eqref{eqn:seshadri-except} is an equality
(including possibly $X=V$).
An irreducible subvariety $V$ is called {\em Seshadri exceptional} (with respect to $x$ and $L$)
if \eqref{eqn:seshadri-except} is an equality,
and if $V$ is not properly contained in a larger subvariety having the same property.
Condition (b) of Theorem \ref{thm:inequal-is-equal} is that $X$ itself is Seshadri exceptional.

\medskip
\noindent
{\bf Further properties of $\Ar_{x}(L)$.}
As in previous sections, it is interesting to work out some formal properties of $\Ar_{x}$,
in particular to ask whether the list of properties in Propositions \ref{prop:alpha} and \ref{prop:ex} hold.
We do not know the status of all the properties listed there, and simply record some elementary observations.
(The letters match those of Propositions \ref{prop:alpha} and \ref{prop:ex}.)

\begin{proposition}
$x\in X$, $L$ an ample line bundle on $X$, then

\begin{enumerate}
\item[(a)] $\Ar_{x}(mL) = m\Ar_{x}(L)$.
\item[(c)] If $Z$ is a subvariety of $X$, $x\in Z$, it is {\em not} necessarily true
that $\Ar_{x,Z}(L|_{Z})\geq\Ar_{x,X}(L)$.
\item[(d)] If $L$ is ample then $\Ar_{x}(L)> 0$.
\item[(f)] Suppose that $X$ is reducible and let $X_1$,\ldots, $X_r$ be the irreducible components containing $x$.
Then $\Ar_{x,X}(L) = \min(\Ar_{x,X_1}(L|_{X_1}),\ldots, \Ar_{x,X_r}(L|_{X_r}))$.
\end{enumerate}
\end{proposition}

\bpf
Property (f) holds by definition of $\Ar_{x}$ (Definition \ref{def:Ar}), and (d) is clear from the
estimate $\Ar_{x}(L)\geq \frac{n}{n+1}\ep_{x}(L)$ and Proposition \ref{prop:ex}(d).
For part (a), fix $m>0$ and
let $f_{L}(\gamma)$ and $f_{mL}(\gamma)$ be the functions $f_{L}(\gamma) = \Vol(L_{\gamma})/\Vol(L)$ and
$f_{mL}(\gamma) = \Vol((mL)_{\gamma})/\Vol(mL)$ respectively.
On an $n$-dimensional variety one has $\Vol(mM) = m^n\Vol(M)$ for every big line bundle $M$ and $m>0$ and hence
$$
f_{mL}(m\gamma) = \Vol((mL)_{m\gamma})/\Vol(mL) = \Vol(mL_{\gamma})/\Vol(mL) =
\frac{m^n}{m^n} \Vol(L_{\gamma})/\Vol(L) = f_{L}(\gamma).
$$

It follows from this equation  or directly from the definition  that $\geffx(mL) = m\geffx(L)$.
Integrating (and using the previous equation) we conclude that $\Ar_{x}(mL)=m\Ar_{x}(L)$.

Finally, to see that $\Ar_{x}$ may strictly decrease under restriction,
recall that $\Ar_x(\Osh_{\PP^n(1)})=\frac{n}{n+1}$ for any point $x\in \PP^n$
(see the example on page \pageref{ex:Ar-Pn}).  Hence  if $Z$ is an $m$-dimensional linear subspace
of $X=\PP^n$ passing through $x$ (with $m<n$) and $L=\Osh_{\PP^n}(1)$ then $\Ar_{x,Z}(L|_Z)  < \Ar_{x}(L)$.  \epf

\noindent
{\bf Remark.}
The fact that $\ep_{x}$ is weakly increasing under restriction has been crucial for our inductive arguments.
The fact $\Ar_{x}$ may decrease under restriction to a subvariety is one reason why this article is focussed
on $\ep_{x}$, and why it was important to estimate $\Ar_{x}$ in terms of $\ep_{x}$.

\section{A special case of Vojta's main conjecture}
\label{sec:Vojta}

Vojta's Main Conjecture (Conjecture 3.4.3 of \cite{Vo}) predicts how the height of rational points grow
as they approach a simple normal crossings divisor $D\subset X$.
One can also investigate the prediction for other subvarieties of $X$, with the result being stronger for
larger subvarieties.
Theorem \ref{thm:simul-approx-II}
easily implies many special cases of the Main Conjecture, albeit ones where the subvariety
in question is a collection of points
(this is natural since the results of this paper are geared towards approximating points).
Despite the fact that this is a weaker version than the classical case in which $D$ is a divisor, many of the cases
established below were previously unknown.

We refer the reader to \cite[\S3]{Vo} for a statement and discussion of the Main Conjecture,
and simply state the relevant result in the language of this paper.

\begin{theorem}\label{thm:vojta}
Let $k$ be a number field, $X$ an irreducible $n$-dimensional variety over $\mbox{Spec}(k)$ such that $-K_X$ is ample,
$D$ a finite subset of $X(\Qbar)$, and $S$ a finite set of places of $k$.  
If $\ep_{x}(-K_{X})> \frac{n+1}{n}$ for every $x\in D$,
then Vojta's Main Conjecture is true for $X$ and $D$.
Specifically, for every $\delta>0$ and any big divisor $A$, there is a closed subset $Z\subset X$ such that
for all $k$-rational points $P\in X(k)\setminus Z(k)$, we have:
\[\sum_{v\in S,\, x\in D} -\log\sepv(x,P)+h_{K_X}(P) < \delta h_A(P) +O(1)\]
\end{theorem}

\noindent
{\em Proof:} \/ If we can show the inequality for one big divisor $A$, then it will immediately follow for an arbitrary
big divisor $A$, by adjusting $\delta$ and $Z$.  Thus, we may assume that $A=-K_X$.  Furthermore, note that for any 
place $v\in S$, there is at most one point in $D$ for which $-\log\sepv(x,P)$ contributes
more than a bounded amount to the sum.  Therefore, we may apply Theorem \ref{thm:simul-approx-II}(b) (in the equivalent form of \eqref{eqn:D})
with $R_v=1$ for each $v\in S$ to see that
there is a proper subset $Z$ so that for any $\delta>0$ the equation
$$\prod_{v\in S,\, x\in D} \sepv(x,P) > H_{-K_X}(P)^{-(1+\delta)}$$
holds for all but finitely many $P\in X(k)\setminus Z(k)$.
Taking $\log$ then gives the result. \epf

\noindent
{\bf Remark.} For $k$-points $x$ of $D$ it is sufficient that the weaker condition $\ep_{x}(-K_X)\geq 1$ hold.
One uses the Liouville bound $\alpha_{x}(-K_X)\geq\ep_{x}(-K_X)$ (valid for points of $X(k)$ -- see \cite{McKR})
in a simultaneous approximation version similar to Corollary \ref{cor:alpha-hypersurface}.

\vspace{.1in}

There are many examples of varieties $X$ which satisfy the criterion of the theorem.
For example for any variety of the form $X=G/P$ where $G$ is a semi-simple algebraic group and $P$ is a parabolic
subgroup (e.g., $\PP^n$ or Grassmannians) one has $\ep_{x}(-K_X)\geq 2$ for all points $x\in X(\Qbar)$.

\end{document}